% processed by citealice (August 1999) on Tue Oct 12 10:35:30 CDT 1999

% To: Saharon Shelah <shelah@math.huji.ac.il>,        Andrzej Roslanowski <roslanowski@unomaha.edu>
% Subject: latest revisions
% Date: Mon, 11 Oct 1999 13:07:31 -0500
% From: Alice Leonhardt <leonhard@math.rutgers.edu>
% Mime-Version: 1.0
% X-sliced-and-diced-by: 'savemail' 0.3, Feb 1999

\ifx\shlhetal\undefinedcontrolsequence\let\shlhetal\relax\fi

\input amstex
% % \input mathdefs
\expandafter\ifx\csname mathdefs.tex\endcsname\relax
  \expandafter\gdef\csname mathdefs.tex\endcsname{}
\else \message{Hey!  Apparently you were trying to
  \string\input{mathdefs.tex} twice.   This does not make sense.} 
\errmessage{Please edit your file (probably \jobname.tex) and remove
any duplicate ``\string\input'' lines}\endinput\fi

%mathdefs.tex v1.3.2

%%% Changes from v1.0: footnote macros, warning for duplicated tags,
%%%   control sequences \( and \verbatimtags.
%%% From v1.2: \pretags, redefinition of \( using \ifinner, multi-part
%%%   equation numbering, control sequences \[, \references, and
%%%   \resetbracket. 
%%% From v1.3: \rm in \lastpart; write root of multi-part tag to .tgs 

%See file texdefs.doc for documentation.

\catcode`\X=12\catcode`\@=11

%Minor control sequences:
\def\n@wcount{\alloc@0\count\countdef\insc@unt}
\def\n@wwrite{\alloc@7\write\chardef\sixt@@n}
\def\n@wread{\alloc@6\read\chardef\sixt@@n}
\def\r@s@t{\relax}\def\v@idline{\par}\def\@mputate#1/{#1}
\def\l@c@l#1X{\firstpart.#1}\def\gl@b@l#1X{#1}\def\t@d@l#1X{{}}

%Creation of tag families and output of assignments and citations:
\def\crossrefs#1{\ifx\all#1\let\tr@ce=\all\else\def\tr@ce{#1,}\fi
   \n@wwrite\cit@tionsout\openout\cit@tionsout=\jobname.cit 
   \write\cit@tionsout{\tr@ce}\expandafter\setfl@gs\tr@ce,}
\def\setfl@gs#1,{\def\@{#1}\ifx\@\empty\let\next=\relax
   \else\let\next=\setfl@gs\expandafter\xdef
   \csname#1tr@cetrue\endcsname{}\fi\next}
\def\m@ketag#1#2{\expandafter\n@wcount\csname#2tagno\endcsname
     \csname#2tagno\endcsname=0\let\tail=\all\xdef\all{\tail#2,}
   \ifx#1\l@c@l\let\tail=\r@s@t\xdef\r@s@t{\csname#2tagno\endcsname=0\tail}\fi
   \expandafter\gdef\csname#2cite\endcsname##1{\expandafter
     \ifx\csname#2tag##1\endcsname\relax?\else\csname#2tag##1\endcsname\fi
     \expandafter\ifx\csname#2tr@cetrue\endcsname\relax\else
     \write\cit@tionsout{#2tag ##1 cited on page \folio.}\fi}
   \expandafter\gdef\csname#2page\endcsname##1{\expandafter
     \ifx\csname#2page##1\endcsname\relax?\else\csname#2page##1\endcsname\fi
     \expandafter\ifx\csname#2tr@cetrue\endcsname\relax\else
     \write\cit@tionsout{#2tag ##1 cited on page \folio.}\fi}
   \expandafter\gdef\csname#2tag\endcsname##1{\expandafter
      \ifx\csname#2check##1\endcsname\relax
      \expandafter\xdef\csname#2check##1\endcsname{}%
      \else\immediate\write16{Warning: #2tag ##1 used more than once.}\fi
      \multit@g{#1}{#2}##1/X%
      \write\t@gsout{#2tag ##1 assigned number \csname#2tag##1\endcsname\space
      on page \number\count0.}%
   \csname#2tag##1\endcsname}}
\def\multit@g#1#2#3/#4X{\def\t@mp{#4}\ifx\t@mp\empty%
      \global\advance\csname#2tagno\endcsname by 1 
      \expandafter\xdef\csname#2tag#3\endcsname
      {#1\number\csname#2tagno\endcsnameX}%
   \else\expandafter\ifx\csname#2last#3\endcsname\relax
      \expandafter\n@wcount\csname#2last#3\endcsname
      \global\advance\csname#2tagno\endcsname by 1 
      \expandafter\xdef\csname#2tag#3\endcsname
      {#1\number\csname#2tagno\endcsnameX}
      \write\t@gsout{#2tag #3 assigned number \csname#2tag#3\endcsname\space
      on page \number\count0.}\fi
   \global\advance\csname#2last#3\endcsname by 1
   \def\t@mp{\expandafter\xdef\csname#2tag#3/}%
   \expandafter\t@mp\@mputate#4\endcsname
   {\csname#2tag#3\endcsname\lastpart{\csname#2last#3\endcsname}}\fi}
\def\t@gs#1{\def\all{}\m@ketag#1e\m@ketag#1s\m@ketag\t@d@l p
   \m@ketag\gl@b@l r \n@wread\t@gsin
   \openin\t@gsin=\jobname.tgs \re@der \closein\t@gsin
   \n@wwrite\t@gsout\openout\t@gsout=\jobname.tgs }
\outer\def\localtags{\t@gs\l@c@l}
\outer\def\globaltags{\t@gs\gl@b@l}
\outer\def\newlocaltag#1{\m@ketag\l@c@l{#1}}
\outer\def\newglobaltag#1{\m@ketag\gl@b@l{#1}}

%Reading in tag information:
\newif\ifpr@ 
\def\m@kecs #1tag #2 assigned number #3 on page #4.%
   {\expandafter\gdef\csname#1tag#2\endcsname{#3}
   \expandafter\gdef\csname#1page#2\endcsname{#4}
   \ifpr@\expandafter\xdef\csname#1check#2\endcsname{}\fi}
\def\re@der{\ifeof\t@gsin\let\next=\relax\else
   \read\t@gsin to\t@gline\ifx\t@gline\v@idline\else
   \expandafter\m@kecs \t@gline\fi\let \next=\re@der\fi\next}
\def\pretags#1{\pr@true\pret@gs#1,,}
\def\pret@gs#1,{\def\@{#1}\ifx\@\empty\let\n@xtfile=\relax
   \else\let\n@xtfile=\pret@gs \openin\t@gsin=#1.tgs \message{#1} \re@der 
   \closein\t@gsin\fi \n@xtfile}

%Sections and subsections; local numbering:
\newcount\sectno\sectno=0\newcount\subsectno\subsectno=0
\newif\ifultr@local \def\ultralocal{\ultr@localtrue}
\def\firstpart{\number\sectno}
\def\lastpart#1{\ifcase#1 \or a\or b\or c\or d\or e\or f\or g\or h\or 
   i\or k\or l\or m\or n\or o\or p\or q\or r\or s\or t\or u\or v\or w\or 
   x\or y\or z \fi}

\def\resetall{\global\advance\sectno by 1\subsectno=0
   \gdef\firstpart{\number\sectno}\r@s@t}
\def\resetsub{\global\advance\subsectno by 1
   \gdef\firstpart{\number\sectno.\number\subsectno}\r@s@t}
\def\newsection#1\par{\resetall\vskip0pt plus.3\vsize\penalty-250
   \vskip0pt plus-.3\vsize\bigskip\bigskip
   \message{#1}\leftline{\bf#1}\nobreak\bigskip}
\def\subsection#1\par{\ifultr@local\resetsub\fi
   \vskip0pt plus.2\vsize\penalty-250\vskip0pt plus-.2\vsize
   \bigskip\smallskip\message{#1}\leftline{\bf#1}\nobreak\medskip}

%Verbatim tags:
\def\t@gsoff#1,{\def\@{#1}\ifx\@\empty\let\next=\relax\else\let\next=\t@gsoff
   \def\@@{p}\ifx\@\@@\else
   \expandafter\gdef\csname#1cite\endcsname##1{\zeigen{##1}}
   \expandafter\gdef\csname#1page\endcsname##1{?}
   \expandafter\gdef\csname#1tag\endcsname##1{\zeigen{##1}}\fi\fi\next}
\def\verbatimtags{\ifx\all\relax\else\expandafter\t@gsoff\all,\fi}
\def\zeigen#1{\hbox{$\langle$}#1\hbox{$\rangle$}}

%Equation numbering:
\def\(#1){\edef\dot@g{\ifmmode\ifinner(\hbox{\noexpand\etag{#1}})
   \else\noexpand\eqno(\hbox{\noexpand\etag{#1}})\fi
   \else(\noexpand\ecite{#1})\fi}\dot@g}

%Reference numbering:
\newif\ifbr@ck
\def\eat#1{}
\def\[#1]{\br@cktrue[\br@cket#1'X]}
\def\br@cket#1'#2X{\def\temp{#2}\ifx\temp\empty\let\next\eat
   \else\let\next\br@cket\fi
   \ifbr@ck\br@ckfalse\br@ck@t#1,X\else\br@cktrue#1\fi\next#2X}
\def\br@ck@t#1,#2X{\def\temp{#2}\ifx\temp\empty\let\neext\eat
   \else\let\neext\br@ck@t\def\temp{,}\fi
   \def\teemp{#1}\ifx\teemp\empty\else\rcite{#1}\fi\temp\neext#2X}
\def\resetbr@cket{\gdef\[##1]{[\rtag{##1}]}}
\def\references{\resetbr@cket\newsection References\par}

%Footnotes:
\newtoks\symb@ls\newtoks\s@mb@ls\newtoks\p@gelist\n@wcount\ftn@mber
    \ftn@mber=1\newif\ifftn@mbers\ftn@mbersfalse\newif\ifbyp@ge\byp@gefalse
\def\defm@rk{\ifftn@mbers\n@mberm@rk\else\symb@lm@rk\fi}
\def\n@mberm@rk{\xdef\m@rk{{\the\ftn@mber}}%
    \global\advance\ftn@mber by 1 }
\def\rot@te#1{\let\temp=#1\global#1=\expandafter\r@t@te\the\temp,X}
\def\r@t@te#1,#2X{{#2#1}\xdef\m@rk{{#1}}}
\def\b@@st#1{{$^{#1}$}}\def\str@p#1{#1}
\def\symb@lm@rk{\ifbyp@ge\rot@te\p@gelist\ifnum\expandafter\str@p\m@rk=1 
    \s@mb@ls=\symb@ls\fi\write\f@nsout{\number\count0}\fi \rot@te\s@mb@ls}
\def\byp@ge{\byp@getrue\n@wwrite\f@nsin\openin\f@nsin=\jobname.fns 
    \n@wcount\currentp@ge\currentp@ge=0\p@gelist={0}
    \re@dfns\closein\f@nsin\rot@te\p@gelist
    \n@wread\f@nsout\openout\f@nsout=\jobname.fns }
\def\m@kelist#1X#2{{#1,#2}}
\def\re@dfns{\ifeof\f@nsin\let\next=\relax\else\read\f@nsin to \f@nline
    \ifx\f@nline\v@idline\else\let\t@mplist=\p@gelist
    \ifnum\currentp@ge=\f@nline
    \global\p@gelist=\expandafter\m@kelist\the\t@mplistX0
    \else\currentp@ge=\f@nline
    \global\p@gelist=\expandafter\m@kelist\the\t@mplistX1\fi\fi
    \let\next=\re@dfns\fi\next}
\def\symbols#1{\symb@ls={#1}\s@mb@ls=\symb@ls} 
\def\bigsymbol{\textstyle}
\symbols{\bigsymbol\ast,\dagger,\ddagger,\sharp,\flat,\natural,\star}
\def\ftnumbers{\ftn@mberstrue} \def\ftsymbols{\ftn@mbersfalse}
\def\paginal{\byp@ge} \def\resetftnumbers{\ftn@mber=1}
\def\ftnote#1{\defm@rk\expandafter\expandafter\expandafter\footnote
    \expandafter\b@@st\m@rk{#1}}

%Miscellaneous macros:
\long\def\jump#1\endjump{}
\def\ssum{\mathop{\lower .1em\hbox{$\textstyle\Sigma$}}\nolimits}

\def\qed{\nobreak\kern 1em \vrule height .5em width .5em depth 0em}
\def\newneq{\hbox{\rlap{\hbox to 1\wd9{\hss$=$\hss}}\raise .1em 
   \hbox to 1\wd9{\hss$\scriptscriptstyle/$\hss}}}
\def\subsetne{\setbox9 = \hbox{$\subset$}\mathrel{\hbox{\rlap
   {\lower .4em \newneq}\raise .13em \hbox{$\subset$}}}}
\def\supsetne{\setbox9 = \hbox{$\subset$}\mathrel{\hbox{\rlap
   {\lower .4em \newneq}\raise .13em \hbox{$\supset$}}}}

%Blackboard bold:
\def\vbar{\mathchoice{\vrule height6.3ptdepth-.5ptwidth.8pt\kern-.8pt}
   {\vrule height6.3ptdepth-.5ptwidth.8pt\kern-.8pt}
   {\vrule height4.1ptdepth-.35ptwidth.6pt\kern-.6pt}
   {\vrule height3.1ptdepth-.25ptwidth.5pt\kern-.5pt}}
\def\f@dge{\mathchoice{}{}{\mkern.5mu}{\mkern.8mu}}
\def\b@c#1#2{{\rm \mkern#2mu\vbar\mkern-#2mu#1}}
\def\b@b#1{{\rm I\mkern-3.5mu #1}}
\def\b@a#1#2{{\rm #1\mkern-#2mu\f@dge #1}}
\def\bb#1{{\count4=`#1 \advance\count4by-64 \ifcase\count4\or\b@a A{11.5}\or
   \b@b B\or\b@c C{5}\or\b@b D\or\b@b E\or\b@b F \or\b@c G{5}\or\b@b H\or
   \b@b I\or\b@c J{3}\or\b@b K\or\b@b L \or\b@b M\or\b@b N\or\b@c O{5} \or
   \b@b P\or\b@c Q{5}\or\b@b R\or\b@a S{8}\or\b@a T{10.5}\or\b@c U{5}\or
   \b@a V{12}\or\b@a W{16.5}\or\b@a X{11}\or\b@a Y{11.7}\or\b@a Z{7.5}\fi}}

\catcode`\X=11 \catcode`\@=12

% % \input citeadd
%   citeadd -- a few additions for 
% files from alice that were procesed with "citealice"

\expandafter\ifx\csname citeadd.tex\endcsname\relax
\expandafter\gdef\csname citeadd.tex\endcsname{}
\else \message{Hey!  Apparently you were trying to
\string\input{citeadd.tex} twice.   This does not make sense.} 
\errmessage{Please edit your file (probably \jobname.tex) and remove
any duplicate ``\string\input'' lines}\endinput\fi

\sectno=-1   % start with sect 0
\localtags
%\verbatimtags
\NoBlackBoxes
\newbox\noforkbox \newdimen\forklinewidth
\forklinewidth=0.3pt   %%  maybe 0.6?  0.7 ?? 
\setbox0\hbox{$\textstyle\bigcup$}
\setbox1\hbox to \wd0{\hfil\vrule width \forklinewidth depth \dp0
                        height \ht0 \hfil}
\wd1=0 cm
\setbox\noforkbox\hbox{\box1\box0\relax}
\def\unionstick{\mathop{\copy\noforkbox}\limits}
\def\nonfork#1#2_#3{#1\unionstick_{\textstyle #3}#2}
\def\nonforkin#1#2_#3^#4{#1\unionstick_{\textstyle #3}^{\textstyle #4}#2}     
%
%%%   and the ``does fork'' symbol:   
\setbox0\hbox{$\textstyle\bigcup$}
%%%%%\setbox1\hbox to \wd0{\hfil$\nmid$\hfil}
\setbox1\hbox to \wd0{\hfil{\sl /\/}\hfil}
\setbox2\hbox to \wd0{\hfil\vrule height \ht0 depth \dp0 width
                                \forklinewidth\hfil}
\wd1=0cm
\wd2=0cm
\newbox\doesforkbox
\setbox\doesforkbox\hbox{\box1\box0\relax}
\def\nunionstick{\mathop{\copy\doesforkbox}\limits}
\def\fork#1#2_#3{#1\nunionstick_{\textstyle #3}#2}
\def\forkin#1#2_#3^#4{#1\nunionstick_{\textstyle #3}^{\textstyle #4}#2}     
\define\mr{\medskip\roster}
\define\sn{\smallskip\noindent}
\define\mn{\medskip\noindent}
\define\bn{\bigskip\noindent}
\define\ub{\underbar}

\define\ermn{\endroster\medskip\noindent}

\define\dbcu{\dsize\bigcup}
\define \nl{\newline}
%%%\magnification=\magstep 1
\documentstyle {amsppt}
\topmatter
\title{On what 
I do not understand (and have something to say), model theory \\
Sh702} \endtitle
\rightheadtext{Rutgers Seminars 1997}
\author {Saharon Shelah \thanks {\null\newline I would like to thank 
Alice Leonhardt for the beautiful typing. \null\newline
 Work done: mainly Fall '97 \null\newline
 First Typed - 97/Sept/12 \null\newline
 Latest Revision - 99/Oct/11} \endthanks} \endauthor 
% Previous Version - 99/Oct/9
\affil{Institute of Mathematics\\
 The Hebrew University\\
 Jerusalem, Israel
 \medskip
 Rutgers University\\
 Mathematics Department\\
 New Brunswick, NJ  USA} \endaffil
\mn
\abstract  This is a non-standard paper, containing some 
problems, mainly in model theory, which I have, in 
various degrees, been interested in.  
Sometimes with a discussion on what I have to say; sometimes, of what makes
them interesting to me, sometimes the problems are presented with a discussion
of how I have tried to solve them, and sometimes with failed tries,
anecdote and opinion.  So the discussion is quite personal, in other words,
egocentric and somewhat accidental.  As we discuss many problems, history
and side references are erratic, usually kept at a minimum (``See..." means:
see the references there and possibly the paper itself). \nl
The base were lectures in Rutgers Fall '97 and reflect my knowledge then.
The other half, concentrating on set theory, is in print \cite{Sh:666}, but
the two halves are independent.
We thank A. Blass, G. Cherlin and R. Grossberg 
for some corrections. \endabstract
\endtopmatter
\document  
% % \input alice2jlem
%% # Keywords  Input file to be used for texing Alice's files

\expandafter\ifx\csname alice2jlem.tex\endcsname\relax
  \expandafter\gdef\csname alice2jlem.tex\endcsname{}
\else \message{Hey!  Apparently you were trying to
\string\input{alice2jlem.tex}  twice.   This does not make sense.}
\errmessage{Please edit your file (probably \jobname.tex) and remove
any duplicate ``\string\input'' lines}\endinput\fi

% % \input bib4plain
\expandafter\ifx\csname bib4plain.tex\endcsname\relax
  \expandafter\gdef\csname bib4plain.tex\endcsname{}
\else \message{Hey!  Apparently you were trying to \string\input
  bib4plain.tex twice.   This does not make sense.}
\errmessage{Please edit your file (probably \jobname.tex) and remove
any duplicate ``\string\input'' lines}\endinput\fi

%  This file should be inputted if you want to use 
%  bibtex fom within plain TeX. 
      % Not really need for standard
       % bibtex files, but these commands
\def\renewcommand{\newcommand}	       % are used in our literal-unsrt.bst
\edef\cite{\the\catcode`@}%
\catcode`@ = 11
\let\@oldatcatcode = \cite
\chardef\@letter = 11
\chardef\@other = 12
%
%
% Next come some things that will be useful later.
%
% Make an outer definition into an inner one (due to Chris Thompson).
% The arguments should be the control sequence to be defined, and the
% new of the \outer control sequence, as characters; the control
% sequence #1 is defined to be just the same as \csname#2\endcsname, but
% not \outer.  For example, \@innerdef\innernewcount{newcount} would
% define \innernewcount to be a non-outer version of \newcount.
%
\def\@innerdef#1#2{\edef#1{\expandafter\noexpand\csname #2\endcsname}}%
%
% We use \@innerdef to make some of our allocations local, because
% Eplain includes our code inside a conditional.  We put @'s in the
% names to minimize the (already small) chance of conflicts.
%
\@innerdef\@innernewcount{newcount}%
\@innerdef\@innernewdimen{newdimen}%
\@innerdef\@innernewif{newif}%
\@innerdef\@innernewwrite{newwrite}%
%
%
% Swallow one parameter.
%
\def\@gobble#1{}%
%
%
% Use TeX 3.0's \inputlineno to get the line number, for better error
% messages, but if we're using an old version of TeX, don't do anything.
%
\ifx\inputlineno\@undefined
   \let\@linenumber = \empty % Pre-3.0.
\else
   \def\@linenumber{\the\inputlineno:\space}%
\fi
%
%
% The following macro \@futurenonspacelet (from the TeXbook) behaves
% essentially like \futurelet except that it discards any implicit or
% explicit space tokens that intervene before a nonspace is scanned:
%
\def\@futurenonspacelet#1{\def\cs{#1}%
   \afterassignment\@stepone\let\@nexttoken=
}%
\begingroup % The grouping here avoids stepping on an outside use of `\\'.
\def\\{\global\let\@stoken= }%
\\ % now \@stoken is a space token (\\ is a control symbol, so that
   % space after it is seen).
\endgroup
\def\@stepone{\expandafter\futurelet\cs\@steptwo}%
\def\@steptwo{\expandafter\ifx\cs\@stoken\let\@@next=\@stepthree
   \else\let\@@next=\@nexttoken\fi \@@next}%
\def\@stepthree{\afterassignment\@stepone\let\@@next= }%
%
%
% \@getoptionalarg\CS gets an optional argument from the input, enclosed
% in brackets, then expands \CS.  We set \@optionalarg to \empty if we
% don't find one, otherwise to the text of the argument.  This assumes
% the brackets don't have a funny category code.
%
\def\@getoptionalarg#1{%
   \let\@optionaltemp = #1%
   \let\@optionalnext = \relax
   \@futurenonspacelet\@optionalnext\@bracketcheck
}%
%
% The \expandafter's in this macro let us avoid the use of \aftergroup,
% which is somewhat more expensive.
%
\def\@bracketcheck{%
   \ifx [\@optionalnext
      \expandafter\@@getoptionalarg
   \else
      \let\@optionalarg = \empty
      % We can't do the \temp after the \fi, because then the \temp gets
      % in the way of reading the optional argument from the input, if
      % we do have one.
      \expandafter\@optionaltemp
   \fi
}%
\def\@@getoptionalarg[#1]{%
   \def\@optionalarg{#1}%
   \@optionaltemp
}%
%
%
% From LaTeX.
%
\def\@nnil{\@nil}%
\def\@fornoop#1\@@#2#3{}%
\def\@for#1:=#2\do#3{%
   \edef\@fortmp{#2}%
   \ifx\@fortmp\empty \else
      \expandafter\@forloop#2,\@nil,\@nil\@@#1{#3}%
   \fi
}%
\def\@forloop#1,#2,#3\@@#4#5{\def#4{#1}\ifx #4\@nnil \else
       #5\def#4{#2}\ifx #4\@nnil \else#5\@iforloop #3\@@#4{#5}\fi\fi
}%
\def\@iforloop#1,#2\@@#3#4{\def#3{#1}\ifx #3\@nnil
       \let\@nextwhile=\@fornoop \else
      #4\relax\let\@nextwhile=\@iforloop\fi\@nextwhile#2\@@#3{#4}%
}%
%
%
% This macro tests if a file \jobname.#1 exists, and sets \if@fileexists
% appropriately.  If an optional argument is given, it is used as the
% root part of the filename instead of \jobname.
%
\@innernewif\if@fileexists
\def\@testfileexistence{\@getoptionalarg\@finishtestfileexistence}%
\def\@finishtestfileexistence#1{%
   \begingroup
      \def\extension{#1}%
      \immediate\openin0 =
         \ifx\@optionalarg\empty\jobname\else\@optionalarg\fi
         \ifx\extension\empty \else .#1\fi
         \space
      \ifeof 0
         \global\@fileexistsfalse
      \else
         \global\@fileexiststrue
      \fi
      \immediate\closein0
   \endgroup
}%
%
%
%% [[[start of BibTeX-specific stuff]]]
%
% Now come the four main LaTeX commands and their associated .aux
% commands.  Just as in LaTeX, \bibliographystyle defines the BibTeX
% style name (.bst file, that is), and \bibliography defines the
% database (.bib) file(s).  The corresponding .aux-file commands are
% \bibstyle and \bibdata, which are there only for BibTeX's (but not
% LaTeX's) use.
%
\def\bibliographystyle#1{%
   \@readauxfile
   \@writeaux{\string\bibstyle{#1}}%
}%
\let\bibstyle = \@gobble
%
% As well as writing the \bibdata command to tell BibTeX which .bib
% files to read, we read the .bbl file that BibTeX (or a person,
% conceivably) has produced.  We use \bblfilebasename as the root of the
% filename to read; this defaults to \jobname.
%
\let\bblfilebasename = \jobname
\def\bibliography#1{%
   \@readauxfile
   \@writeaux{\string\bibdata{#1}}%
   \@testfileexistence[\bblfilebasename]{bbl}%
   \if@fileexists
      % We just output a non-discardable item (the `whatsit' with the
      % \bibdata command).  This means that the glue that will be
      % inserted next (\parskip or \baselineskip, most likely) will be a
      % legal breakpoint.  Most likely, this is after some kind of
      % heading, where we don't want to allow a page break.  So:
      \nobreak
      \@readbblfile
   \fi
}%
\let\bibdata = \@gobble
%
% The \nocite{label,label,...} command writes its argument to \@auxfile,
% unless instructed not to, but produces no text in the document.  Both
% the \nocite and \cite commands produce \citation commands in the .aux file.
%
\def\nocite#1{%
   \@readauxfile
   \@writeaux{\string\citation{#1}}%
}%
\@innernewif\if@notfirstcitation
%
% \cite[note]{label,label,...} produces the citations for the labels as
% well.  If the optional argument `note' is present, it's added after
% the labels.  Since \cite calls \nocite to do its .aux-file writing,
% \cite doesn't need to call \@readauxfile (\nocite does).
%
\def\cite{\@getoptionalarg\@cite}%
%
% Typeset the citations for the labels in #1, followed by the note, if
% it exists.  To change the citation's format in the text, redefine one
% or more `\print...' macros, whose defaults appear later in this file.
%
\def\@cite#1{%
   % Remember the optional argument, in case one of the macros we call
   % below ends up looking for an optional argument itself.  For
   % example, if a \cite[note] triggers reading the .aux file, then the
   % [note] would be clobbered, since \@testfileexistence looks for an
   % optional arg.
   \let\@citenotetext = \@optionalarg
   % Start printing the text, beginning with a left bracket by default.
   \printcitestart
   % It's complicated, but because \nocite puts a `whatsit' onto the list,
   % \nocite should follow \printcitestart.  It's conceivable, but very
   % unlikely, that this `whatsit' will cause a problem (glue that doesn't
   % disappear when you want it to is the most likely symptom), requiring
   % a change either to \printcitestart or to the label that the .bst file
   % produces.
   \nocite{#1}%
   \@notfirstcitationfalse
   \@for \@citation :=#1\do
   {%
      \expandafter\@onecitation\@citation\@@
   }%
   \ifx\empty\@citenotetext\else
      \printcitenote{\@citenotetext}%
   \fi
   \printcitefinish
}%
\def\@onecitation#1\@@{%
   \if@notfirstcitation
      \printbetweencitations
   \fi
   \expandafter \ifx \csname\@citelabel{#1}\endcsname \relax
      \if@citewarning
         \message{\@linenumber Undefined citation `#1'.}%
      \fi
      % Give it a dummy definition:
      \expandafter\gdef\csname\@citelabel{#1}\endcsname{%
% Change: marginal remark added, goldstrn@math.huji.ac.il, 
% goldstern@tuwien.ac.at, May 1996 mg
%  !!! change !!!
\strut
\vadjust{\vskip-\dp\strutbox
\vbox to 0pt{\vss\parindent0cm \leftskip=\hsize 
\advance\leftskip3mm
\advance\hsize 4cm\strut\openup-4pt 
\rightskip 0cm plus 1cm minus 0.5cm ?  #1 ?\strut}}
         {\tt
            \escapechar = -1
            \nobreak\hskip0pt
            \expandafter\string\csname#1\endcsname
            \nobreak\hskip0pt
         }%
      }%
   \fi
   % Now produce the text, whether it was undefined or not.
   \csname\@citelabel{#1}\endcsname
   \@notfirstcitationtrue
}%
%
% Given a label `foo', the macro `\b@foo' is supposed to
% hold the text that should be produced.
%
\def\@citelabel#1{b@#1}%
%
% So, how does a citation label get defined?  When we read the .bbl file
% (below), a \bibitem writes out a \@citedef command.  And when we read
% the \@citedef, we define \@citelabel{#1}, where #1 is the user's
% label.
%
\def\@citedef#1#2{\expandafter\gdef\csname\@citelabel{#1}\endcsname{#2}}%
%
%
% Reading the .bbl file also produces the typeset bibliography.  Please
% notice, however, that we do not produce the title for the references
% (e.g., `References'), as LaTeX does.  The formatting and spacing of
% that title, whether it should go into the headline, and so on, are all
% things determined by your format.  We cannot know those things in
% advance.  If you wish, you can define \bblhook to produce the title.
% Or just do it before the \bibliography command.
%
\def\@readbblfile{%
   % Define a counter to tell us which item number we are on, unless
   % we've already defined it (because the document has more than one
   % bibliography).
   \ifx\@itemnum\@undefined
      \@innernewcount\@itemnum
   \fi
   \begingroup
      \def\begin##1##2{%
         % ##1 is just `thebibliography'.
         % ##2 is the widest label.
         % We set (new dimen) \biblabelwidth based on the widest label
         \setbox0 = \hbox{\biblabelcontents{##2}}%
         \biblabelwidth = \wd0
      }%
      \def\end##1{}% ##1 is `thebibliography' again.
      %
      % Here we have two possibilities:
      % \bibitem[typesetlabel]{citationlabel}
      % \bibitem{citationlabel}
      % If we have the second of these, the citations are numbered, starting
      % from one; we use our own count register \@itemnum for this.
      %
      \@itemnum = 0
      \def\bibitem{\@getoptionalarg\@bibitem}%
      \def\@bibitem{%
         \ifx\@optionalarg\empty
            \expandafter\@numberedbibitem
         \else
            \expandafter\@alphabibitem
         \fi
      }%
      \def\@alphabibitem##1{%
         % Need \xdef here for various reasons.
         \expandafter \xdef\csname\@citelabel{##1}\endcsname {\@optionalarg}%
         % Left-justify alpha labels, unless \biblabel{pre,post}contents
         % are already defined.
         \ifx\biblabelprecontents\@undefined
            \let\biblabelprecontents = \relax
         \fi
         \ifx\biblabelpostcontents\@undefined
            \let\biblabelpostcontents = \hss
         \fi
         \@finishbibitem{##1}%
      }%
      \def\@numberedbibitem##1{%
         \advance\@itemnum by 1
         \expandafter \xdef\csname\@citelabel{##1}\endcsname{\number\@itemnum}%
         % Right-justify numeric labels, unless \biblabel{pre,post}contents
         % are already defined.
         \ifx\biblabelprecontents\@undefined
            \let\biblabelprecontents = \hss
         \fi
         \ifx\biblabelpostcontents\@undefined
            \let\biblabelpostcontents = \relax
         \fi
         \@finishbibitem{##1}%
      }%
      \def\@finishbibitem##1{%
         \biblabelprint{\csname\@citelabel{##1}\endcsname}%
         \@writeaux{\string\@citedef{##1}{\csname\@citelabel{##1}\endcsname}}%
         \ignorespaces
      }%
      %
      % Do the printing (we're producing the bibliography, remember).
      %
      \let\em = \bblem
      \let\newblock = \bblnewblock
      \let\sc = \bblsc
      % Punctuation won't affect spacing;
      \frenchspacing
      % the penalties below are from LaTeX's [article,book,report].sty;
      \clubpenalty = 4000 \widowpenalty = 4000
      % the next two values come from LaTeX's \sloppy command;
      \tolerance = 10000 \hfuzz = .5pt
      \everypar = {\hangindent = \biblabelwidth
                      \advance\hangindent by \biblabelextraspace}%
      \bblrm
      % the \parskip is a guess at what looks good;
      \parskip = 1.5ex plus .5ex minus .5ex
      % and the space between label and text comes from LaTeX's \labelsep.
      \biblabelextraspace = .5em
      \bblhook
      \input \bblfilebasename.bbl
   \endgroup
}%
%
% The widest label's width is useful for redefining \biblabelprint;
% you redefine \biblabelwidth, in effect, by redefining the
% \biblabelcontents macro that appears below.  And \biblabelextraspace,
% which is redefinable inside \bblhook, is added to \biblabelwidth to
% determine the amount of hanging indentation.
%
\@innernewdimen\biblabelwidth
\@innernewdimen\biblabelextraspace
%
% Now come the main macros that are related to the printing of the
% bibliography.  Since you might want to redefine them, they are given
% default definitions outside of \@readbblfile.
%
% The first one controls the printing of a bibliography entry's label.
% If you change it, make sure that it starts with something like
% \noindent or \indent or \leavevmode that puts TeX into horizontal mode
% (even if the label itself is empty); otherwise, the hanging
% indentation will get messed up in certain circumstances.
%
\def\biblabelprint#1{%
   \noindent
   \hbox to \biblabelwidth{%
      \biblabelprecontents
      \biblabelcontents{#1}%
      \biblabelpostcontents
   }%
   \kern\biblabelextraspace
}%
%
% If you are using numeric labels, and you want them left-justified
% (numeric labels by default are right-justified), do something like:
%     \def\biblabelprecontents{\relax}
%     \def\biblabelpostcontents{\hss}
%
% By default the labels are typeset in \bblrm, and enclosed in brackets.
%
\def\biblabelcontents#1{{\bblrm [#1]}}%
%
% The main text, too, is typeset using \bblrm, which is \rm by default.
%
\def\bblrm{\rm}%
%
% Emphasis for producing, e.g., titles, is done with \it by default.
%
\def\bblem{\it}%
%
% Some styles use a caps-and-small-caps font for author names.  LaTeX
% defines an \sc command but plain TeX doesn't, so we need one here.
% The definition below doesn't load the font unless it's needed, but it
% tries to load only the 10pt version, because it might not exist at
% other point sizes.
%
\def\bblsc{\ifx\@scfont\@undefined
              \font\@scfont = cmcsc10
           \fi
           \@scfont
}%
%
% The major parts of an entry are separated with \bblnewblock.  The
% numbers below are taken from LaTeX's `article' style.
%
\def\bblnewblock{\hskip .11em plus .33em minus .07em }%
%
% Here's where you stick any other bibliography-formatting goodies, or
% redefine the values above.
%
\let\bblhook = \empty
%
%
% Here are the four default definitions for formatting the in-text
% citations.  These are what you redefine (after your \input btxmac but
% before your \bibliography) to get parens instead of brackets, or
% superscripts, or footnotes, or whatever.
%
\def\printcitestart{[}%         left bracket
\def\printcitefinish{]}%        right bracket
\def\printbetweencitations{, }% comma, space
\def\printcitenote#1{, #1}%     comma, space, note (if it exists)
%
% That scheme is pretty flexible.  For example you could use
%     \def\printcitestart{\unskip $^\bgroup}
%     \def\printcitefinish{\egroup$}
%     \def\printbetweencitations{,}
%     \def\printcitenote#1{\hbox{\sevenrm\space (#1)}}
%     \font\eighttt = cmtt8
%     \scriptfont\ttfam = \eighttt
% to get superscripted in-text citations.  (The scriptfont stuff
% exists only to print an undefined citation; it's in cmtt8 because
% there is no cmtt7.)  To get something radically different, however,
% you'll have to define your own \cite command.
%
% When we read `\citation' from the .aux file, it means nothing.
%
\let\citation = \@gobble
%
%
% Now comes the stuff for dealing with LaTeX's \newcommand.  As
% mentioned earlier, this \newcommand will redefine a preexisting
% command; that's different from how LaTeX's \newcommand behaves.
%
\@innernewcount\@numparams
%
% \newcommand{\foo}[n]{text} defines the control sequence \foo to have
% n parameters, and replacement text `text'.
%
\def\newcommand#1{%
   \def\@commandname{#1}%
   \@getoptionalarg\@continuenewcommand
}%
%
% Figure out if this definition has parameters.
%
\def\@continuenewcommand{%
   % If no optional argument, we have zero parameters.  Otherwise, we
   % have that many.
   \@numparams = \ifx\@optionalarg\empty 0\else\@optionalarg \fi \relax
   \@newcommand
}%
%
% \@numparams is how many arguments this command has.  The name of the
% command is \@commandname.  The replacement text for the new macro is #1.
%
\def\@newcommand#1{%
   \def\@startdef{\expandafter\edef\@commandname}%
   \ifnum\@numparams=0
      \let\@paramdef = \empty
   \else
      \ifnum\@numparams>9
         \errmessage{\the\@numparams\space is too many parameters}%
      \else
         \ifnum\@numparams<0
            \errmessage{\the\@numparams\space is too few parameters}%
         \else
            \edef\@paramdef{%
               % This is disgusting, but \loop doesn't work inside \edef,
               % because \body isn't defined.
               \ifcase\@numparams
                  \empty  No arguments.
               \or ####1%
               \or ####1####2%
               \or ####1####2####3%
               \or ####1####2####3####4%
               \or ####1####2####3####4####5%
               \or ####1####2####3####4####5####6%
               \or ####1####2####3####4####5####6####7%
               \or ####1####2####3####4####5####6####7####8%
               \or ####1####2####3####4####5####6####7####8####9%
               \fi
            }%
         \fi
      \fi
   \fi
   \expandafter\@startdef\@paramdef{#1}%
}%
%
%% [[[end of BibTeX-specific stuff]]]
%
%
% Names of references (arguments given in the \cite and \nocite
% commands) and file names (arguments given in the \bibliography and
% \bibliographystyle commands) are recorded in \jobname.aux, called the
% \@auxfile in these macros.  Here's how they get read in.
%
\def\@readauxfile{%
   \if@auxfiledone \else % remember: \@auxfiledonetrue if \noauxfile is defined
      \global\@auxfiledonetrue
      \@testfileexistence{aux}%
      \if@fileexists
         \begingroup
            % Because we might be in horizontal mode when \@readauxfile
            % is called (if it's in response to a \cite or \nocite), we
            % want to ignore all the would-be spaces at the ends of
            % lines in the aux file.  Fortunately, it's highly unlikely
            % an end-of-line might actually be desired.
            % And because we don't change the category code of anything
            % but @, primitives like \gdef can't be used to define labels
            % in the aux file.  The solution adopted by btxmac.tex is to
            % write `\@citedef{LABEL}{DEFINITION}' to the aux file, and
            % use \csname on LABEL.
            \endlinechar = -1
            \catcode`@ = 11
            \input \jobname.aux
         \endgroup
      \else
         \message{\@undefinedmessage}%
         \global\@citewarningfalse
      \fi
      \immediate\openout\@auxfile = \jobname.aux
   \fi
}%
%
% The \@readauxfile macro does all that work the first time it's called.
% Since it's called once for every \cite, \nocite, \bibliography, and
% \bibliographystyle command that the user issues, we need to remember
% whether the work's been done.  It's considered done if we're not to do
% it---that is, if \noauxfile is defined.
%
\newif\if@auxfiledone
\ifx\noauxfile\@undefined \else \@auxfiledonetrue\fi
%
% It's conceivable you'd want to change how other characters are read;
% to do that, change their category code before doing \input btxmac.
%
%
% After reading the .aux file, \@readauxfile opens it for writing.
% The \@writeaux macro does the actual writing (as long as
% \noauxfile is undefined).
%
\@innernewwrite\@auxfile
\def\@writeaux#1{\ifx\noauxfile\@undefined \write\@auxfile{#1}\fi}%
%
%
% A macro package that uses btxmac.tex might define
% \@undefinedmessage (before doing an \input btxmac).
%
\ifx\@undefinedmessage\@undefined
   \def\@undefinedmessage{No .aux file; I won't give you warnings about
                          undefined citations.}%
\fi
%
% Even if citations are undefined, we want to complain only if
% \@citewarningtrue.  The default is to set \@citewarningtrue unless
% \noauxfile is defined.  Again, a macro package that uses
% btxmac.tex might want to redefine this.
%
\@innernewif\if@citewarning
\ifx\noauxfile\@undefined \@citewarningtrue\fi
%
%
% Finally, before leaving we restore @'s old category code.
%
\catcode`@ = \@oldatcatcode

  % This will define \cite and make sure it works as in latex

\def\widestnumber#1#2{}
  % Our amstex-ppt style does not know about \widestnumber

\def\rm{\fam0 \tenrm}

\def\fakesubhead#1\endsubhead{\bigskip\noindent{\bf#1}\par}

% % \input rsfs

% # Keywords: Script or Calligraphic (Caligraphic) letters with the RSFS Font

% The story so far:    July 1998 -- Saharon would like to have a
% ``nicer'' calligraphic font. In particualr, the leters S and P in
% the usual calligraphic font do not look ``special'' enough. 
% 
% I found out that ``rsfs'' (``Ralph Smith Formal Script'') may be
% what he wants.   I installed the mf file, the .tfm file, as well as
% a few pk files in ~/TeX/rsfs.    Let's hope that this is enough.
% Using amstex, all you have to do is to \input rsfs.tex 
% Files prepared with citealice willdothis automatically. 
%
%  Note:  for some reason xdvi calls MakeTeXpk, then Maketexpk
%  complains about wrong resolution, but still writes commands to
%  missfont.log...  
%

% we redefine a macro inside amstex's \Cal command , so that it calls
% our nice font ``rsfs'' rather than the usual calligraphic font. 
% Note thisworks for amstex only.   
% In plain tex, would have to add definitions of \Cal
% in latex... we should insteaduse mathrsfs.sty
% 

\font\textrsfs=rsfs10
\font\scriptrsfs=rsfs7
\font\scriptscriptrsfs=rsfs5

\newfam\rsfsfam
\textfont\rsfsfam=\textrsfs
\scriptfont\rsfsfam=\scriptrsfs
\scriptscriptfont\rsfsfam=\scriptscriptrsfs

\edef\oldcatcodeofat{\the\catcode`\@}
\catcode`\@11

\def\Cal@@#1{\noaccents@ \fam \rsfsfam #1}

\catcode`\@\oldcatcodeofat

\newpage

\head {Content \\
Part II} \endhead  \resetall 
 % \resetall 
\bn
\S1 $\qquad$ Two cardinal theorems and partition theorems
\mr
\item "{{}}"  [$n$-cardinal theorems (\scite{tcr.11} - \scite{tcr.16}),
$\lambda$-like for $n$-Mahlo (\scite{tcr.16A}), from finite models
(\scite{tcr.17} - \scite{tcr.20}), omitting types and Borel squares and
rectangles (\scite{tcr.21} - \scite{tcr.26}), Hanf numbers connected to
$L_{\lambda^+,\omega}$ (\scite{tcr.26} - \scite{tcr.30}).]
\endroster
\mn
\S2 $\qquad$ Monadic Logic and indiscernible sequences 
\mr
\item "{{}}"  [Monadic logic for linder orders (\scite{Mon.1} - \scite{Mon.5}
+ \scite{Mon.13}), classification by unary expansion (\scite{Mon.6} +
\scite{Mon.11} + \scite{Mon.12}), classifying by existence of 
indiscernibles and
generalizations, the properties and relatives $(k,c)-*$-stable
(\scite{Mon.7} - \scite{Mon.11c}), Borel theory (\scite{Mon.13}).]
\endroster
\mn
\S3 $\qquad$ Automorphisms and quantifiers 
\mr
\item "{{}}"  [Compact second order quantifiers (\scite{aqt.1} -
\scite{aqt.4a}).  Rigid and strongly rigid theories, pseudo decomposable
theories (\scite{aqt.5} - \scite{aqt.7}), ``all automorphisms extended", on
characterization (\scite{aqt.11}) interpreting in the automorphism group
of (free) algebra in a variety (\scite{aqt.11} - \scite{aqt.12}), properties
of abstract logics (\scite{aqt.14} - \scite{aqt.15} + \scite{aqt.17}),
second order quantifiers like (aaX)$\varphi$ 
(\scite{aqt.15a}, \scite{aqt.16}).]
\endroster
\mn
\S4 $\qquad$ Relatives of the main gap
\mr
\item "{{}}"  [Generally does the main gap characterize theories with models
characterized by invariants (\scite{cft.3a},\scite{cft.5} - \scite{cft.10});
classifying will not die.  Variation of the main gap for stable countable
theories (\scite{cft.0} - \scite{cft.1a}), pseudo elementary classes
(\scite{cft.2}) uncountable theories (\scite{cft.4}). \nl
Minimal models under embeddability (\scite{cft.12} - \scite{cft.14}), can
forcing make models isormorphic (\scite{cft.15} - \scite{cft.16}), models
up to $L_{\infty,\kappa}$-equivalence and on Karp height
(\scite{cft.17} - \scite{cft.19}).]
\endroster
\mn
\S5 $\qquad$ Classifying unstable theories
\mr
\item "{{}}"  [Dividng lines, poor man ZFC-answer (beginning +
\scite{cus.1},\scite{cus.2}, \scite{cus.2a}, \scite{cus.33}, \scite{cus.40}),
SP$(T)$ and simplicity (\scite{cus.3}), NIP, generalizing universality
spectrum.  NSOP$_n$ and tree coding (\scite{cus.13} - \scite{cus.16}).  We
look at classifying such properties (\scite{cus.17} - \scite{cus.24}; again
universality spectrum (\scite{cus.30} - \scite{cus.36}, \scite{cus.33})
about test problems, NIP (\scite{cus.41} - \scite{cus.44}) earlier
(\scite{cus.10a}).]
\endroster
\mn
\S6 $\qquad$ Classification theory for non-elementary classes
\mr
\item "{{}}"  [We ask about stability for $K_D$ (\scite{nec.1} - 
\scite{nec.5}), categoricity for $\psi \in L_{\lambda^+,\omega}$
(\scite{nec.6} - \scite{nec.12}) classification for such $\psi$
(\scite{nec.13}), $\Phi$ (\scite{nec.14}); instead of $\psi \in
L_{\lambda^+,\omega}$ we usually deal with a.e.c. (abstract elementary
classes).]
\endroster
\mn
\S7 $\qquad$ Finite model theory
\mr
\item "{{}}"  [Finding a logic (\scite{fin.1}), model theoretic content of
some 0-1 laws (\scite{fin.2}), looking for dichotomies (\scite{fin.4}),
generalized quantifiers.]
\endroster
\mn
\S8 $\qquad$ More on finite partition theorems
\mr
\item "{{}}"  [Relatives of Halse-Jewett are considered.]
\endroster
\newpage

\head {\S1 Two cardinal theorems} \endhead  \resetall \sectno=1
 % \resetall 
\bn
During the 1960's, two cardinal theorems were popular among model theorists.
\definition{\stag{tcr.1} Definition}  $(\lambda_1,\dotsc,\lambda_n)
\rightarrow_\kappa (\mu_1,\dotsc,\mu_n)$ holds if whenever $T$ is a set of
$\kappa$ f.o. sentences with unary predicates $P_1,\dotsc,P_n$ and every
finite subset of $T$ has a model $M$ such that $|P^{M_i}| = \lambda_i$ for
$i=1,\dotsc,n$, \ub{then} $T$ has a model $N$ such that $|P^{N_i}| = \mu_i$
for $i=1,\dotsc,n$.  If $\kappa$ is omitted we mean $\kappa = \aleph_0$.  For 
notational simplicity we always assume $\lambda_1 \ge \ldots \ge \lambda_n,
\mu_1 \ge \ldots \ge \mu_n$.

We shall usually speak on the case $n=2$; we have, for general discussion,
ignore the possibility of adding cardinality quantifier 
$(\exists^{\le \lambda}x)$.  
Later the subject becomes less popular; 
Jensen complained ``when I start to deal with gap $n$ 2-cardinal theorems, 
they were the epitome of model theory and as I finished, it stopped to be of
interest to model theorists". \nl
I sympathize, though model theorists has reasonable excuses: one is that they
want ZFC-provable theorems or at least semi-ZFC ones (see \cite[1.20t]{Sh:666})
the second is that it has not been clear if there were any more.
\enddefinition
\bn
\ub{\stag{tcr.2} Question}:  Are there more nontrivial $n$-cardinal ZFC
theorems, or only assuming facts on cardinal arithmetic (i.e. semi ZFC ones).
\bn
Maybe I better recall the classical ones.
\proclaim{\stag{tcr.3} Theorem}  [Vaught] $(\lambda^+,\lambda) \rightarrow
(\aleph_1,\aleph_0)$.
\endproclaim
\bigskip

\proclaim{\stag{tcr.4} Theorem}  [Chang] $\mu = \mu^{< \mu} \Rightarrow
(\lambda^+,\lambda) \rightarrow (\mu^+,\mu)$.
\endproclaim
\bigskip

\proclaim{\stag{tcr.5} Theorem}  [Vaught] $(\beth_\omega(\lambda),\lambda) 
\rightarrow (\mu_1,\mu_2)$ when $\mu_1 \ge \mu_2$. 
\endproclaim
\mn
But there were many independence results and positive theorems for $V=L$
(see \cite{Ho93}, \cite{CK}, \cite{Sch85}, \cite{Sh:18}).
%\endproclaim
\bn
After several years of drawing a blank, I found a short and easy proof of
\proclaim{\stag{tcr.6} Claim}  $(\aleph_\omega,\aleph_0) \rightarrow
(2^{\aleph_0},\aleph_0)$. \nl
In fact $(\lambda^{+ \omega},\lambda) \rightarrow (\mu_1,\mu_2)$ if
{\rm Ded}$'(\mu_2) > \mu_1 \ge \mu_2$ where
\endproclaim
\bigskip

\definition{\stag{tcr.7} Definition}  Ded$'(\mu) = \text{ Min}\{\lambda:
\text{ if } {\Cal T}$ is a tree with $\lambda$ nodes and $\delta \le \lambda$
levels, then the number of its $\delta$-branches is $< \lambda\}$.  This is
essentially equal to

$$
\align
\text{Ded}(\mu) = \text{ Min} \bigl\{\lambda:&\text{if } I
\text{ is a linear order of cardinality } \mu \\  
  &\text{then } I \text{ has } < \lambda \text{ Dedekind cuts} \bigr\}.
\endalign
$$
\enddefinition
\bn
See \cite{Sh:49}.  Considering the many 
independence proofs and natural limitations, one may ask (\cite{CK}) \nl
\ub{\stag{tcr.8} Question}:  Assume $\lambda = \lambda^{\beth_\omega(\mu)}$
and $\lambda_1 = (\lambda_1)^{< \lambda_1} \ge \mu$,
do we have $(\lambda^+,\lambda,\mu) \rightarrow (\lambda^+_1,\lambda_1,
\mu_1)$.  \nl
Things are not commutative, if $\mu = \mu^{< \mu}$ then
$(\beth_\omega(\mu^+),\mu^+,\mu)
\rightarrow (\lambda,\mu^+,\mu)$ is easy and well known (a consequence of
\scite{tcr.4} + \scite{tcr.5}). \nl
In fact, the impression this becomes set theory has some formal standing: we
know that all such theorems are provably equivalent to suitable partition 
theorems, for formalizing this we need the following definition.  
\bigskip

\definition{\stag{tcr.9} Definition}  1) Let $E$ be an equivalence relation
on ${\Cal P}(n)$ preserving cardinality; we call such a pair $(n,E)$ an 
identity.  Let $\lambda \rightarrow (n,E)_\mu$ mean that if 
$F_\ell:[\lambda]^\ell
\rightarrow \mu$ for $\ell \le n$, \ub{then} we can find $\alpha_0 < \ldots 
< \alpha_{n-1} < \lambda$ such that for any $u,v \in [n]^k,k \le n$ we have:

$$
u E v \Rightarrow F_k(\dotsc,\alpha_\ell,\dotsc)_{\ell \in u} =
F_k(\dotsc,\alpha_\ell,\dotsc)_{\ell \in v}
$$
\mn
we call $(n,E)$ an identity of $(\lambda,\mu)$. \nl
2) Id$(\lambda,\mu) =: \{(n,E):(n,E) \text{ is an identity of } 
(\lambda,\mu)\}$.
\enddefinition
\bn
Now
\proclaim{\stag{tcr.10} Claim}  Essentially assuming
$\lambda > \mu,\lambda_1 \ge \mu_1 \ge \kappa$ we have:
$(\lambda,\mu) \rightarrow_\kappa (\lambda_1,\mu_1)$ iff Id$(\lambda,\mu) 
\supseteq \text{ Id}(\lambda_1,
\mu_1)$ (see \cite{Sh:8}, \cite{Sh:E17}). 
\sn
Fully: if $\mu_1 = \mu^{\aleph_0}$ or just $(\lambda_1,\mu_1) 
\rightarrow_{\aleph_0} (\lambda_1,\mu_1)$ then the equivalence holds; the
implication $\Rightarrow$ holds always.
\endproclaim
\bn
This leaves open: \nl
\ub{\stag{tcr.11} Question}:  Prove the consistency of the existence of
$\lambda \ge \mu$ such that $(\lambda,\mu) \nrightarrow (\lambda,\mu)$
(another formulation is: $(\lambda,\mu)$ is not $\aleph_0$-compact). 
\bn
\ub{\stag{tcr.11a} Discussion}:  I am sure that the statement in
\scite{tcr.11} is consistent.  Note that all the cases we mention gives
the $\aleph_0$-compactness (and a completeness theorem).
\sn
Originally the theorems quoted above were not proved in this way.

Vaught proved \scite{tcr.3} by (sequence)-homogeneous models.  Chang proved
\scite{tcr.4} by saturated models of suitable expansion of $T$. \nl
Vaught \scite{tcr.5} finds a consistent expansion $T_1$ of $T$ which has a
built-in elementary extension increasing $P_1$ but preserving $P_2$.  Morley
used Erd\"os Rado theorem to give an alternative proof.  Now $(\aleph_\omega,
\aleph_0) \rightarrow (2^{\aleph_0},\aleph_0)$ was proved this way.  It took 
me some effort to characterize the identities for the pair 
$(\aleph_1,\aleph_0)$, see \cite{Sh:74}, so it gives an alternate proof.
\sn
Surely Jensen's proof of his 2-cardinal theorems can be analyzed in this way,
but I have not looked at this.
\bn
Now Jensen's proofs in this light, say 
\proclaim{\stag{tcr.12} Theorem}  1) Fixing $n$, if we look at what can be
Id$(\lambda^{+n},\lambda)$, when $V=L$, it is minimal. \nl
2) If $V=L$, then Id$(\mu^+,\mu),\mu$ singular, (e.g. Id$(\aleph_{\omega+1},
\aleph_\omega))$ is equal to Id$(\aleph_1,\aleph_0)$ hence is minimal. \nl
\endproclaim
\bn
This fits the intuition that $L$ tends to have objects.  So there are many 
colorings in this case. \nl
So we can ask \nl
\ub{\stag{tcr.13} Question}:  Fixing a pattern of cardinal arithmetic, what is
the minimal possible set Id$(\lambda,\mu)$ if it exists?  (Minimal: varying
on forcing extensions giving such patterns).  As equivalent formulation is:
what identities are provable?  E.g. $\mu = \mu^{< \mu},2^\mu = \mu^{++},
\lambda = \mu^{+3} = 2^{\mu^+}$.
\mn
The idea is:  if we 
lose hope that all such pairs have the same set of identities,
resolvable in ZFC, can we at least find minimal pairs.  We may instead of
cardinal arithmetic use e.g. ``there is a kurepa tree" or whatever, but this
is less appealing to me.
\bn
It is natural to ask also: \nl
\ub{\stag{tcr.14} Question}:  Fixing a pattern of cardinal arithemtic, what
is the maximal possible for set Id$(\lambda,\mu)$?
\sn
A very natural case is $\lambda = \beth_n(\mu),n \ge 2$.  In fact, I think
it is almost sure that the following case gives it.  Let $\lambda_0 >
\lambda_1 > \ldots > \lambda_n$, with each $\lambda_\ell$ is supercompact 
Laver indestructible, now force by
$\dsize \prod_{\ell < n} P_\ell$ where $P_\ell$ is adding $\lambda_{\ell +1}$
Cohens subsets of $\lambda_\ell$.  I think $(\lambda_n,\lambda_0)$ in this 
model has
a maximal set of identities.  The point is that each $\lambda_{\ell +1}$ 
satisfies a
generalization of Halpern-Lauchli theorem (see \cite[\S4]{Sh:288}).
\bn
\ub{\stag{tcr.15} Question}:  Assume GCH, $\mu$ singular limit of 
supercompacts.  Is Id$(\mu^+,\mu)$ maximal?  \nl
Jensen had found the minimal \scite{tcr.2}; now see \cite{MgSh:324}, there
is no $\mu^+$-tree for $\mu$ as above, so it is a natural candidate for 
maximality.
\bn
\ub{\stag{tcr.16} Question}:  1) What is the maximal set of identities
Id$(\lambda^{+n},\lambda)$ under GCH? \nl
2) Can we have a universe of set theory satisfying 
GCH $+ \dsize \bigwedge_\lambda \text{ Id}
(\lambda^{+n},\lambda)$ maximal? \nl
3) Similarly for $(\beth_n(\lambda),\lambda))$.  \nl
For (2),(3) we need ``GCH fails everywhere (badly)", see
Foreman Woodin \cite{FW}.
\nl
Generally, our knowledge on the family of forcing doing something for all 
cardinals seems not to be developed flexibly enough now (see \cite{Sh:666}).
\bn
\ub{\stag{tcr.16A} Question}:  1) If $\lambda$ is strongly inaccessible
cardinal and $\lambda_1$ is an inaccessible non-Mahlo cardinal which has a
square (or even $V=L$) and the first order $\psi$ has a $\lambda$-like
model \ub{then} $\psi$ has a $\lambda$-like model? \nl
2) Similarly with $\lambda$ being $n$-Mahlo, $\lambda_1$ being not 
$(n+1)$-Mahlo (and $V=L$) (see \cite{Sch85}). \nl
We know it is surely true (at least if $V=L$), but this is not a proof.
The singular case is Keisler \cite{Ke68} (and
more in \cite{Sh:18}).
\bn
\centerline{$* \qquad * \qquad *$}
\bn
We can ask \nl
\ub{\stag{tcr.17} Problem}:  When do we have $\kappa \in \{2,\aleph_0\}$ and
$\{(m_{i,1},\dotsc,m_{i,n}):i < \omega\} \rightarrow_\kappa (\mu_1,\dotsc,
\mu_n)$ with $m_{i,\ell} < \omega$ which means: 
\mr
\item "{$(*)$}"  if $T$ is first order theory of cardinality $\le \kappa$ and
every finite $T' \subseteq T$ for $i$ large enough have a model $M_i$
such that $|P^{M_i}_\ell| = m_{i,\ell}$ for $\ell = 1,\dotsc,n$, 
\ub{then} $T$ has a model $M,|P^M_\ell| = \mu_\ell$ for $\ell = 1,\dotsc,n$.  
\nl
We know something, see \cite{Sh:37}, 
\cite[p.250-1]{Sh:18}.  (We ignore here the can of disjunctions; for every $i$
large enough for some $j,T'$ has a model $M,|P^{M_\ell}| = m_{i,j,\ell}$). 
\endroster
\bn
\proclaim{\stag{tcr.18} Claim}:  If 
$i \le m_{i,1}$ and $(m_{i,1})^i \le m_{i,2}$ then

$$
\{(m^i_1,m^i_2):i < \omega\} \rightarrow (\aleph_0,2^{\aleph_0}).
$$
\endproclaim
\bn
\ub{\stag{tcr.19} Question}:  For $m_{i,1},m_{i,2}$ as in \scite{tcr.18} do
we have always (i.e. for every $\lambda$, provably in ZFC)
$\{(m^i_1,m^i_2):i < \omega\} \rightarrow (2^\lambda,\lambda)$?  Or at least
$\{(2^m,m):m < \omega\} \rightarrow (2^\lambda,\lambda)$.
\bn
(The problem is when Ded$(\lambda) < (2^\lambda)^+$).  Those problems
(\scite{tcr.17} - \scite{tcr.20}) are involved with problems in (finitary)
Ramsey theory.
Natural (and enough) to try
to show consistency of (for $T$ with Skolem functions)

$$
\align
T_{Sk} &\cup \{x_\eta \ne x_\nu \and P_1(x_\eta):\eta \in {}^\lambda 2,
\eta \ne \nu \in {}^\lambda 2\} \\
  &\cup \{P_2(\sigma(x_{\eta_1},\dotsc,x_{\eta_n})) \rightarrow \sigma
(x_{\eta_1},\dotsc,x_{\eta_n}) = 
 \sigma(x_{\nu_1},\dotsc,x_{\nu_n}): \\
  &n < \omega,\sigma \text{ a term and }
\langle \eta_1,\dotsc,\eta_n \rangle \approx \langle \nu_1,\dotsc,\nu_n 
\rangle\}
\endalign
$$
\mn
where $\langle \eta_0,\dotsc,\eta_{n-1} \rangle \approx \langle \nu_0,\dotsc,
\nu_{n-1} \rangle$, for $\eta_\ell,\nu_\ell \in {}^\lambda 2$, means
\mr
\widestnumber\item{$(*)(a)$}
\item "{$(*)(a)$}"  $\eta_\ell <_{lex} \eta_k \equiv \nu_\ell <_{lex} \nu_k$,
of course $<_{\ell ex}$ is lexicographic order
\sn
\item "{$(b)$}"  $\ell g(\eta_{\ell_1} \cap \eta_{k_1}) < \ell g(\eta_{\ell_2}
\cap \eta_{k_2}) \Leftrightarrow \ell g(\nu_{\ell_1} \cap \nu_{k_1}) < 
\ell g(\nu_{\ell_2} \cap \nu_{k_2})$
\sn
\item "{$(c)$}"  $\eta_m(\ell g(\eta_\ell \cap \eta_k)) = 
\nu_m(\ell g(\nu_\ell \cap \nu_k))$.
\endroster
\bn

(Main Point: level of the splitting not important, unlike the proof of the
previous theorem \scite{tcr.18}). \nl
This approach tells us to find more identities for the relevant finite pairs.
We can, on the other hand, try to exploit that ``Id$(2^\lambda,\lambda)$ is
smaller than suggested by the above approach" (see \cite[3.4,6.3]{Sh:430}).
\bn
\ub{\stag{tcr.20} Question}:  Does, for $W \subseteq \omega$ infinite,
$n < \omega$

$$
\{(\beth_n(i),i):i \in W\} \rightarrow (\beth_n(\lambda),\lambda)?
$$
\mn
or even

$$
\{ \bigl( (\beth_n(i))^i,i\bigr):i \in W\} \rightarrow 
(\beth_{n+1}(\lambda),\lambda)?
$$
\mn
Some of the theorems above have also parallel with omitting types.  So
considering some parallelism it is very natural to ask \nl
\ub{\stag{tcr.21} Question}:  If $\psi \in L_{\omega_1,\omega}$ has a model
of cardinality $\ge \aleph_{\omega_1}$ does it have a model of cardinality
continuum? (well assuming $2^{\aleph_0} > \aleph_{\omega_1}$).
\bn
This is connected to the problem of Borel squares, a problem I had heard
from Harrington about.
\definition{\stag{tcr.22} Definition}  1) A set $B \subseteq
{}^\omega 2 \times {}^\omega 2$ contains a $\lambda$-square if for some
$A \subseteq {}^\omega 2$ of cardinality $\lambda$ we have $A \times A
\subseteq B$ i.e. $\eta,\nu \in A \Rightarrow (\eta,\nu) \in B$. \nl
2) A set $B \subseteq {}^\omega 2 \times {}^\omega 2$ contains a
perfect square if there is a perfect set ${\Cal P} \subseteq {}^\omega 2$
such that ${\Cal P} \times {\Cal P} \subseteq B$. \nl
3) A set $B \subseteq {}^\omega 2 \times {}^\omega 2$ contains a
$\lambda$-rectangle if for some $A_1,A_2 \subseteq {}^\omega 2$ of
cardinality $\lambda$ we have $A_1 \times A_2 \subseteq B$.  We add perfect
if $A_1,A_2$ are perfect.
\enddefinition
\bn
The connection is (see \cite{Sh:522}).
\proclaim{\stag{tcr.23} Claim}  Assume $MA + 2^{\aleph_0} > \aleph_{\omega_1}$,
for some cardinal $\lambda^*$ we have
\mr
\item "{$(a)_1$}"  if $\psi \in L_{\omega_1,\omega}$ has a model of
cardinality $\ge \lambda^*$ then it has a model of cardinality continuum
\sn
\item "{$(a)_2$}"  for no $\lambda' < \lambda^*$ does $(a)_1$ hold
\sn
\item "{$(b)_1$}"  if $\lambda^* < 2^{\aleph_0}$ and ${\Cal B}$ 
is a Borel subset of ${}^\omega 2 \times
{}^\omega 2$ and it has a $\lambda^*$-square, \ub{then} ${\Cal B}$ contains
a perfect set
\sn
\item "{$(b)_2$}"  for no $\lambda' < \lambda^*$ does $(b)_1$ hold
\sn
\item "{$(c)$}"  if $\lambda^* < 2^{\aleph_0}$ then $\lambda^*$ is a limit
cardinal of cofinality $\aleph_1$.
\endroster
\endproclaim
\bn
In fact this $\lambda^*$ essentially
can be defined as $\lambda_{\aleph_1}(\aleph_0)$ where
\definition{\stag{tcr.24} Definition}  1) For a model $M$ with countable
vocabulary, we define

$$
\text{rk}_\mu:\{w \subseteq M:w \text{ finite nonempty}\} \rightarrow
\text{ Ord } \cup \{\infty\}
$$
\mn
(really rk$_{M,\mu}$) by

$$
\align
\text{rk}_\mu(w) \ge \alpha +1 &\text{ iff, for any enumeration }
\langle a_\ell:\ell < |w| \rangle \text{ of } w \\
  &\text{ and first order formula } \varphi(x_0,x_1,\dotsc,x_{n-1}) \in
L_{\tau(M)} \text{ such that} \\
  &M \models \varphi[a_0,a_1,\dotsc,a_{n-1}] \text{ we can find } 
\ge \mu \text{ members} \\
  &a'_0 \in M \backslash \{a_0\} \text{ such that} \\
  &M \models \varphi[a'_0,a_1,a_2,\dotsc,a_{n-1}] \text{ and rk}(w \cup
\{a'_0\}) \ge \alpha.
\endalign
$$
\mn
2) $\lambda_{\mu,\alpha}(\aleph_0) = \text{ Min}\{\lambda:\text{ if } M
\text{ is a model of cardinality } \lambda$ and countable vocabulary then
$\alpha \le \text{ sup}\{\text{rk}_\mu(w)+1:w \subseteq M$ finite
nonempty$\}$.  We may omit $\mu$ if $\mu = 1$.
\enddefinition
\bn
So question \scite{tcr.21} can be rephrased as  \nl
\ub{\stag{tcr.25} Question}:  If $\lambda_{\omega_1}(\aleph_0) =
\aleph_{\omega_1}$?

It is harder but we can deal similarly with rectangles and with equivalence
relations (see \cite{Sh:522} and hopefully \cite{Sh:532}); 
so e.g.
\mn
\ub{\stag{tcr.26} Question}: If a Borel set ${\Cal B} \subseteq {}^\omega 2
\times {}^\omega 2$ contains an e.g. $\aleph_{\omega_1}$-rectangle (i.e. a
$A_1 \times A_2,|A_1| = |A_2| = \aleph_{\omega_1}$) then does it contain a
perfect rectangle? 
\bn
\centerline {$* \qquad * \qquad *$}
\bn
On Hanf numbers of omitting types and relatives see Grossberg Shelah
\cite{GrSh:259}.  Let
$\delta_2(\lambda,\kappa)$ be the minimal ordinal 
$\delta$ such that if $\psi \in
L_{\kappa^+,\omega}$ has a model $M$, otp$(M,<^M) \ge \delta,|P^M| =
\lambda$, \ub{then} $\psi$ has a non-well ordered model $N$ such that
$N \restriction P^N \prec M \restriction P^M$.
\bn
\ub{\stag{tcr.27} Question}:  If 
$\lambda > 2^\kappa$, cf$(\kappa) \ge \aleph_0$ do we
have $\delta_2(\lambda,\kappa) = (\text{cov}(\lambda,\kappa)+ 2^\kappa)^+$?
\bn
\ub{\stag{tcr.30} Question}:  Let cf$(\kappa) > \aleph_0$; is 
$\delta_2(\kappa,\kappa) < (\text{sup}\{\text{rk}_D(f):D \text{ an }
\aleph_1$-complete filter on $\kappa,f \in {}^\kappa \kappa\})^+$?
\newpage

\head {\S2 Monadic Logic and Indiscernible} \endhead  \resetall \sectno=2
 % \resetall 
\bn
On monadic logic generally see Gurevich \cite{G} (till '81).

We almost know how complicated the monadic theory of the real line is:
of course, it is interpretable in the 2nd order theory of $2^{\aleph_0}$,
while we can interpret in it the second order theory of $2^{\aleph_0}$ in
$V^{\text{Cohen}}$ (Boolean interpretation - probably the reason it 
(the undecidability of the monadic theory of $(\Bbb R,<)$) was
difficult is that first order interpretation was expected; but it takes 
more years to see that this speaks on forcing.  We cannot represent 
syntactically $\Bbb N$, but we can represent Cohen names of natural numbers), 
see latest version \cite{Sh:284a}.
\bn
\ub{\stag{Mon.1} Question}:  1) Can we
\mr
\item "{$(a)$}"  interpret the monadic theory of (the order) $\Bbb R$ in
(second order theory of $2^{\aleph_0})^{V^{\text{Cohen}}}$? \nl
\ub{or} just show 
\sn
\item "{$(b)$}"  Turing degree(monadic theory of $\Bbb R$) $\le$ Turing 
degree(second order theory of $2^{\aleph_0})^{V^{\text{Cohen}}}$?
\ermn
There are many variants.
\bigskip

\definition{\stag{Mon.1A} Definition}  1) For a logic ${\Cal L}$, 
Th$_{\Cal L}(M)$, the
${\Cal L}$-theory of the structure $M$ in the universe $V$ is $\{\varphi:
\varphi \in {\Cal L}$ in the vocabulary of $M$ and in $V$ we have
$M \models \varphi\}$.
\nl
2) When ${\Cal L}$ is a logic, ${\Cal L}(Q_t)_{t \in I}$ means we add the
quantifiers $Q_t,{\Cal L}_{\lambda,\kappa}$ means we allow (forming the
formulas) take conjunctions on $< \lambda$ formulas and use a string of
$< \kappa$ quantifiers.  But we may use 
$L = L_{\omega,\omega}$ for first order,
so $L_{\lambda,\kappa},L(Q_t)$ for the expansions as above.
\enddefinition
\bn
You may ask: \nl
\ub{\stag{Mon.1B} Question}:  How are the $L$(2nd)-th theory of 
$2^{\aleph_0}$ in $V$ and in $V^{\text{Cohen}}$ related? 
Of course, 2-nd stand for the quantifier on say arbitrary binary relations.
\nl
This is a different question - how many times are they equal, e.g. if
$V = V^{\text{Cohen}}_0$, then they are equal.
\mn
\relax From the point of view of monadic logic, the question I think is: can we
``eliminate quantifiers" using names, and the answer ``they are equal" to 
the second question (\scite{Mon.1}(b)) may be accidental, in the sense that
does not answer ``can monadic formulas say more than the appropriate forcing 
statements".  (They may be one definable from the other...)
\bn
We may also ask, (more specifically than in \scite{Mon.1B}) \nl
\ub{\stag{Mon.2} Question}:  Can the monadic theory of $\Bbb R$ be changed by
adding Cohen?  What if we assume $V=L$?

As indicated, the hope is a ``meaningful" reduction of monadic formulas to
relevant forcing statement.  If we try for other direction, it is natural to
try to interpret the second order theory of $2^{\aleph_0}$ in $V^Q$ for 
$Q$ another forcing, e.g. Sacks forcing.

It is reasonable to try to deal with a similar problem where the upper and
lower bounds are further apart.  Consider $M_\lambda = ({}^{\omega >}
\lambda,\triangleleft)$ in the logic $L(Q_{pr})$, where $Q_{pr}$ is
the quantifier over pressing down unary function $f$, where pressing down
means ``$f(x)$ is an initial segment of $x$". \nl  
Alternatively, ask on the monadic theory of

$$
M_\lambda = ({}^{\omega >}\lambda,\triangleleft,+,\times)
$$

$$
\triangleleft = \{(\eta,\nu):\eta \text{ an initial segment of } \nu
\text{ both in } {}^{\omega >} \lambda\}
$$

$$
+ = \{(\eta,\nu,\rho):\eta \triangleleft \rho,\nu \triangleleft \rho,
\text{ all three in } {}^{\omega >} \lambda \text{ and }
\ell g(\eta) + \ell g(\nu) = \ell g(\rho)\}
$$

$$
\align
\times = \{(\eta,\nu,\rho):&\,\,\eta,\nu,\rho \text{ belongs to} \\
  &\{\eta^* \restriction n:n < \omega\} \text{ for some } \eta^* \in
{}^\omega \lambda \text{ and} \\ 
  &\ell g(\eta) \times \ell g(\nu) = \ell g(\rho)\}
\endalign
$$
\bn
Now (see \cite{Sh:205})
\proclaim{\stag{Mon.3} Theorem}  In the $L(Q_{pr})$-theory of 
$M_\lambda$, we can interpret the
Levy$(\aleph_0,\lambda)$-Boolean valued second order theory of $\lambda =$
second order theory of $\aleph_0$ in $V^{\text{Levy}(\aleph_0,\lambda)}$.
So the complexity of the $L(Q_{pr})$-theory of $M_\lambda$ is at 
most that of the second order theory of $\lambda$ and at least that of
the second order theory of
$\lambda$ in $V^{\text{Levy}(\aleph_0,\lambda)}$.  
\endproclaim
\bn
(Note: this is just second order theory of $\aleph_0$ which stabilize under
large cardinals).  This depends on $\lambda$ because in second order theory
$\lambda^{V^{\text{Levy}(\aleph_0,\lambda)}}$ we can e.g. 
interpret f.o. theory
of $(L_{\lambda^+},\in)$.  So not unnatural to assume that the same is true on
the $L(Q_{pr})$-theory of $M_\lambda$, this is true, e.g. if 
Th$_{L(Q_{\text{pr}})}(M_\lambda)$ is
interpretable in the Levy$(\aleph_0,\lambda)$-Boolean valued second order
theory of $\lambda$, that is $\aleph_0$.
%\endproclaim
\bn
\ub{\stag{Mon.4} Problem}:  The parallel of \scite{Mon.1}(b), \scite{Mon.3}
for $L(Q_{pr})$.

We know that the monadic theory of linear order is complicated, exactly as
second order theory (so they have the same Lowenheim number).  Is there a
sizable class where we can have simple monadic theory?  
\bn
\ub{\stag{Mon.5} Problem}:  Can the monadic theory of well orders be
decidable?  And/or has a small Lowenheim number?  Even $\aleph_\omega?$
\nl
(Why ``can" not ``is"?  Consistently monadic theory of $(\omega_2,<)$ is as
complicated as you like (\cite{GMSh:141}, \cite{LeSh:411}).  
Note that the statement ``every stationary $S
\subseteq S^2_0$ reflect" can be expressed in monadic logic on $(\omega_2,<)$,
hence the theory is ``set theoretically sensitive".  There are theorems
saying that there is a strong connection.)

There is a natural candidate for such a 
model of set theory, but it is not known if it works.  
The consequence will be that also the Lowenheim number of well ordering and
the Lowenheim numbers of the class of linear orders are small. 

The candidate we mention is: let $\bold V_0$ satisfies GCH, we shall force
with $P_\infty = \dbcu_\alpha P_\alpha$ where we use an iterated forcing
$\langle P_\alpha,{\underset\tilde {}\to Q_\alpha}:\alpha$ an ordinal
$\rangle$ with full support with ${\underset\tilde {}\to Q_i}$ defined as
$Q^{\lambda_i}$ in $\bold V = \bold V^{P_i},\lambda_i =$ the ith regular
uncountable cardinal in $\bold V_0$, defined as below.  
In universe $\bold V$ with a
cardinal $\lambda = \lambda^{< \lambda}$, let $Q^\lambda$ be the result of
iteration of length $\lambda^+,\langle Q_{\lambda,i},
{\underset\tilde {}\to R_{\lambda,i}}:i < \lambda^+ \rangle,Q^\lambda =
\dbcu_{i < \lambda^+} Q_{\lambda,i},R_{\lambda,i}$ has cardinality $\lambda$
and has an extra partial order $\le_{pr} = \le^{R_{\lambda,i}}_{pr}$ such
that $p \le_r q \Rightarrow p < q$ and if $\delta < \lambda^+$ is limit,
$\langle p_i:i < \delta \rangle$ is $\le^{R_{\lambda,i}}_{pr}$-increasing
continuous then it has a $\le_{pr}$-lub and for every dense open ${\Cal I}
\subseteq R^{\lambda,i}$ and $p \in R^{\lambda,i}$ there is $q$
satisfying $p \le^{R_{\lambda,i}}_{\text{pr}} q \in {\Cal I}$.  
This forcing is easy to handle and
add e.g. many non reflecting stationary sets (e.g. use for regular $\lambda
> \aleph_0,R = \{h:h \text{ is a function from some } \alpha_h < \lambda$
to $h^{-1}\{1\}$ do not reflect$\},h_1 \le h_2 \Leftrightarrow h_1 \subseteq
h_2$ and $h_1 \le_{pr} h_2 \Leftrightarrow h_1 = h_2 \vee (h_1 \subseteq
h_2 \and h_2(\alpha_{h_1})=0$).

The analysis of the monadic theory I expect uses the lemmas (and notions)
of \cite[\S4]{Sh:42}.
\bn
\centerline {$* \qquad * \qquad *$}
\bn
Suppose we fix a first order theory $T$ (e.g. countable),
look at monadic logic on its class of models.  There was much research on the
monadic theory of linear orders and trees.  Why?  Just accident?
(see Baldwin Shelah \cite{BlSh:156}).
\bn
\ub{\stag{Mon.6} Problem}:  Let $T$ be first order complete.
If we cannot (f.o.) interpret second order theory in the monadic theory of
model of $T$, \ub{then} models of $T$ are not much more complicated 
than trees.
\mn
Note: if in some model $M$ of $T$ expanded 
by unary predicates call it $M^+$, we
can interpret a one to one function $H:A \times B \rightarrow 2$ where
$A,B$ are infinite, \ub{then} the theory is at least as complicated as
second order logic, so those are hopelessly complicated for the purpose of 
our present investigation.  Assume not, that is
\mr
\item "{$(*)$}"  for any $M^+,Th_L(M^+)$ does not have the independence
property.
\ermn
So we feel the cut is meaningful, a dividing line.  We shall return to this
later (\scite{Mon.11}) because this connects somehow to another problem
also on classifying f.o. theories suggested by Grossberg and Shelah
(observing $(*)$ below):
\bn
\ub{\stag{Mon.7} Problem}:  Investigate $\rightarrow_T$ according to
properties of $T$, where $T$ is a complete first order theory, where
\bigskip

\definition{\stag{Mon.8} Definition}  1) Let $\lambda \rightarrow_T 
(\mu)_\kappa$ mean that:  
if $M \models T,A \in [M]^\kappa,\bar a_i \in {}^k M$ 
for $i < \lambda$, \ub{then} for some $Y \subseteq \lambda,|Y| = \mu$, the
sequence $\langle \bar a_i:i \in Y \rangle$ is
indiscernible over $A$ in $M$. \nl
2) Let $\lambda \rightarrow^{\text{loc}}_T(\mu)_\kappa$ mean that for any
finite set $\Delta$ of formulas, we get above $\Delta$-indiscernibility. \nl
3) We may replace $T$ by $K$ for a class of models, or by $M$ if $K = \{M\}$.
\enddefinition
\bigskip

\definition{\stag{Mon.9} Definition}  $T$ has the $\omega$-independence
property if there are $k < \omega$ and formula $\varphi_n(\bar x_1,\dotsc,
\bar x_n,\bar y_n)$ for $n < \omega$ where $\ell g(\bar x_i) = k$ 
such that for
every $\lambda$ and $F:[\lambda]^{< \aleph_0} \rightarrow 2$ there are
$M \models T,\bar a_i \in {}^k M$ and $\bar b_n \in {}^{\ell g(\bar y_n)}M$ 
such that: $M \models \varphi_n(\bar a_{i_0},\dotsc,\bar a_{i_{n-1}},
\bar b_n)$ iff $F(\{i_0,\dotsc,i_{n-1}\}) = 1$ (see \cite{LwSh:560}).
\sn
So
\mr
\item "{$(*)$}"   for $T$ with the $\omega$-independence property
$\lambda \rightarrow_T (\mu)_{\aleph_0}$ is equivalent to 
$\lambda \rightarrow (\mu)^{< \omega}_{\aleph_0}$ so they are maximally
complicated under this test.
\endroster
\enddefinition
\bn
\ub{\stag{Mon.10} Problem}:  If $T$ doesn't have the independence property
(= have NIP), $\rightarrow_T$ is ``nice" because supposedly the prototypes 
of the class of unstable theories having NIP is linear order, for which 
$\rightarrow_T$ has a nice theory
(as we can go down to well ordering).
\bn
We expect a nice solution.  
The problem (\scite{Mon.7}) may be partially resolved by 
an answer to \scite{cus.42}.  Though the last two 
problems remain open, we can use a 
weak answer to the last to give some information on the earlier one.
\bigskip

\definition{\stag{Mon.10A} Definition}  1) Let $\lambda \rightarrow_{T,m}
(\mu)_\kappa$ be defined as in \scite{Mon.8} restricting ourselves to $k$ 
such that $k < 1+m$ (so for $m = \omega$ we get \scite{Mon.8}). \nl
2) Let $\lambda \rightarrow^{\text{loc}}_{T,m}(\mu)_\kappa$ means that for any
finite set $\Delta$ of formulas, we get above $\Delta$-indiscernibility.

Well \scite{Mon.10A}(2) is, of course, interesting only when the Erd\"os-Rado
Theorem does not give the answer. Now you may ask: will it make a
difference to demand $k=1$. Surprisingly there is: it suffices to have
``no $\varphi(x,y;\bar z)$ has the order property in $M$" to get
strong results on $\rightarrow_M$ (see \cite[Ch.I,\S4]{Sh:300}). More
elaborately, the surprise for me was that the condition like ``no
$\varphi(\bar x,\bar y;\bar z)$ has the order property" when
restricting $\ell g(\bar x) = \ell g(\bar y) = k$ but not $\ell g(\bar
z)$ has any consequences (some readers missed the point that the model
was not required to be stable), even $T$ was not required to be
stable, but it is less interesting (\cite[np1.11t]{Sh:715}). 
\enddefinition
\bigskip

\definition{\stag{Mon.10B} Definition}  1) We say $T$ is $(k,r)-*$-NIP if
every formula $\varphi = \varphi(\bar x,\bar y,\bar z)$ with $\ell g(\bar x) 
= k,\ell g(\bar y) = r$ is $(k,r)-*$-NIP which means that:
for no $\bar a_u,\bar b_\ell,\bar c$ for $u \subseteq \omega,\ell < \omega$ 
do we have ${\frak C} \models \varphi[\bar a_u,\bar b_k,\bar c]$ iff 
$\ell \in u$ (so $\ell g(\bar a_\ell) = k,\ell g(\bar b_k)
=r$, can be phrase by a variant of $|S^k_\varphi(A)|$ small).  We may 
replace $(\ell,m)$ by a set of such pairs.  \nl
2) Similarly for other ``straight" properties, see \scite{cus.17},
particularly part (4), \scite{cus.18}, \scite{cus.21}.
\enddefinition
\bn
Note that we have considered $\varphi(\bar x,\bar y,\bar z)$ as the quadruple
$\langle \varphi,\bar x,\bar y,\bar z \rangle$ with $\bar x \char 94 \bar y
\char 94 \bar z$ a sequence with no repetitions of variables, including every
variable which occurs freely in $\varphi$. \nl
On the relationships of those properties, the independence
property and the strict order property see \cite{Sh:715}.
%\enddefinition
\bn
\ub{\stag{Mon.10C} Problem}:  Is there a reasonable theory for the family of
$(k,m)-*$-NIP first order theories (complete) $T$?  Or for the family of
first order $T$ without the $\omega$-independence property?  Certainly
this is hopeful. \nl
A ``theory" here means say as in \cite{Sh:c} for the class of superstable
(complete first order theories) $T$.
\bn
\ub{\stag{Mon.10D} Question}:  1) Prove that for any $k < \omega$, 
for some $\ell,m$ (in fact, quite low) we hope that any complete 
first order $T$ we have:
$T$ is $(\ell,m)-*$-non-independence, \ub{iff}
$\lambda \rightarrow_{T,k}(\mu)_\chi$ under reasonable conditions on 
$\lambda,\mu,\chi,|T|$ as in the Erd\"os Rado Theorem (rather than large
cardinals). You may use a set of $(\ell,m)$'s.   
\nl
2) For $k = \omega$ we similarly consider the failure of the 
$\omega$-independence property.  This will prove that the
$\omega$-independence property is a real dividing line for \scite{Mon.7},
but I have no reasonable speculations on what a theory for this property will
say.
\bn
What we can get (see \cite{Sh:197})
\proclaim{\stag{Mon.11} Theorem}  If every monadic expansion of $T$ does not
have the independence property, \ub{then}

$$
\beth_\omega(\kappa + |T|)^+ \rightarrow_T (\kappa^+)_\kappa
$$
\endproclaim
\mn
(the property in the assumption is very strong, \ub{but} it is reasonable in
context of ``why the research on monadic logic concentrates on 
trees + linear orders"?  How is this proved?  We 
can decompose any model to a tree sum starting by \scite{Mon.11} with a
large sequence of indiscernible, extend it to a decomposition, so the tree
has 2 levels.  However, the cardinality of the ``leaves" have no apriori upper
bound.  But as there are many leaves such that the model is their sum we 
can show that the model, if it is not too little can be extended to all 
larger cardinalities retaining its monadic theory.

This proves that the dividing line (mentioned in \scite{Mon.11} and
discussed earlier) is real.  

Macintyre had said that cardinals appearing in a theorem make it uninteresting
(though he has moderate lately).  I think inversely and find 
fascinating theorems
showing that for the family of models of $T$ of cardinality $\lambda$ having
a property is equivalent to an inside ``syntactical" property of $T$.  Also,
I think it is a good way of discovering a worthwhile property of $T$ which
should be persuasive even for those who unlike myself do not see their beauty.
Macintyre supposedly is even less friendly toward infinitary logics; but
\bn
\ub{\stag{Mon.11a} Thesis}  We use infinitary logic to ``drown the noise";
only from the distance you see the major outlines of the landscape clearly,
so for many purposes; e.g.  
examining the $L_{\infty,\kappa}$-theory of models of
$T$ will give a more coherent and interesting picture whereas probably
$L_{\infty,\aleph_0}$-theory gives an opaque one.
It probably is not accidental that superstability
was discovered looking at behaviour in cardinals like $\beth_\omega$ and not
from investigating countable models.
\bn
We may consider more complicated partition relations
\definition{\stag{Mon.11b} Definition}  1) Let $\lambda \rightarrow_T
(\mu)^n_\kappa$ for first order complete $T$ means; letting ${\frak C} = 
{\frak C}_T$, the monster model for $T$:
\mr
\item "{$\boxtimes$}"  if $\bar{\bold a}_u \in {}^{\kappa >}{\frak C}$ for
$u \in [\lambda]^{\le n}$, \ub{then} we can find $W \in [\lambda]^\mu$ and
$\bar{\bold b}_u \in {}^{\kappa >}{\frak C}$ for $u \in [W]^{\le n}$ such
that:
{\roster
\itemitem{ $(\alpha)$ }  $\bar{\bold a}_u$ is an initial segment of
$\bar{\bold b}_u$
\sn
\itemitem{ $(\beta)$ }  $v \subseteq u \in [W]^{\le n} \Rightarrow
\text{ Rang}(\bar{\bold b}_v) \subseteq \text{ Rang}(\bar{\bold b}_u)$
\sn
\itemitem{ $(\gamma)$ }  if $m < \omega$ and $i_0 < \ldots < i_{m-1}$ and
$j_0 < \ldots < j_{m-1}$ from $W$ then
there is a ${\frak C}$-elementary mapping
$f$ (even an automorphism of ${\frak C}$) such that: \nl
$v \in [\{0,\dotsc,m-1\}]^{\le n} \and u_1 = \{i_\ell:\ell \in v\}
\and u_2 = \{j_\ell:\ell \in v\}$ implies that 
$f$ maps $\bar{\bold b}_{u_1}$ to $\bar{\bold b}_{u_2}$.
\endroster}
\ermn
2) Similarly for $\rightarrow_K,\rightarrow_M$.\enddefinition

In \cite[Ch.VI]{Sh:h} we get that for stable $T$, no large cardinal is
needed: cardinal bounds which is essentially $(\beth_n(\mu^-))^+$ suffice
(in fact this is done in a much more general framework, and also for trees
$({}^{\omega >}\lambda$).
\bn
\ub{\stag{Mon.11c} Problem}:  For first order theories $T$ for every
$\mu,\kappa$ how large is Min$\{\lambda:\lambda \rightarrow_T
(\mu)^n_\kappa\}$?

We expect a dichotomy: either suitable large cardinal are needed, so
$\beth_k(\mu + \kappa + |T|)$ for $k = k_n$ large enough suffice.
\bn
\centerline{$ * \qquad * \qquad *$}
\bn
Returning to classifying first order theories $T$ by the monadic 
logic, the case of $T$ stable is reasonably analyzed 
(\cite{BlSh:156}, \cite{Sh:284c}), still there is a troublesome dividing line.
\bn
\ub{\stag{Mon.12} Problem}  Assume any model of $T$ is a non-forking sum
of $\langle M_\eta:\eta \in \bold T \rangle$ where 
$\bold T \subseteq {}^{\omega >}
\lambda$ (closed under initial segments).  In some cases the 
${\Cal L}$(mon)-theory is essentially exactly as complicated as that of 
$({}^{\omega >}\lambda,\triangleleft)$, in other cases we can interpret 
$Q_{pr}$.  Can we prove the 
dichotomy, i.e. that always at least one of those holds.
\bn
Probably not so characteristic of me, but I asked \nl
\ub{\stag{Mon.13} Question}:  Is the monadic-Borel theory of the real line
decidable?  \nl
2) Is the monadic theory of $({}^{\omega \ge}2,\triangleleft)$ undecidable?

The meaning of monadic-Borel is that we interpret the monadic quantifier
$(\exists X)\varphi$ by ``there is a Borel set $X$" such that $\varphi$.

The choice of Borel is just a family of subsets of $\Bbb R$ 
(or ${}^{\omega \ge}2)$ which is closed under reasonable operations and 
do not contain subsets
gotten by diagonalization on the continuum.  So ${\Cal P}(\Bbb R) \cap
L[\Bbb R]$ assuming AD is okay, too.  If we try the $({}^{\omega \ge}2,
\triangleleft)$ version, Borel determinacy + Rabin machines looks the 
obvious choice for trying to prove a decidability answer.  
For $(\Bbb R,<)$ it is reasonable
to try to get elimination of quantifiers, i.e. an appropriate version of
$UTh^n(\Bbb R,\bar Q)$ should be enough (\cite[\S4]{Sh:42}).
\newpage

\head {\S3 Automorphisms and quantifiers} \endhead  \resetall \sectno=3
 % \resetall 
\bn
\ub{\stag{aqt.0} Discussion}:
As known for long: for first order complete theory $T$ there are lots of
models with lots of automorphisms (in the direction of saturated ones or
EM ones).  To
build models with no nontrivial ones is hard (even in special cases - there
is literature).  Ehrenfeucht conjectures that the classes

$$
\{\lambda:\psi \text{ has a rigid model in } \lambda,\lambda >
\aleph_0\}
$$
\mn
are simple (like omitting types, in particular: initial segments);
``unfortunately", essentially any $\sum^1_2$ class of cardinals may occur
(see \cite{Sh:56}).
So set theoretically we understand what these families of cardinals are,
but model theoretically the answer is considered negative.  
We may try to change the question, so that we can say something interesting.
\bigskip

\definition{\stag{aqt.1} Definition}  Let $\psi = \psi(\bar R)$ 
be a first order
sentence on the finite sequence $\bar R$ of predicates and function symbols
(with $R_0$, i.e. $\{\bar x_0:R_0(\bar x_0)\}$ being
 ``the universe", so unambiguous and for simplicity each 
$R_\ell$ a predicate; in general $\bar x^0$ is not a singleton, and we may let
$R_1$ be equality).  Consider enriching first order logic by quantifiers 
$Q^\psi = Q^{\text{aut}}_\psi$ 
which means that we can apply $(Q^\psi_{\bar \varphi}f)$ 
to a formula where $\bar \varphi = 
\langle \varphi_\ell (\bar x^\ell,\bar z):\ell < \ell g(\bar R) \rangle,
\ell g(\bar x^\ell) =$ arity of $R_\ell$, and in the inductive definition
of satisfaction $M \models (Q^\psi_{\bar \varphi,\bar a}f)\vartheta$ holds
when: if $\langle \varphi_i(\bar x^\ell,\bar a):\ell < \ell g(\bar R) 
\rangle$ defines in $M$ an $\bar R$-model 
$M_{\bar \varphi,\bar a}$ of $\psi$ then there is an automorphism $f$ of 
$M_{\bar \varphi,\bar a}$ such that $\vartheta$ holds.  So syntactically $f$ 
is a variable on partial unary functions.
\mn
Note: those quantifiers (\cite{Sh:43},\cite{Sh:e}; really more general
there, see \scite{aqt.16}) do not exactly fit
``Lindstrom quantifiers".  They can be expressed artificially by having
many Lindstrom quantifiers and 
each Lindstrom quantifier is a case of this.  But those
are naturally second order quantifiers and e.g.
adding two such quantifiers is more than adding the
cases for each.  So for a vocabulary $\tau$ in the language

$$
L_{\omega,\omega}(Q^\psi)(\tau)
$$
\mn
we have variables: individual variables and unary partial function variables, 
we can form $(Q^\psi_{\langle \psi_\ell:\ell < \ell g(\bar R) \rangle 
{\bar y}})\vartheta$ if 
$\varphi_\ell,\vartheta$ are already in 
$L_{\omega,\omega}(Q^\psi)(\tau)$ and
satisfaction is defined as above.
We may allow such quantifier to act only on models $M_{\bar \varphi,\bar a}$
whose universe is $\subseteq M$ \ub{or} to allow the set of elements of 
$M_{\bar \varphi,\bar a}$ (equivalently $\bar x^0$) to be e.g. the set
automorphisms of $M_{\varphi^*,{\bar b}^*}$ for any $\bar b^*$ satisfying
say $\theta^*(\bar y,c)$ where $\psi^*,\vartheta^*$ are formulas in our logic
of smaller depth, etc.
For compactness this does not matter.
\enddefinition
\bn
\ub{\stag{aqt.2} Problem}:  For which $\psi$ is $L(Q^\psi)$ a compact logic?
\bn
\ub{\stag{aqt.3} Example}:  If $\bar R = \langle R_0 \rangle,\psi = \forall x
R_0(x)$, then we have quantifications on unary functions varying on
permutations, so the quantifier $Q^\psi$ gives second order logic (on
nontrivial structures).  So in this problem even 
though the models of $\psi$ can be written as $M_0 +M_1$ or
``degenerated", we get second order logic.
\bn
\ub{Note}:  So for this classification a sentence $\psi$ 
which says ``the model (of $\psi$) is trivial" gives a complicated logic
$L_{\omega,\omega}(Q^\psi)$. 

If $\psi$ has only finite models, the logic is compact in a dull way.
You may wonder if compactness 
holds for any sentence $\psi$ at all, as this looks like a
second order logic.
However, there are interesting sentences $\psi$ with $L(Q^\psi)$ compact:
\mr
\item "{$(a)$}" $\psi =$ the axiom of the theory of Boolean Algebras, (i.e.
conjunction of the axioms in standard axiomatization)
\sn
\item "{$(b)$}"  the axiom of the theory of ordered fields.
\ermn
We expect that if the models of $\psi$ are complicated enough, the logic
will be compact.
We may also have applications to the compactness: it was known

CON(there is $1$-homogenous\footnote{A Boolean Algebra $B$ is 
1-homogeneous if for any $x,y \in B \backslash \{0_B,1_B\}$, some 
automorphism of $B$ map $x$ to $y$}
atomless Boolean Algebra $B$ \nl

$\qquad \quad$ such that Aut$(B)$ is not simple) 
\mn
\ub{even}: Con$(\exists G \triangleleft \text{ Aut}(B)$(Aut$(B)/G
\text{ commutative}))$
\mn
(see \cite{Sh:384}; it was known that 
Aut$(B)'$ (= commutator subgroup) is simple).
\mn
So the compactness and completeness theorems show: ZFC $\vdash$ ``there are
such Boolean Algebras".  So considering the success of the compactness and
completeness theorems having such quantifier will be plausably in addition
to being good by itself, also applicable.
\bn
So we are interested in: \nl
\ub{\stag{aqt.4} Problem}:  Find more such quantifiers (homomorphisms of
embeddings instead of automorphisms are welcomed, see \cite{Sh:e} on the
cases above).
\bn
The proof gives more examples but we like to have: \nl
\ub{\stag{aqt.4a} Problem}:  Characterize the $\psi$ for which we have
a compact $L(Q^\psi)$ or at least find:
\mr
\item "{$(a)$}"  general criterion
\sn
\item "{$(b)$}"  natural examples rather than those which look to the 
proofs one has.
\endroster
\bn
We may consider also: \nl
\ub{\stag{aqt.5} Problem}:  Characterize the strongly rigid first order theory
$T$ and the rigid ones where, 
\definition{\stag{aqt.5a} Definition}  1) First order 
$T$ is called strongly rigid if:
for every theory $T_1 \supseteq T$ there is a theory $T_2 \supseteq T_1$ such
that the pair $(T_2,T)$ is rigid which means that $T_2$ has a model $M_2$, 
such that every $f \in \text{ Aut}(M_2 \restriction \tau(T))$ is first 
order definable with parameters in $M_2$.  
We say $T$ is super rigid if above $T_2 = T_1$.
We say $T$ is essentially rigid if $(T,T)$ is rigid.  
We say $(T_1,T)$ is rigid for
$\varphi(M)$ if $\bold \varphi(-)$ is a property of models of $T_1$ and 
$M_2 \restriction \tau(T)$ satisfies $\varphi$ (e.g. $|T_1|^+$-saturated).
We add ``in $\lambda$" if the models is required to be of cardinality
$\lambda$.  \nl
2) We add the adjective everywhere if we omit the demand ``$T \subseteq T_1$"
and replace $f \in \text{ Aut}(M_2 \restriction \tau_T)$ by $f$ an
automorphism $M$ of $T$ which is interpreted in $M_2$ by (first order)
formulas with parameters (as in \scite{aqt.1}, of course the model $M'$ of
$T$ has the vocabulary of $T$). \enddefinition

This is a way to classify $T$'s.

Those are relatives of having rigid models.  The definable automorphisms 
are the parallel of inner automorphisms of a group. 
Note that all those notions do not imply that $T$ has a rigid model; if $M$
is a complete, say infinite, non-abelian group
(i.e. any automorphism is inner) then
Th$(M)$ is essentially rigid but has no rigid model.  The version with
$T_1 = T_2$ is the best case.  
If we replace $T$ by any model of $T$ interpretable in $M_1$ (as in
\scite{aqt.5a}(2)) and allow $T_1$
to have parallels of Skolem functions we are approximating the compactness
and completeness problem discuss above.  We may even let 
$T_1 = Th({\Cal H}(\chi),\in),\chi$
strong limit, and consider interpretation of $T$ on ``sets" of the model
$M_2$ of $T_2$ rather than classes. \nl
Why have we concentrated on ordered fields and Boolean Algebras? \nl
The point is that e.g. for a dense partial order we can get a model where for
every partial order definable in it, every automorphism of it as a partial 
order is, for a dense set of intervals, definable with parameters. 
(If the partial order is not dense, consider ``infinite intervals").  Why
``ordered field"?  Only as in this case there any automorphism is determined 
by its action
on any interval.  Concerning Boolean Algebras, the underpinning point is that
we consider structures $(A,B,R),R \subseteq A \times B$ which satisfy
comprehension, that is:

$$
(\forall y_1 \ne y_2 \in B)(\exists x \in A)(xRy_1 \equiv xRy_2)
$$
\mn
and have the strong independence, that is,

$$
\align
(\forall x_1,\dotsc,&x_n \in A)(\forall y_1,\dotsc,y_n \in A)(\exists z) \\
  &(\dsize \bigwedge_{\ell,k} x_\ell \ne y_k \Rightarrow 
(\dsize \bigwedge^n_{i=1} x_i Rz \and \neg y_iRz)).
\endalign
$$
\mn
(An abvious exampls is an atomic Boolean Algebra 
$B,A =$ atoms$(B),B=B$ are okay). \nl
For some of the readers a bell may ring.  A theory $T$ is 
\ub{unstable}:  iff it has the strict order property (that is some 
$\varphi(\bar x,\bar y)$ is a partial order with infinite chains) \ub{or} 
has the
independence property (a relative of the strong independence property).  This
does not say any unstable theory will do but indicates that an unstable
theory at least locally will do.
\mn
\ub{Note}:  For the theory of linear orders, for $(A,<)$, if $E$ is a convex
equivalence relation with classes $A_i$ for $i < i^*$ and 
$f_i \in \text{ Aut}(A,<)$ maps 
$A_i$ to itself, then $\dsize \bigcup_i (f_i \restriction A_i) \in \text{ Aut}
(A,<)$.  (We can express that informally as ``models of $T$ are, in general,
decomposable"; to avoid trivialities we restrict ourselves to 
uncountables ones).
So for $T$ any theory of infinite linear orders, $T$ is not strongly rigid.  
We need $\psi$ (or $T$) to say that the model is not decomposable.
\bn
Generally,
\definition{\stag{aqt.6} Definition}  1) We say $\psi$ (or $T$) is pseudo
decomposable when:  if for every $n$, there are a model $M$ of 
$\psi$ (or of $T$), $M$ the disjoint union of the nonempty sets 
$A_i$ (for $i < n$) and $f^1_i \ne f^2_i$ from $\text{Aut}(M)$ such that

$$
f^1_i \restriction (M \backslash A_i) = f^2_i \restriction (M \backslash
A_i)
$$
\mn
and $\dbcu^{n-1}_{i=0} (f^{\eta(i)}_i \restriction A_i) \in
\text{ Aut}(M)$ for every $\eta \in {}^n 2$; in other words, $M$ has a
nontrivial automorphism over $M \backslash A_i$ for each $i$. \nl
2) We say $\psi$ (or $T$) is semi-decomposable if for every $n$ we can find a
model $M$ of $\psi$ and partition $\langle A_\ell:\ell < n \rangle$ of $M$
to infinite subsets such that:
\mr
\item "{$(*)$}"  for every finite set $\Delta_1$ of formulas in 
$L(\tau_T)$ there
is a finite set $\Delta_2$ of formulas in $L(T)$ such that
\sn
\item "{$(**)$}"  if for $\ell < n,k^\ell < \omega$ and $\bar a_\ell,
\bar b_\ell \in {}^{k_\ell}(A_\ell)$ and $\bar a_\ell,\bar b_\ell$ realize
the same $\Delta_2$-type in $M$ and $\ell = 0 \Rightarrow \bar a_\ell =
\bar b_\ell$, \ub{then} $\bar a_0 \char 94 \bar a_1 \char 94 \ldots 
\bar a_{n-1},\bar b_0 \char 94 \bar b_1 \char 94 \ldots \bar b_{n-1}$
realize the same $\Delta_1$-type in $M$. 
\ermn
3) We say almost decomposable if the function $\Delta_1 \mapsto \Delta_2$
does not depend on $n$.
\enddefinition
\bigskip

\proclaim{\stag{aqt.7} Claim}  1) If $T$ is pseudo decomposable, \ub{then}
we can find $T_1 \supseteq T$ such that:
\mr
\item "{$(a)$}"  for any model $M_1$ of $T_1$ we have: $(\text{Aut}(M_1
\restriction \tau(T))$ has cardinality $2^{\|M_1\|}$ hence some $f \in
\text{ Aut}(M_1 \restriction \tau(T))$ 
is not definable in $M_1$ even with parameters
\sn
\item "{$(b)$}"  if $T = \{\psi\}$ then for models of $T_1$, in the logic
$L(Q^\psi)$ we can interpret second order logic on $M_1$
\sn
\item "{$(c)$}"  we can embed also some product $\Pi\{G_i:i < \|M_1\|\},G_i$ a
nontrivial group.
\ermn
2) If $T$ is semi-decomposable then $T$ is pseudo decomposable. \nl
3) If $T$ is almost decomposable, \ub{then} it is semi-decomposable and
for any
saturated model of cardinality $\lambda$ (or just $\lambda^+$-resplended
model of $T$), $\lambda > |T|$, we can find $\langle A_i:i < \lambda \rangle$
as in \scite{aqt.7}, in fact:
\mr
\item "{$(a)$}"  if $\bar a_i,\bar b_i \in {}^{\omega >}(A_i)$ for
$i < \lambda,\bar a_0 = b_0$, tp$(\bar a_i,\emptyset,M) = \text{ tp}(\bar b_i,
\emptyset,M)$ then $\bar a_{i_0} \char 94 \ldots \char 94 \bar a_{i_n},
\bar b_{i_0} \char 94 \ldots \char 94 \bar b_{i_n}$ realizes the same type
in $M$ for any $i_0 < \ldots < i_n < \lambda$
\sn
\item "{$(b)$}"   $(M,A_i)_{i \in \omega}$ is $\lambda$-saturated for
$\omega \in [\lambda]^{< \lambda}$.
\endroster
\endproclaim
\bn 
For Boolean Algebras we can decompose the set of atoms, but the image of an
element is not deciphered so this theory is not even pseudo decomposable.
\sn
Be careful, the statement ``$B$ is the Boolean Algebra generated by the
close-open intervals of a linear order $I$" is not first order (this follows
by the compactness so if $T_1$ extend the 
theory of Boolean algebras then it has
a model with no undefinable automorphism).  Now for the first problem,
\scite{aqt.2}, the hope is that failure pseudo indecomposability is enough
for compactness, it is of course necessarily by \scite{aqt.7}.
\bn
\ub{\stag{aqt.8} Question(Cherlin)}:  What occurs for
vector spaces over finite fields?

Let $\Bbb F$ be a (fixed) finite field and let $\psi_{\Bbb F}$ be 
the conjunction of the
axioms of vector spaces over the field $\Bbb F$ 
(we have binary function symbols for
$x+y,x-y$, individual constant 0 and unary functions $F_c$ for
$c \in \Bbb F$ to denote multiplication by $c$).  There is $T_1 \supseteq
\{\psi_{\Bbb F}\}$ such that for models of $T_1$, in the logic
$L(Q^{\text{aut}}_{\psi_{\Bbb F}})$ we can interpret second order logic
on $M_1$ (similarly for a finitely generated field).
\bigskip

\remark{Remark}  Also for a general field this works, except that we 
do not have the
quantifier as $\psi_{\Bbb F}$ is an infinite conjunction of first order
formulas. \nl
Why?  Enough to have $T_1$ such that; for $M^*$ a model of $T_1$:
\mr
\widestnumber\item{$(iii)$}
\item "{$(i)$}"    $P^{M^*}_1,P^{M^*}_2$ are disjoint subsets 
of the vector space
\sn
\item "{$(ii)$}"  $P^{M^*}_1 \cup P^{M^*}_2$ is an independent set in the
vector space
\sn
\item "{$(iii)$}"  $F_1$ is a unary function, $F_1 \restriction
P^{M^*}_\ell$ is one-to-one onto $M^*$ for $\ell=1,2$
\sn
\item "{$(iv)$}"  $F_2$, a two-place function, is a pairing function.
\ermn
How is the interpretation?  For any function $g$ from $P^{M^*}_1$
to $P^{M^*}_2$ there is an automorphism $f$ of the vector space such that:

$$
x \in P^M_1 \Rightarrow f(x) = x + g(x)
$$

$$
x \in P^M_2 \Rightarrow f(x) = x.
$$
\endremark 
\bn
My impression is that any reasonable 
example will fall easily one way or the other by existing methods.  
\bn
\centerline{$* \qquad * \qquad *$}
\bigskip

\definition{\stag{aqt.9} Definition}  Let $M$ be a model of $T,P \subseteq M$.
We say $T' = Th(M,P)$ has the automorphic embeddability property
over $P$ if for every model $(M',P')$ of $T'$, every automorphism $f$ 
of $M' \restriction P'$ can be extended to an automorphism of $M$.
\enddefinition
\bn
\ub{\stag{aqt.10} Question}  Characterize the theories $T' = Th
(M,P)$ which has the automorphic embeddability property over $P$.

This looks hard on us as characterization of this would probably involve
${\Cal P}^-(n)$-diagrams as in classification over a predicate;
on the case with no two cardinal models (i.e. $\|M\| > |P^M|$ assuming there
is $\lambda = \lambda^{< \lambda} \ge |T|$), see \cite{Sh:234}.  The general
case is, unfortunately, still in preparation (\cite{Sh:322})); see 
end of \S6.
\bn
\centerline{$* \qquad * \qquad *$}
\bn
There are other ways to consider quantification over automorphisms: 

For a model $M$ let $(M,\text{Aut}(M))$ be the two sorted model, one sort
is $M$, the other is the group Aut$(M)$, with the application function, that
is in the formulas, we allow forming $f(x)$ for $x$ of first sort, $f$ of
second sort.  We may replace Aut$(M)$ by the semi-group of endomorphisms or
one-to-one endomorphisms.
\bn
Now for a variety ${\Cal V}$, the complicatedness of the first order theory
of the endomorphism semi-group End$(F_\lambda)$ of the free algebra with
$\lambda$ generators is reasonably understood (see \cite{Sh:61}) but so far
not the automorphism group in the general case though 
several specific cases were analyzed, (see \cite{ShTr:605}, \cite{BTV91}).
\bn
\ub{\stag{aqt.11} Problem}:  For which varieties ${\Cal V}$, 
letting $F_\alpha$ be the free algebra in ${\Cal V}$ generated by
$\{x_i:i < \alpha\}$, can we in Aut$(F_\lambda)$ (first order) interpret 
second order theory of $\lambda$?  We hope for a solution which depend 
``lightly" on $\lambda$ (like Aut$(\lambda,=))$?
We may allow quantification on elements or even use $(F_\lambda,\text{Aut}
(F_\lambda))$; but, of course, better is if we succeed to regain it.
\bn
The following property looks like a relevant dividing line
\definition{\stag{aqt.12} Definition}  We say the variety 
${\Cal V}$ is Aut-decomposable if:
\mr
\item "{{}}"  if $F_{\omega 2}$ is the algebra generated freely 
by $\{x_i:i < \omega + \omega\}$ for ${\Cal V}$ and 
$f \in \text{ Aut}(F_{\omega 2})$ satisfies $f(x_n) = x_n$ for
$n < \omega$, \ub{then} we have: \nl
$f$  maps $\langle x_{\omega +n}:n < \omega \rangle_{F_{\omega 2}}$ to
itself.
\ermn
Why?  For varieties ${\Cal V}$ with this property we can repeat the analysis
of Aut$(\lambda,=)$ which is the group of permutations of $\lambda$; though 
first order interpretation of elements has to be
reconsidered.  But this is not needed in generalizing the ``upper bound",
the equivalences.  That is for proving, say ${\Cal V}$ with
countable vocabulary for simplicity, that Th(Aut$(F_{\aleph_\alpha}))$ depend
``lightly" on $\alpha$; i.e. if for $\ell = 1,2,\alpha_\ell = \delta_\ell +
\gamma_\ell,\gamma_\ell < ((2^{\aleph_0})^+)^\omega$ (ordinal exponentiations)
and Min$\{\text{cf}(\alpha),(2^{\aleph_0})^+\} = \text{ Min}\{\text{cf}
(\alpha_2),(2^{\aleph_0})^+\},\gamma_1 = \gamma_2$ then 
Th(Aut$(F_{\aleph_1})) = \text{ Th(Aut}(F_{\aleph_{\alpha_2}}))$.  
On the other hand, if it fails an automorphism of
$F_\lambda$ code a complicated subset of $\lambda$. 
\enddefinition
\bn
\centerline {$* \qquad * \qquad *$}
\bn
We may look at questions on ``are there logics with specified properties?"
\nl
An old problem (see \cite{BF}): \nl
\ub{\stag{aqt.13} Question}:  Is there an 
$\aleph_1$-compact extension of $L(Q)$ which has interpolation (Craig)? 
\bn
I prefer \nl
\ub{\stag{aqt.14} Question}:  Is there in addition to first order logic
a compact logic which has interpolation?

Barwise prefers to look at definability properties of logics (e.g.
characterizing $L_{\infty,\omega}$) but my taste goes to:
\bn
\ub{\stag{aqt.15} Problem}:  Find 
(nontrivial) implications between properties of logics. \nl
See for example \cite{Mw85}, \cite{Sh:199};  interpolation and 
Beth theorems are, under reasonable assumption, equivalent; 
and amalgamation essentially implies
compactness.  After great popularity in the seventies, the interest has gone
down, a contributing factor may have been the impression that there are mainly
counterexamples.  This seems to me too early to despair. \nl
However, Vannanen's book \cite{Va9x} should appear.
\bn
\ub{\stag{aqt.15a} Discussion}:
So we are interested in enriching first order logic by additional quantifiers
preserving compactness and getting interpolation. \nl
A natural play is to allow second order variables $X$ but restrict the
existential quantifier to cases when the relation $P$ (or function) satisfies
some first order sentence $\psi$ with some specific old $R$ as a
parameter (e.g. $P$ is $f$, an automorphism of a model of $\psi$, a case
discussed above).  Another way (\cite{Sh:43}) is to replace ``exist" by 
``the family of those satisfying it belongs to a family $\Bbb D(M)$ of 
such relations over $M$".  An example introduced in \cite{Sh:43} is the
case of a unary predicate, with $\Bbb D(M)$ being the club filter on 
$[|M|]^{\aleph_0}$; or equivalently for the strength of the logic, the 
family of stationary subsets of $[|M|]^{\aleph_0}$.  Those quantifiers are 
$(aaP),(stP)$, respectively.  This logic 
has many properties like $L(Q)$, see \cite{BKM78}, some like 
second order, \cite{ShKf:150}, \cite{Sh:199}. \nl
Now interpolation holds for the pair of logic 
\footnote{$Q^{\text{cf}}_{\le \aleph_0}$ tells us the cofinality of a
linear order is $\le \aleph_0$}
$(L(Q^{\text{cf}}_{\le \aleph_0}),L(aa))$ which means: if $\varphi_\ell$ is
a sentence in $L(Q^{cf}_{\le \aleph_0})(\tau_\ell)$ for $\ell = 1,2$ and
$\vdash \varphi_1 \rightarrow \varphi_2$ then for some $\psi \in L(aa)
(\tau_1 \cap \tau_2)$ we have $\vdash \varphi_1 \rightarrow \psi$ and
$\vdash \psi \rightarrow \varphi_2$.  Also the Beth closure of
$L(Q^{\text{cf}}_{\le 2^{\aleph_0}})$ is compact so there is a compact logic
which satisfies the ``implicit definability implies the explicit
definability"; moreover, is reasonably natural (at least in my eyes).  
Seems near the mark but not in it.  Consider (see \cite{HoSh:271}) the 
following logic: let $\aleph_0 < \kappa < 
\lambda$ and $\kappa,\lambda$ are compact cardinals, and expand first order
logic by all the connectives of the form
$\dsize \bigwedge_D \varphi_i$ where $D$ is a 
$\kappa$-complete ultrafilter on some $\theta \in [\kappa,\lambda)$, 
meaning naturally $\{i < \theta:\varphi_i\} \in 
D$.  It has interpolation but not full compactness (only $\mu$-compact for
$\mu < \kappa$). 
\bn
\centerline {$* \qquad * \qquad *$}
\bn
More formally and fully
\definition{\stag{aqt.16} Definition}  1) Assume
\mr
\item "{$(a)$}"  $\psi(\bar R,S)$ be a sentence (usually in first order) in
the vocabulary $\bar R$, i.e. $\bar R$ is a list with no repetitions of
predicates and function symbols, $S$ an additional predicate, each have a
given arity (for notational simplicity $R_0$ is a unary predicate for
``the universe" each $R_\ell$ is a predicate)
\sn
\item "{$(b)$}"  $\Bbb D$ is a function, its domain is $K_{\bar R} =
\{M:M \text{ a model, the vocabulary of } M$ \nl
$\text{is } \{R_i:1 + i < \ell g(\bar R)\},R^M_0 
= |M|\},\Bbb D(M)$ is a family of subsets of \nl
$\{N:N$ is an expansion of $M$ by $S^N\}$ and, of course, if $f$ is an
isomorphism from $M_1$ onto $M_2$ then $f$ maps $\Bbb D(M_1)$ onto
$\Bbb D(M_2)$.
\ermn
The quantifier $Q_{\bar \varphi,\Bbb D}$, syntactically acts as
$(Q_{\bar \varphi,\Bbb D}S) \vartheta$ where $S$ is a variable on
$n$-place relations, $\bar \varphi = \langle \varphi_i(\bar x_i,\bar z_i):
i < \ell g(\bar R) \rangle$ and $\ell g(\bar x_i) = \text{ arity}(R_i)
\rangle$ and $M_{\bar \varphi,\bar a} = M_{\bar \varphi,\bar a}[{\frak A}]$
is defined as in Definition \scite{aqt.1} and 

$$
\align
{\frak A} \models &(Q_{\bar \varphi,\Bbb D}S)\vartheta
(S,\bar a) \text{ iff} \\
  &\{S:S \text{ an } n\text{-place relation on } \{x:\varphi_0(x,\bar z)\} \\
  &\text{ and } (M_{\bar \varphi,\bar a}[{\frak A}],S) \models \vartheta
(S,\bar a)\} \\
  &\text{ belongs to } \Bbb D(M_{\bar \varphi,\bar a}[{\frak A}]).
\endalign
$$
\mn
2) Similarly when for defining $M_{\bar \varphi,\bar a}[{\frak A}]$ we
replace equality by an equivalence relation $R_1$.
\enddefinition
\bn
A variant of \scite{aqt.14} is \nl
\ub{\stag{aqt.17} Question}:  Is there a reasonably defined such quantifier
$Q$ such that $L(Q)$ is compact and has interpolation? or at least has the
Beth property?
\newpage

\head {\S4 Relatives of the main gap} \endhead  \resetall \sectno=4
 % \resetall 
\bn

A main gap theorem here means, for a family of classes of models, that for
each class $K$ \ub{either} we have a complete set of invariants for models
of $K$ (presently, which are basically just sets) \ub{or} it has quite 
complicated models
(see below after \scite{cft.8}).

This seems obviously a worthwhile dichotomy, if it occurs indeed, and have
been approached as a dichotomy on the number of models (but see below).  

We know for a countable first order, $T$ complete for simplicity,  that
$I(\aleph_\alpha,T) =: \{M/\cong :M \models T \and \|M\| = \aleph_\alpha\}$ 
behave nicely (either $I(\aleph_\alpha,T) = 2^{\aleph_\alpha}$ for
every $\alpha > 0$ or $< \beth_{\omega_1}(|\alpha|+\aleph_0)$ for every 
$\alpha$).
But many relatives of this question are open.
\bn
I thought a priori on several of them that they will be easier, but have
worked more on the case of models so the earlier solution in \cite{Sh:c} does
not prove this thought wrong; still this a priori opinion is not necessarily
true. \nl
\ub{\stag{cft.0} Question}  1) Prove the main gap for the class of 
$\aleph_1$-saturated models. \nl
2) Prove the main gap for the class of $\aleph_0$-saturated models. 

Now \scite{cft.0}(1) have looked a priori relatively not hard, 
in fact the work in \cite{Sh:c} seems
to solve it ``except" for 
lack of regular types, so in the decomposition theorem
we are lacking how to exhaust the model.  
\bn
Another direction is: \nl
\ub{\stag{cft.1} Problem}:  Let $T$ be countable stable complete first order
theory.  Show that if $\neg(*)$, \ub{then} $T$ has otop (or dop; for otp we
allow types over countable sets), where
\mr
\item "{$(*)$}"  if $M_0 \prec M_\ell \prec {\frak C}$ for 
$\ell = 1,2$ and \nl
$M_1 \nonfork{}{}_{M_0} M_2$, (i.e. tp$_*(M_2,M_2)$ does not fork over $M_0$)
\ub{then} there is a prime model (even $F_{\aleph_0}$-prime)
over $M_1 \cup M_2$.
\ermn
\ub{Note}:  For superstable this is true 
(this was the main last piece for the main gap, see \cite[Ch.XII]{Sh:c}).
\bn
\ub{\stag{cft.1a} Discussion}: 
Our problem is that the proof there uses induction on ranks, and generally
stable theories have less well understood theory of types (not enough regular
types exist), just as in \scite{cft.0}.  However, if we assume $T$ superstable
without DOP, then every
regular type is either trivial (= the dependence relation
is) or of depth zero (\cite[Ch.X,\S7]{Sh:c}).  
There is some parallel theorem for stable theories without DOP, it 
may be helpful.
\sn
Maybe relevant is the theory of types for stable $T$ in \cite[Ch.V,\S5]{Sh:c},
\cite{Sh:429} and Hernandez \cite{He92} which proved that if 
$\bold I_0,\bold I_1$ are 
indiscernible not orthogonal then for some indiscernible $\bold J,\bold I_\ell
\le_s \bold J$ (\cite[Ch.V,\S1]{Sh:c}), but in spite of early expectation
this has not been enough to solve \scite{cft.0}(1). \nl
Where could \scite{cft.1} help?  For the theories which are ``low" 
for the main gap, a model is characterized up to isomorphism by its 
$L_{\infty,\aleph_1}$(dimensional quantifier) theory.  
But we may look at logics allowing e.g. a sequence of quantifiers with 
countable length (even $\omega_1$), as investigated by the Finnish school.  
We know that for unstable theories, 
and for stable theories with DOP we have the nonstructure, see
Hyttinen Shelah \cite{HySh:676}.  It seems that \scite{cft.1} would 
complete a piece in finding another dividing line here.  
Some stable, unsuperstable
theories become low.  Essentially, the hope is that 
\ub{either} every model of $T$ can 
be coded by trees with
at most $\delta^*$ levels, $\delta$ fixed, even countable or $\le \omega_1$
\ub{or} we have the order property (even independence property) in a stronger
logic (in NOPOT or NDOP holds).  However, \scite{cft.1} is not enough, we
need also a decomposition theorem.
\bn
\ub{\stag{cft.1b} Question}:  If $T$ is countable stable with NDOP and NOPOT
and $(*)$ of \scite{cft.1} holds, does the decomposition theorem hold at
least for shallow $T$?  \nl
Interpretation of 
groups may be relevant, particularly non-isolated types, because
non-orthogonal, weakly orthogonal types tend to involve groups. \nl
Note that here the existnece of $\bold F^\ell_{\aleph_0}$-primary model on
$N \cup \{a\}$ included inside a given 
$M \supseteq N \cup \{a\}$ is not assured.
\bn
Another problem is \nl
\ub{\stag{cft.2} Problem}:  Prove the main gap for $K_T =: \cap\{M_1
\restriction T:M_1 \text{ a model of } T_1,T \subseteq T_1,|T_1| \le
2^{\aleph_0}\}$. \nl
Note that if $(\forall \lambda < 2^{\aleph_0}) 2^\lambda = 2^{\aleph_0}$,
then we can find one $T_1$ which suffices (as by Robinson lemma we have
``amalgamation" for theories, so there is a universal (oven ``saturated") 
$T_1$, i.e. if $T'_1 \supseteq T$ is complete, $|T'_1| \le 2^{\aleph_0}$ 
by changing names of predicates not in $\tau(T)$ we can embed $T'_1$ 
into $T_2$ over $T$.
\nl
Like all these problems, possibly a large part of the work is 
already done, but though a priori I thought this was easier, it is not 
necessarily true.  The natural hidden order property is by 
$\exists^{\ge \lambda} x \varphi(x,y,z)$ (cardinality quantifiers) (maybe on
the number of equivalence classes or dimension for $\Delta$-indiscernible
sets, $\Delta$ finite),  we hope there will not be a need to consider several
cardinality quantifiers simultaneously.  
If $M$ is a model of $T$ which looks like
$(A,0,P,E),A = \omega \cup \{(n,m,\ell):\ell < k(n,m)\},P = \omega,F_1,F_2$
unary functions, $F_\ell(n_1,n_2,k) = n_\ell$ for $\ell=1,2,
E = \{((n_1,m_1,k_1),(n_2,m_1,k_2)):n_1 = n_2,m_1 = m_2$ and $k_1,k_2 <
k(n_1,m_1)\}$ and the function $k(n,m)$ random enough, $T$ has a hidden
order property, that is, the formula $\varphi(x,y) =: (\exists z)(F_1(z) = x
\and F_2(z) = y \and (\exists^{\ge \aleph_1}z')(z'Ez))$. \nl
We phrase it appropriately (and there are fewer divisions).
\sn
The very low parts of the hierarchy have been analyzed, i.e. the bottom part:
categorical

$$
\aleph_\alpha > 2^{\aleph_0} \Rightarrow I(\aleph_\alpha,K_T) =1 \text{ or }
I(\aleph_\alpha,K_T) \ge 2^{|\alpha|}.
$$
\mn
For the main gap, we can assume $T$ is superstable and we 
should analyze for $M \in K_T$, which we know is
$\aleph_\varepsilon$-saturated and it is natural to analyze the different
dimensions.
\mn
\ub{Note}:   If $T$ is a theory of one equivalence relation $E$ 
saying there is an
equivalence class with $n$ elements, for infinitely many $n$, it is 
not in the lowest
class, but still we understand it.  For the $T$ above, if for every $n$ we
have $\aleph_0$ classes with $n$ elements then $K_T = \{M:$ there are $\|M\|$
classes of cardinality $\|M\|$, for each $x \in M$, has $\|M\|$-classes with
$(x/E)$-elements and $\forall n \exists x(|X/E|=n)\}$. \nl
In the first case (i.e. $I(\aleph_\alpha,K_T) =1$) every model is saturated.
We expect that if $I(\aleph_\alpha,K_T) = 2^\alpha$, then for $M \in K_T$,
there is an equivalence relation between indiscernible sets but on the
set of equivalence classes, there is no further structure (well, maybe unary
``predicates") such that no two equivalence classes have the same dimension.
\mn
In general, the theme taken for granted is: \nl
\ub{\stag{cft.3a} Thesis}:  If $K$ is a reasonable \footnote{but see on
rigidity!} class of models then:
\mr
\item "{{}}"  the behaviour of the function $I(\aleph_\alpha,K)$ is uniform,
the ``same" for all relevant cardinals.
\ermn
Have not really looked.  The expectation is that after this level, we'll
have $2^{2^\alpha}$, (and also $2^{|\alpha|^{\aleph_0}},
2^{|\alpha|^{2^{\aleph_0}}}$, etc), and $\beth_n(|\alpha|)$, (or
$\beth_n(|\alpha|^{\aleph_0}),\beth_n(|\alpha|^{2^{\aleph_0}}))$ and
$\beth_\zeta(|\alpha|)$ for each $\zeta \in [\omega,\omega_1)$ and then
after our $\omega_1$ steps we have
$\aleph_\alpha \mapsto 2^{\aleph_\alpha}$.
\bn
The main \footnote{why I have been feeling so?  As for almost all this book,
countability plays a minor role} question left in \cite{Sh:c} is \nl
\ub{\stag{cft.4} Problem}:  Prove the main gap for uncountable $T$.

The problem in proving is the lack of primary models, particularly over
nonforking triples of models.   Maybe more interpretation of groups
will help in solving this.  
Maybe replacing ``primary models" by prime models, and 
isolated types by unavoidable ones may help.
\sn
(Recall that $B$ is primary over $A$ if $B = A \cup \{a_i:i < \alpha\}$ and
tp$(a_i,A \cup \{a_j:j < i\}$) is isolated for $i < \alpha$.) \nl
Isolated types have been great (for $\aleph_0$-categoricity, no $T$ with
exactly two countable models, Morley
theorem), but for an uncountable theory they are not sufficient, the lack of
them does not witness much.  Still there can be prime models.

Maybe we should look at derived non-elementary classes, where we look for
hidden order and if there is none we get nicer properties.  Maybe even
define such classes inductively on $\alpha < \omega_1$ (or even
$D(x=x,L,\infty)$, but carry enough connection to the original $T$
to be able to finish soon (and carry enough to continue, see
\cite{Sh:h}, \cite{Sh:600}).

It may be reasonable to start with analyzing unidimensional $T$
(concerning \scite{cft.4}).
\bn
\ub{\stag{cft.5} Thesis}:  All such problems have a ``good" solution, (unlike
Ehrenheuft Conjecture, see \cite{Sh:54}, see \S3).
\bn
The audience asked \nl
\ub{\stag{cft.6} Question}:  Can a theory $T$ be ``nice" in 
spite of having many models, maybe still models of $T$ can be 
understood by invariants.
\bn
\ub{\stag{cft.7} Answer}:  ``Nice" certainly yes (see \S5 as you may choose to
consider say linear order as reasonable invariants and so ask for which first
order theories such invariant suffice).  But not true, 
if you define a generalized cardinal
invariant as follows (for simplicity $|T| = \aleph_0$).
\mn
Depth zero: cardinal invariant is a cardinal
\mn
Depth $\alpha+1$: cardinal invariant are sets of sequences of length $\le
2^{\aleph_0}$ of cardinal invariants of depth $\alpha$ \ub{or} a cardinal
invariant of depth $\alpha$
\mn
\ub{depth $\delta$ for $\delta$ limit}:  depth $\alpha$ for some $\alpha <
\delta$
\bigskip

\proclaim{\stag{cft.8} Claim}  If models of $T$ of cardinality $\aleph_\alpha$
are characterized up to isomorphism by generalized cardinal of depth
$\le \gamma_T$, \ub{then} $I(\aleph_\alpha,T) \le \beth_{\gamma_T+1}
(2^{\aleph_0})$ (see \cite{Sh:200}).
\endproclaim
\bn
Really, the main gap for countable complete $T$'s is a division to three
cases.  If $T$ is in the upper case, model of $T$ codes 
stationary sets; if $T$ is in the lower case, a
model can be described by a tree with $\le \omega$ levels and depth
$\le \gamma_T < \omega_1$; and if $T$ is in the middle case, a model can be
described by a tree with $\le \omega$ levels, but can have depth an
arbitrarily large ordinal.  The first case is $T$ unsuperstable or is NDOP or
NDTOP, the second case are the deep theories (which are superstable, NDOP,
NOTOP) and the third are the rest.

A theological question is which of those two dividing lines is the more
striking dividing line.  Probably between the upper case and the rest.
Clearly the fact that from the isomorphism type of a model of $T$ we can
naturally compute a stationary set modulo a club (see \scite{cft.10} below),
getting any such set, say that the class of models of $T$ is very 
complicated, whereas
a tree with $\omega$ levels seems reasonably understood though their number
(up to isomorphism) is large.   We can look at it in another way: if we
``understood" the isomorphism types of $M$, forcing notions ``which do not
do much damage" (including preserving inequality of cardinality of the 
relevant sets), preserve non isomorphism of models if $T$ is in the lower
or middle case.  E.g. if $\lambda = \text{ cf}(\lambda) > |T| = \aleph_0$ 
and $\Bbb P$ is a forcing notion not adding $\omega$-sequences to $\lambda$ 
preserving cardinalities $\le \lambda$ \ub{then} $\Bbb P$ preserves non 
isomorphism of models of $T$ of cardinality
$\le \lambda$ iff $T$ is in the lower of middle case.  It seems a very weak
demand of a complete set of invariants to be preserved by such a change in
the universe.  This is the intended meaning of the word (main) gap here,
though to say that the isomorphism types of models of $T$ are all ``simple",
``well understood" is open to variations, here the ``good, well understood"
case is very good, and the ``bad" are so bad, that it is an evidence to this
dividing line to be a major natural division (on c.c.c. forcing - see
\scite{cft.15} below).  E.g. we may above require the forcing to add no
$(< \lambda)$-sequences getting the same division.

The audience asked \nl
\ub{\stag{cft.9} Question}:  Can we assign stationary sets as invariants?
\bn
\ub{\stag{cft.10} Answer}:  In restricted classes of models it works but the
question is what the connection should be between the model and the
stationary set.  That is, generally, there are enough stationary sets to
code models in cardinality $\lambda$, so we have to say $M,S(M)/{\Cal D}
_\lambda$ (or ${\Cal D}_{\le \kappa}(\lambda)$ or whatever) should be nicely
connected.  Hence this remains vague.  Note that if we aim not at a complete
set of invariants but as an evidence for nonstructure, then we can.  That is,
for any $T$ in the upper case we can naturally assign a stationary subset 
of $\lambda$ modulo 
${\Cal D}_\lambda$ as an invariant to models of $T$ of cardinality 
$\lambda = \text{ cf}(\lambda)$ such that any 
stationary subset of $\lambda$ (or of $\{\delta < \lambda:
\text{cf}(\delta) = \aleph_0\}$) appears.  E.g. let $T$ be unsuperstable.
If say $M$ has universe
$\lambda = \text{ cf}(\lambda) > \aleph_1$, use $\{\delta < \lambda:M
\restriction \delta \prec M$ and for every $\bar b \in {}^{\omega >}M$ every
countable subtype of tp$(b,M_\delta,M)$ is realized in $M_\delta\}$. \nl
[Why?  Let $\lambda = \text{ cf}(\lambda) > |T| + \mu$, where $\mu =
\text{ cf}(\mu) > \aleph_0$.  Let $\Phi$ be proper template for
$K^{\text{tr}}_\omega$, EM$_{\tau(T)}(I,\Phi)$ a model of $T$, witnessing
unsuperstability, let $I$ be a linear order of the form $\lambda + J,J$
isomorphic to the inverse of $\mu$ and for 
$\delta \in S^* = \{\delta < \lambda:\text{cf}(\delta)
= \aleph_0\}$ let $\eta_\delta$ be an increasing $\omega$-sequence of ordinals
with limit $\delta$.  Now for $S \subseteq S^*$ let $I_S = {}^{\omega >}
(\lambda + J) \cup \{\eta_\delta:\delta \in S\} \cup \{\eta:\eta \in {}^\omega
(\lambda + J)$ is eventually zero$\}$.  We can check that the invariant of
EM$(I_S,\Phi)$ is $S/{\underset\tilde {}\to D_\lambda}$.]

However, this is not an invariant which characterizes up to isomorphism.
The cases of NDOP, NOTOP are in face easier (can use the end of
\cite[Ch.III,\S3]{Sh:e}).
\bn
Classifying will not die as \nl
\ub{\stag{cft.11} Thesis}:  In any reasonable classification (in the present
sense) there are examples of the ``complicated" class which are actually
well understood so should be prototypes of another class which is
analyzable.
\bn
Hodges in his thesis had asked about \nl
\ub{\stag{cft.12} Question}:  When does a first order 
theory have a $\prec$-minimal model in $\lambda$?  What can be PrSp$(T)$?  
\bigskip

\definition{\stag{cft.13} Definition}  1) $M$ is a $\prec$-minimal 
model in $\lambda$ if it can be elementary embeddable into any other 
model of $N$ of Th$(M)$ of cardinality $\lambda$. \nl
2) PrSp$(T) = \{\lambda:M \text{ has a } \prec \text{-minimal model of } T\}$.
\enddefinition
\bn
We may consider
\definition{\stag{cft.14} Definition}  PrSp$'(T) = 
\{(\lambda,\mu):M \text{ a model of }
T \text{ of cardinality } \mu \text{ which is}$ \nl
$(\lambda,\prec)\text{-embeddable}\}$ where 
$M$ is $(\lambda,\prec)$-embeddable if it is embeddable into every model of
Th$(M)$ of cardinality $\lambda$.
\sn
Hodges gave some examples of PrSp$(T)$ and then I add a few others.  
Hodges showed
that if $T$ is the theory of infinite atomic Boolean Algebra, then PrSp$(T) =
\{\lambda:\lambda \text{ is strong limit}\}$.  Also if $I$ is a linear 
order with no monotonic sequence of elements of length
$\lambda,\lambda = \text{ cf}(\lambda)$, then in $EM(I,\Phi)$ there is 
no formula
defining when restricted to some set, a well order of length $\lambda$. 
\enddefinition
\bn
My old remarks and a theorem of Hrushovski that ``unidimensional stable
theory is superstable" gives \nl
\ub{Fact}:  Assume $T$ is stable, cf$(\lambda) > |T|$; if $\lambda \in
PrSp(T)$ then $T$ is unidimensional (hence superstable) and cf$(\mu) > |T|
\Rightarrow \mu \in PrSp(T)$. 
\mn
[why? \nl
\ub{Case 1}:  $T$ not unidimensional. \nl
As in \cite[Ch.V,\S2]{Sh:c}.
\mn
\ub{Case 2}:  $T$ is unidimensional, superstable.  As in 
\cite[Ch.IX,\S2]{Sh:c}.
\bn
\centerline {$* \qquad * \qquad *$}
\bn
As discussed above we 
know that for complicated theories (say unstable or unsuperstable or ones
with DOP or OTOP), models can code stationary sets hence ``isomorphism types
are very set theoretic sensitive".  E.g. changing by forcing shooting a club 
of $\lambda = \text{ cf}(\lambda)$ disjoints to some $S$.  It is natural to
consider: ``nice forcing notion can make non-isomorphic models to
isomorphic".
\bn
\ub{\stag{cft.15} Problem}:  For which 
first order countable $T$ is there a c.c.c. forcing
notion making non-isomorphic isomorphic.

Now (see \cite{BLSh:464}) for unsuperstable $T$ the answer is no and so with
superstable with DOP or OTOP.

By \cite{LwSh:518}, there is $T$ among the remaining with the answer no,
after preliminary c.c.c. forcing.  Laskowski and me agree there is a serious
unsatisfactory point in the paper, but do not agree on its identity.  
He thinks it is the preliminary c.c.c. forcing.  However,
I think that as anyhow we deal with forcing this is minor, but the
restriction on $T$ major (whereas he think not).
\sn
Newelski told me that for $T$ superstable countable,
$p \in S^m(A)$ with uncountably many stationarization (\cite[Ch.III]{Sh:c}), 
he considered the meagre ideal on the set of stationarizations of
$p$ (in connection with Vaught conjecture for superstable $T$).
Subsequently in \cite{LwSh:560} we used the null ideal on this space for
proving in the proof (could have used the meagre).
He asked
\mn
\ub{\stag{cft.16} Question}:  Are 
there $T,A,p$ as above such that the ideal of null sets
and of meagre sets are different?
\bn
\centerline{$* \qquad * \qquad *$}
\bn
We can measure the number of models in other ways.
\definition{\stag{cft.17} Definition}:  Let 
$I'_\kappa(\lambda,K) = \{M / \equiv_{\infty,\kappa}:M \in K,\|M\|=\lambda\}$
where $K$ is a class of models of a fixed vocabulary, $\tau(K) \equiv_{\infty,
\kappa}$ is the equivalence relation of having the same $L_{\infty,\kappa}$-
theory.  If $\kappa = \aleph_0$ we may omit it, if $K = \{M:M \models T\}$
we may write $T$.
\enddefinition
\bn
Starting my Ph.D. studies, I note (concentrating on 
$\kappa = \aleph_0$ but Rabin was not enthusiastic)
\proclaim{\stag{cft.18} Theorem}  If $K$ is elementary 
(or defined by $\psi \in
L_{\lambda^+,\kappa}$), $\lambda = \lambda^{< \kappa} \ge |\tau(K)|$ and
$I_\kappa(\lambda,K) \le \lambda$, \ub{then} $\lambda \le \mu
\Rightarrow I_\kappa(\mu,K) \le I_\kappa(\lambda,K)$. \nl
(Later I understand that this is easy by Levy absoluteness; see
\cite{Sh:11} and see Nadel's thesis). 
\endproclaim
\bn
So \nl
\ub{\stag{cft.19} Question}:  For first order $T$, what can

$$
\text{Min}\{\lambda:I_\kappa(\lambda,T) \le \lambda = \lambda^\kappa\}
$$
\mn
be when 
it is $< \infty$, i.e. can you give a better bound than the Hanf number
of $L_{(\tau(T)+\aleph_0)^+,\omega}$ (well ordering) $L_{\infty,\omega}$.

Lately, Laskowski and me investigate what can be the supremum of the
$L_{\infty,\kappa}$ Karp height for models of $T$, so a theory is considered
complicated if this is not bounded; this is closely related, see 
\cite{LwSh:560}.  The point is that while case $\kappa = \aleph_0$ is
opaque the cases of many bigger $\kappa$ is at least at present, not a dead
end, supporting \scite{Mon.11a}.
\newpage

\head {\S5 Unstable first order theory} \endhead  \resetall \sectno=5
 % \resetall 
\bn
The major theme of classification theory has been for me, since \cite{Sh:1}:
\nl
\ub{\stag{cus.1} Meta Problem}:  Find worthwhile dividing lines on the family
of (complete first order) theories.  \nl
A dividing line is not just a good,
interesting property, it is one for which we have something 
to say on both sides;
so for some problems naturally a solution goes by working on each side
separately.

Of course phrased as ``find dividing lines among the possible mathematical
theories" this is too general and too vague to lead to mathematical theorems.
But it is quite natural to restrict ourselves to the family of classes of
models of first order theories (complete, and even countable).

Almost by definition, a dividing line is an interesting property (though not
inversely: the class of non groups among $(A,F),F$ an associative two-place
function or non 0-minimal first order theories are not so remarkable), but
it is remarkable that, for our contexts, there are some.  I have changed the
name of \cite{Sh:a} from ``stability and the number of nonisomorphic models"
to ``classification theory"in order to stress its aim - finding meaningful
dividing lines.

We believe good test problems are needed and, of course, problems on the
number of non-isomorphic models were inherently interesting and serve well.
But they could not serve for unstable theories.  We shall see below how some
problems succeed or fails in this role, but sometime we do not know of a good
candidate.  I have considered at various times $\lambda \mapsto
\text{ sup}\{|S(A)|^+:|A| \le \lambda\}$, Keisler's order (i.e. saturation
of ultraproducts), SP (see below) and later $\triangleleft^*$ 
(\cite{Sh:500}) and
the existence of universal models.
Sometime getting a full ZFC answer (on which I work hard in \cite{Sh:c}) seems
too much so decide that it is reasonable to content myself with: \nl
\ub{\stag{cus.2} Half ZFC or Poor Man ZFC Answer}:  The result on the lower
half of a dividing line will be ZFC (or semi-ZFC, i.e. depending on cardinal
arithmetic in relevant cardinals), whereas in the complicated, upside we 
allow consistency results (in semi-ZFC: 
you may distinguish between cases to high consistency
strength and those really consistent you may argue to add diamond, etc.).  

This may help, as getting a too fine division is not our aim.
Also if we are more interested in the dividing lines themselves, consistency 
results should be enough.  This is even more relevant in classifying
non-elementary classes and in classification over a predicate. \nl
Note that if we look at ``having complicated phenomena" as barrier to
positive theorems, clearly a consistency result suffices.
\bn
\centerline{$* \qquad * \qquad *$}
\bn
\ub{\stag{cus.2a} Discussion}:  I find 
it particularly nice if the property have some equivalent definitions
by ``outside notions" and ``by inside notions", some got for dealing with
the ``down side", some with ``the upside".  To clarify consider the example of
stability; unstable theories are characterized by the order property (inside
property for the upside, helpful in proving the class of models of an
unstable $T$ is complicated), stable theories are characterized by 
having finite local ranks $(R^n(p,\varphi,2) < \omega)$ (inside property 
for the downside; helpful in developing stability theory, showing we can in
some senses understand the class of models of a stable $T$),
instability is characterized by ``for every 
$\lambda$ for some $A,|S(A)| > \lambda \ge |A|$" (a weak
outside property for the upside), stability by 
``for every $\lambda,|A| \le \lambda
\Rightarrow |S(A)| \le \lambda^{|T|}$ (a weak outside property for 
the downside); late coming outside property characterizing 
``unstable $T$" is, ``has many $\kappa$-resplended models of 
cardinality $\lambda$" where $\lambda = \lambda^{< \kappa} >
2^{|T|}$, (outside property for the upside), ``stable $T$ has exactly one
$|T|^+$-resplended model of cardinality $\lambda$ when $\lambda = 
\lambda^{|T|}$" (outside property for downside; see \cite[Ch.V]{Sh:e}).
\bn
\centerline{$* \qquad * \qquad *$}
\bn
Considering unstable theories, we knew they have the independence property
or the strict order property, but not necessarily both, so the simplest
prototypes of unstable theories are the T$_{ord}$, the theory of dense 
linear order and the theory $T_{rg}$ of random graph.  
We have earlier in \S2 discussed $(k,m)-*$-NIP and 
it is natural to ask on
the inter-relations of them, the strict order property and the independence
property, see \cite{Sh:715}.
\bn
For the neighborhood of $T_{rg}$, the problem I had chosen as a 
test problem was \nl
\ub{\stag{cus.3} Problem}:  Classify first order theories by

$$
\align
SP(T) = \bigl\{(\lambda,\kappa):&\text{ every model } M \text{ of } T
\text{ of cardinality } \lambda \\
  &\text{ has an elementary extension } N \text{ of cardinality } \lambda \\
  &\text{ which is } \kappa \text{-saturated}\bigr\}
\endalign
$$
\mn
or, for simplicity

$$
SP'(T) = \{(\lambda,\kappa) \in SP(T):\lambda^{2^{|T|}} = \lambda > 2^\kappa
\text{ and } \kappa = \text{ cf}(\kappa) > 2^{|T|}\}.
$$
\mn
[Why ``for every $M$", not just there is $M$?  Because then, letting
$T = Th(M),M = M_1 + M_2,T_\ell = Th(M_\ell)$ and $T_2$ trivial 
(e.g. Th$(\omega,=)$, that is, having infinite models, all
relations empty) we easily can check that $SP'(T)$ is maximal; that is, 
equal to $=\{(\lambda,\kappa):\lambda = \lambda^{2^{|T|}},\lambda > 
2^\kappa,\kappa = \text{ cf}(\kappa) >
2^{|T|}\})$ and the intended intuition is to say that 
$T,T_1$ has the same complicatedness.

Now \cite{Sh:93} give a semi-ZFC answer to the question 
on for which $T$ is SP$'(T)$ is minimal
(i.e. are maximally complicated under this criterion).
\bigskip

\proclaim{\stag{cus.4} Theorem}  If $T$ is not simple, \ub{then} 
$SP'(T)$ is minimal (that is, is equal to $\{(\lambda,\kappa):
\lambda = \lambda^\kappa,\kappa = \text{ cf}(\kappa)
> 2^{|T|}\})$. 
\endproclaim
\bn
The other directions, if $T$ is not simple (hence having the
tree property) then SP$(T)$ is minimal, holds by \cite{Sh:c}.
For this \cite{Sh:93} began the generalization \cite[Ch.II,III]{Sh:a} to simple
theories, I suggested to some to continue but only lately Hrushovski
\footnote{his preprint has not appeared (and, unlike the others, will not), 
it investigates the generalization
of ``geometric stability theory" and group interpretations for theories
minimal for $D(-,L,\infty)$.  He has some theories of fields and
investigating finite models with few types in mind} and 
then Kim, Pillay Laskowski, Buechler, Morgan, Shami
and others use and investigate parallels of
\cite{Sh:c} to simple theories; for surveys see \cite{GIL97x}, \cite{KiPi}.

We expect that \nl
\ub{\stag{cus.5} Conjecture}:  There is a 
finer division of simple theories to $\omega +1$
families, by the properties Pr$_n$ such that if $T$ has the Pr$_n$-th
property for every $n$ then $SP'(T)$ is minimal, two theories of the same
family (i.e. satisfying Pr$_n$ but not Pr$_{n+1}(T)$ and let $n=n(T)$) 
essentially have the same SP$'(T)$, but if two have different $n(T)$
then consistently there is a cardinal separating
them (in the SP$(T') \backslash \text{ SP}(T'')$ if $n(T'') < n(T')$);
this should be the relatively easy part.
A prototype of a counter-example to the $n$-th property, $k \ge 3$ is the
model completion of $T_k$, where $T_k$ say: $(x_1,\dotsc,x_n)$ is symmetric
irreflexive, $R(x_1,\dotsc,x_n) \rightarrow P(x),xSy \rightarrow P(x) \and
\neg P(y)$ and $\neg(\exists x,\dotsc,x_ny)(\neg P(y) \and
R(x_1,\dotsc,x_n) \and \dsize \bigwedge^n_{\ell =1} x_\ell Sy)$.

The intention is that Pr$_n(T)$ is a syntactical property which implies:
\mr
\item "{$\boxtimes$}"  if $M$ is a model of $T,\mu = \mu^{|T|} \le \|M\|
\le 2^\mu$ and $M^* \prec M,\|M^*\| \le |T|,p^* \in S(M^*),\Phi =: \{
\varphi(x,\bar a):\bar \in M \text{ and } \{\varphi(x,\bar a)\} \cup p^*$
does not fork over $M^*\}$, \ub{then} can be represented as
$\dbcu_{i < \mu} \Gamma_i$ and $\varphi_1,\dotsc,\varphi_n \in \Gamma_i 
\Rightarrow p^* \cup \{\varphi_1,\dotsc,\varphi_n\}$ does not fork over
$M^*$.  (For $n=2$, \cite[7.8]{Sh:93} is a version).  We may use \cite{Sh:234}.
\endroster
\bn
What is $SP(T_{rg})$? \nl
For any $\mu$ let log$(\mu) = \text{ Min}\{\lambda:2^\lambda \ge \mu\}$ 
if, now if $\mu \ge (\text{log}(\mu))^{< \kappa}$ then 
$(\lambda,\kappa) \in SP(T_{rg})$ (using \cite{EK}).  If 
(log$(\mu))^{< \kappa} > \mu$ this is conected to SCH.
By \cite{GiSh:597} the answer is independent.  Note that $T_{rg}$ is minimal
among simple theories in the sense that $SP(T)$ is maximal among unstable
theories.
\bn
It is not a priori clear that the answer is so coherent, there may be a myriad
of properties with many independent results; I have not tried this direction.
This will not help us much in
classification.  Here I am not sure if the ``armies of god" will prevail.
In other words, I am not sure it is a good test problem any more.
\bn
\centerline{$* \qquad * \qquad *$}
\bn
As said above any unstable theory has the independence property or strict
order property.  So among unstable theories the theory of random graphs and
the theory of linear order are in some sense the simplest.  So we can expect
to have a theory of some family of first order theories for which linear 
order is a prototype (as discussed earlier for theory of random graph).  Best,
of course, is if we can have something for all $T$ without the independence
property (see after \scite{Mon.7}).  It was encouraging
(\cite[Ch.III,\S7]{Sh:c}).
\bigskip

\proclaim{\stag{cus.10} Theorem}  ($T$ first order 
without the independence property). \nl
If $\kappa > |T|$ is regular $A \subseteq {\frak C}$ \ub{then}
we can find a $\kappa$-saturated
$M \prec {\frak C}$ such that $A \subseteq M$ and $M$ is in some sense
constructible over $A:|M| = A \cup \{a_i:i < \alpha\}$ and tp$(a_i,A \cup
\{a_j:j < i\})$ does not split over some $B_i \subseteq A \cup \{a_j:j < i\}$
which has cardinality $< \kappa$.
\endproclaim
\bn
For long there was no reasonable candidate for test question: 
(the results on $\lambda \mapsto \sup\{|S(A)|:|A| \le \lambda
\}$ were satisfactory but do not lead to something).  
Now \cite{LwSh:560} start to classify
by the $L_{\infty,\kappa}$-Karp height; note that some superstable theories
are maximal there.
\bn
We may look for a parallel of \cite{Sh:93}, e.g.
\definition{\stag{cus.10a} Definition}  Assume 
$T$ without the independence property and $\lambda =
\lambda^{< \lambda} + |T| < \mu,M \in EC(\mu,T),N \prec_{L_{|T|^+,|T|^+}}M,
\|N\| = 2^{|T|},p^* \in S(N)$ and 

$$
\align
{\Cal P}_\ell = {\Cal P}_\ell(p^*,M) =: \biggl\{ p \in S(M):&p 
\text{ in some sense does not fork over } p^* \\
  &\text{ which means that } p^* \subseteq p \text{ and}: \\
  &(a) \quad \text{ if } \ell =1 \text{ letting }
P_\varphi = P_{\varphi(x,\bar y)} = \{\bar c \in {}^{\ell g(\bar y)}M:\varphi
(x,\bar a) \in p\} \\
  &\qquad \text{ we have }
(N,P_\varphi)_\varphi \prec_{L_{|T|^+,|T|^+}} (M,P_\varphi)_\varphi \\
  &(b) \quad \text{ if } \ell = 2, \text{ then for every } \varphi(x,\bar a)
\in p \text{ and for } A \subseteq N,|A| \le |T| \\
  &\qquad \text{ there is } \varphi(x,\bar a') \in p^* \text{ such that }
\bar a,\bar a' \\
  &\qquad \text{ realizes the same type over } A \text{ inside } M \biggr\} 
\endalign
$$

$$
{\Cal F}_\ell(p,M) = \cup\{p:p \in {\Cal P}_\ell(p^*,M)\}.
$$
\mn
We can make ${\Cal F}_\ell(p,M)$ to a Boolean Algebra (as in the later
parts of \cite{Sh:93}).\enddefinition
\bn
\ub{\stag{cus.11} Question}:    Can you force this Boolean Algebra by a 
$\lambda^+$-c.c. $\lambda$-complete forcing notion to be ``simple" in some 
sense? best: union of $\lambda$ subalgebra which are interval Boolean
Algebra.

Probably too much to hope for but the direction may be reasonable, see more
\scite{cus.41} - \scite{cus.44} and \cite{Sh:715}.
%\enddemo
\bn
\centerline{$* \qquad * \qquad *$}
\bn
Not having the strict order property look to me a priori very promising 
dividing line, however, the test problems which look promising lead to smaller
classes (see below on \cite[\S2]{Sh:500}).  This includes
\definition{\stag{cus.12} Definition}  1) The universality spectrum of $T$ is

$$
\align
UvSp(T) = \{\lambda:&T \text{ has a universal model, i.e. every other model 
of} \\
  &T \text{ of cardinality } \lambda \text{ can be elementarily embedded
into it}\}
\endalign
$$
\mn
2)  The pairs-Universality spectrum of $T$ is

$$
\align
\text{UvpSp}(T) = 
\{(\lambda,\mu):&\lambda \le \mu \text{ and there } M \in EC
(\mu,T) \\
  &\text{ into which every } N \in EC(\lambda,T) \text{ can be} \\
  &\text{elementarily embedded}\}
\endalign
$$
\mn
(we can look at the size of a universal family; cov sheds light on the
connections, see \cite{Sh:457}, \cite{DjSh:614}).
\sn
Under GCH for $\lambda > |T|$ the answer is known, so we can look only for
weak solutions involving consistency, ``semi-ZFC solutions" as suggested
in \scite{cus.2} above. \nl
Now the theory of universal graphs consistently has large universal spectrum
even for $\lambda < 2^{< \lambda}$ (\cite{Sh:175a}).  
So once we know (\cite{KjSh:409}) that the theory of linear
order has few (e.g. $2^\lambda > \lambda^{++} \Rightarrow \lambda^{++}
\notin \text{ UvSp}(T_{\text{ord}})$), and that this 
applies to any $T$ with the
strict order property, it raises hope that this is a good test problem for
that property.  
\enddefinition
\bn
Alas, it may be good but not for the strict order property
as (\cite[\S2]{Sh:500}) NSOP$_4$ suffices where (see
\cite{Sh:500}, \cite{DjSh:692}):
\definition{\stag{cus.13} Definition}  1) $T$ has the SOP (the strong order
property) if some type $p(\bar x,\bar y)$ defined in ${\frak C}_T$ a partial
order with arbitrarily long $(< \bar \kappa)$ chains). \nl
2) $T$ has SOP$_n$ (the strong $n$-order property, $n \ge 3$) if for some
formula $\varphi(\bar x,\bar y)$:
\mr
\item "{$(a)$}"  there is an infinite indiscernible sequence ordered by
$\varphi$
\sn
\item "{$(b)$}"  we cannot find $m \le n$ and $\bar a_0,\dotsc,\bar a_{m-1},
\bar a_m = \bar a_0$ such that $\varphi(\bar a_\ell,\bar a_{\ell +1})$ for
$\ell < m$.
\ermn
3) $T$ has SOP$_2$ if some $\varphi(\bar x,\bar y)$ has it which means that
we can find in ${\frak C}_T,\bar a_\eta \in {}^{\ell g(\bar y)}({\frak C}_T)$
for $\eta \in {}^{\omega >}2$ such that:
\mr
\item "{$(a)$}"  if $\eta \char 94 \langle \ell \rangle \trianglelefteq
\eta_\ell \in {}^{\omega >}2$ for $\ell = 0,1$ then $\{\varphi(\bar x,\bar a
_{\eta_1}),\varphi(\bar x,\bar a_{\eta_0})\}$ is inconsistent
\sn
\item "{$(b)$}"  if $\eta \in {}^\omega 2$ then $\{\varphi(\bar x,
\bar a_{\eta \restriction n}:n < \omega\}$ is inconsistent.
\ermn
4) $T$ has SOP$_1$ is defined as in (3) only in clause (a) we demand
$\eta_0 = \eta \char 94 \langle 0 \rangle$.
\enddefinition
\bn
\ub{\stag{cus.14} Problem}:  1) Develop a theory for NSOP $T$'s. \nl
2) Develop a theory for NSOP$_n \, T$'s.  \nl
3) Find additional evidence of complicatedness to the SOP$_n$'s (and SOP).
\nl
Earlier I thought that the most 
promising is the case $n=3$,
a prototypical theory seems $T_{feq}$ (\cite{Sh:457}), but now we know that
$n=2$ is a real dividing line (\cite{DjSh:692}).  However, we have 
SOP$_n \Rightarrow$ SOP$_{n+1}$ and for $n \ge 3$ the inverse implication
fails, but for $n=1,2$?
\bn
Now \cite[\S1]{Sh:457} indicates another direction, see Dzamonja Shelah 
\cite{DjSh:710}; there for theory $T$ with tree coding we prove 
some non-existence of universal models.
\bigskip

\definition{\stag{cus.15} Definition}  1) The formula 
$\varphi(\bar x,\bar y,\bar z)$ is tree coding in $T$, 
if for every (equivalently some) $\lambda \ge \kappa \ge
\aleph_0$ we can find in ${\frak C}_T,\bar c_\nu(\nu \in {}^\kappa \lambda),
\bar b_\eta(\eta \in {}^{\kappa >}\lambda)$, \nl
$\bar a_\alpha(\alpha < \kappa)$
such that:
\mr
\item "{$(a)$}"  $\models \varphi[\bar c_\nu,\bar b_\eta,\bar a_\alpha]$
if $\eta = \nu \restriction \alpha \and \nu \in {}^\kappa \lambda$
\sn
\item "{$(b)$}"  if $\alpha < \kappa$ and $\nu,\rho \in {}^\kappa \lambda$ and
$\nu \restriction \alpha \ne \rho \restriction \alpha$ then $\varphi(
\bar c_\nu,\bar y,\bar a_\alpha),\varphi(\bar c_\rho,\bar y,\bar a_\alpha)$
are contradictory.
\ermn
2) $T$ has tree coding if some $\varphi(\bar x,\bar y,\bar z)$ has. 
\enddefinition
\bn
\ub{\stag{cus.16} Problem}:  Develop the theory of $T$'s without tree coding
(and further nonstructure theorems for those with).
\nl
Clearly in some sense the dividing line stable/unstable is simpler than
superstable/unsuperstable not to mention NDOP/DOP, etc.  
The following definitions tries explicate this.  The point being that 
many properties are properties of a formula $\varphi(\bar x,\bar y)$ in $T$.
\bigskip

\definition{\stag{cus.17} Definition}  Fix $T$ and ${\frak C} =
{\frak C}_T$. \nl
1) For first order formula 
$\varphi = \varphi(\bar x,\bar y)$ and $\bar a_0,\dotsc,
\bar a_{n-1} \in {}^{\ell g(\bar y)}{\frak C}$, let $\bold u_\varphi
(\bar a_0,\dotsc,\bar a_{n-1}) = \{u \subseteq n:\models(\exists \bar x)
\dsize \bigwedge_{\ell < n} \varphi(\bar x,\bar a_\ell)^{\text{cf}(i \in
u)}\}$ where $\varphi^{\text{true}} = \psi,\psi^{\text{false}} = \neg \psi$.
\nl
2) For first order formula $\varphi = \varphi(\bar x,\bar y)$ let

$$
\align
\Gamma_\varphi = \{(n,\bold u):&\text{for some } n < \omega \text{ and }
\bar a_0,\dotsc,\bar a_{n-1} \in {}^{\ell g(\bar y)}{\frak C}; \\
  &\text{we have } \bold u_\varphi(\bar a_0,\dotsc,\bar a_{n-1}) =
\bold u\}.
\endalign
$$
\mn
3) We let ${\bold \Gamma}_T = \{\Gamma_\varphi:\varphi = \varphi(\bar x,
\bar y) \in L(T)\}$. 
A division of first order theories is straightly defined if: for some
${\bold \Gamma}$ it is the family of $T$ such that ${\bold \Gamma} \in 
{\bold \Gamma}_T$ the first order $T$ are divided to those $T$'s that 
${\bold \Gamma} \in {\bold \Gamma}_T$ (the up sets) and those $T$'s that 
${\bold \Gamma} \notin \bold \Gamma_T$ (the down side). \nl
4) Let ${\bold \Gamma}_{T,\varphi(\bar x,\bar y,\bar z)} = 
\{\Gamma_{\varphi(\bar x,\bar y,\bar c)}:\bar c \in 
{}^{\ell g(z)}{\frak C}\}$.
\enddefinition
\bigskip

\definition{\stag{cus.18} Definition}  1) For $\Gamma$ as above we say: $T$
is $\Gamma$-high if $\Gamma \in \Gamma_T$ and $\Gamma$-low otherwise. \nl
2) We say that a class ${\frak T}$ of complete first order theories is
straight if the truth value of $T \in {\frak T}$ is determined by $\Gamma_T$.
\enddefinition
\bn
A variant which seems to capture the main point is:
\definition{\stag{cus.19} Definition}  1) Let $\Gamma^* = \{(n,F_1,F_2):
n < \omega,F_1,F_2$ are disjoint families of sets of partial functions 
from $\{0,\dotsc,n-1\}$ to $\{$true false$\}\}$. \nl
2) We say that $\langle \bar a_0,\dotsc,\bar a_{n-1}\rangle$ does
$\varphi$-realizes 
$(n,F_1,F_2) \in \Gamma^*$ if $f \in F_1 \Rightarrow {\frak C}
\models (\exists \bar x) \dsize \bigwedge_{\ell \in \text{ Dom}(f)}
\varphi(\bar x,\bar a_\ell)^{f(\ell)}$ and $f \in F_2 \Rightarrow {\frak C}
\models \neg(\exists \bar x)\dsize \bigwedge_{\ell \in \text{ Dom}(f)}
\varphi(\bar x,\bar a_\ell)^{f(\ell)}$.  We can apply this to
$\varphi(\bar x,\bar y,\bar c),\bar c$ from ${\frak C}_T$.
\nl
3) For $\Gamma \subseteq \Gamma^*$ we say that $\varphi(\bar x,\bar y)$ has
the weak $\Gamma$-property (in $T$) if any $(n,F_1,F_2) \in \Gamma$ is 
$\varphi$-realized by some $\langle \bar a_0,\dotsc,\bar a_{n-1} \rangle$.
We say that $\varphi(\bar x,\bar y)$ has the strong $\Gamma$-property if for
$(n,F_1,F_2) \in \Gamma^*$ we have $(n,F_1,F_2) \in \Gamma$ iff
$(n,F_1,F_2)$ is $\varphi$-realized by some $\langle \bar a_0,\dotsc,
\bar a_{n-1} \rangle$.
We say $\varphi(\bar x,\bar y,\bar z)$ has such a property if this holds 
for some $\varphi(\bar x,\bar y,\bar c)$.  \nl
4) $T$ has the weak/strong $\Gamma$-property if some 
$\varphi(\bar x,\bar y,\bar z)$ has it.
$T$ has the weak/strong pure $\Gamma$-property if some 
$\varphi(\bar x,\bar y)$ has it.
\nl
5) We say that a class (or property) ${\frak T}$ of complete first order
theories is weakly/strong simply high straight if for some $\Gamma \subseteq 
\Gamma^*$ we have: $T \in {\frak T}$ if $T$ has the weak/strong 
$\Gamma$-property in $T$.  The class
${\frak T}$ is weakly/strongly simply low straight if it is the 
compliment of a simply high straight one. \nl
6) Omitting the ``weak" and ``strong" we shall mean weak(ly). \nl
7) $\boldsymbol\Gamma^{**} = \{\Gamma \subseteq \Gamma^*$: for some $T$, some
$\varphi(\bar x,\bar y)$ has the $\Gamma$-property$\}$.
\enddefinition
\bn
\ub{\stag{cus.20} Fact}:  1) For any 
$\Gamma \subseteq \Gamma^*$, the truth of ``$T$ has the weak
$\Gamma$-property" is determined by ${\bold\Gamma}_T$. \nl
2) Allowing in Definition \scite{cus.18}(10, (2), \scite{cus.19} for the
weak versions, formulas
$\varphi(\bar x,\bar y,\bar c)$, does not make a difference for having the
$\Gamma$-property.
\bn
\ub{\stag{cus.21} Observation}:  The following 
properties can be represented as ``$T$ has the weak $\Gamma$-property". \nl
1) $T$ is unstable. \nl
2) $T$ has the independence property. \nl
3) $T$ has the strict order property. \nl
4) $T$ has the tree property (equivalently, is not simple). \nl
5) $T$ has NSOP$_n$ (the $n$-strong order property)(where $n \ge 3$). \nl
6) $T$ has the NSOP$_2$. \nl
7) $T$ has the NSOP$_1$.
\bigskip

\demo{Proof}  Only 5) is not immediate. 

It suffices to show
\mr
\item "{$(*)$}"  $T$ has NSOP$_n$ \ub{iff} for some $\varphi = \varphi
(\bar x,\bar y)$ we have
\sn
\item "{$(*)_\varphi$}"  there is an indiscernible set $\langle \bar a_\ell:
\ell < \omega \rangle$ such that
{\roster
\itemitem{ $(a)$ }  for every $0 = i_0 < \ldots < i_n = \omega,\{\varphi
(\bar x,\bar a_m)^{\text{if}(\ell \text{ even})}:\ell < n,m \in [i_\ell,
i_{\ell +1})\}$ is consistent
\sn
\itemitem{ $(b)$ }  for no $i_0 < \ldots < i_n$ do we have $\models(\exists
\bar x) \dsize \bigwedge_{\ell \le n} \varphi(\bar x,\bar a_{i_\ell})
^{\text{if}(\ell \text{ even})}$ \nl
${{}}$ \hfill$\square_{\scite{cus.21}}$
\endroster}
\endroster
\enddemo
\bn
In this context, there are naturally the most complex theories:
\definition{\stag{cus.22} Definition}  1) We say that 
$\varphi(\bar x,\bar y)$ straightly
maximal (in $T$) if $\Gamma_\varphi$ is maximal. \nl
2) We say $\varphi(\bar x,\bar y,\bar z)$ is strongly straightly maximal
(in $T$) if ${\bold \Gamma}_{T,\varphi(\bar x,\bar y,\bar z)}$ is maximal.
3) Call $T$ straightly maximal if some $\varphi(\bar x,\bar y)$ is. \nl
4) Call $T$ strongly straightly 
maximal if some $\varphi(\bar x,\bar y,\bar z)$ is.
\enddefinition
\bn
An example is true arithmetic, i.e. $Th(\omega,0,1,+,\times)$ \nl
\ub{\stag{cus.23} Problem}:  1) Develop a theory 
\mr
\item "{$(a)$}"  for non-straightly maximal $T$
\sn
\item "{$(b)$}"  for non-strongly straightly maximal.
\ermn
2) Find natural nonstructure theorem, i.e. witness for having complicated
models
\mr
\item "{$(a)$}"  for straightly maximal $T$'s
\sn
\item "{$(b)$}"  for strongly straightly maximal $T$'s.
\endroster
\bn
\centerline {$* \qquad * \qquad *$}
\bn
Now \scite{cus.23} seems quite persuasive to me, \ub{but} I have to say
I do not know of a test problem, nor what
should we expect of a good theory for the nonmaximal theory.  Note that this
scheme does not include the $(k,n)-*$-NIP where ``arity" is important.

We can easily adapt the definitions to include it, but the present version
is not necessarily a drawback - the present version does not discriminate
elements from seven-tuples, etc., and \nl
\ub{\stag{cus.23a} Thesis}:  It is certainly reasonable to map
the continents and oceans before we look at hills and lakes (if we can,
of course).

Now superstability does not fit this scheme, too, again it is a finer
distinction; yet, we write down this version.
\bigskip
\definition{\stag{cus.24} Definition}  Fix $T$ and ${\frak C} = {\frak C}_T$
and $\alpha$ an ordinal. \nl
1) Let $\bar \varphi = \langle \varphi_i(\bar x,\bar y_i,\bar c_i):i <
\alpha \rangle$ be a sequence with $\bar c_i$ from ${\frak C}_T$. \nl
For $n < \omega$ and $\bar a_{i,\ell} \in {}^{(\ell g(\bar y_i))}{\frak C}$
let

$$
\align
\bold u_{\bar \varphi}( \langle \bar a_{i,\ell}:i < \alpha,0 < n \rangle)
= \biggl\{ u \subseteq \alpha \times n:& \text{ the type } \{\varphi
(\bar x,\bar a_{i,\ell}) {}^{\text{if}((i,\ell) \in u)}: \\
  &i < \alpha,\ell < n\} \text{ is consistent}\biggr\}.
\endalign
$$
\mn
2) For a sequence $\bar \varphi = \langle \varphi_i(\bar x,\bar y_i,
\bar c_i):i < \alpha \rangle$ as above let

$$
\align
\Gamma_{\bar \varphi} = \bigl\{(n,\bold u):&\text{for some } n \text{ and }
\bar a_{i,\ell} \in {}^{(\ell g(\bar y_i)}{\frak C} \\
  &\text{(for } i < \alpha,\ell < n) \text{ we have} \\
  &\bold u_{\bar \varphi}(\langle \bar a_{i,\ell}:i < \alpha,\ell < n \rangle)
= \bold u \bigr\}.
\endalign
$$
\mn
3) Let $\Gamma^\alpha_T = \{\Gamma_{\bar \varphi}:\bar \varphi = \langle
\varphi_i(\bar x,\bar y_i):i < \alpha \rangle\}$. \nl
4) Let $\Gamma^*_2 = \{(n,F_1,F_2):n < \omega \text{ and } F_1,F_2$ are
disjoint families of partial finite functions from $\alpha \times n$ to
$\{$true,false$\}\}$. \nl
5) We say that $\langle \bar a_{i,\ell}:i < \alpha,\ell < n \rangle$ does
$\bar \varphi$-realizes $(n,F_1,F_2) \in \Gamma^*_\alpha$ if is as above
and:
\mr
\item "{$(*)_1$}"  $f \in F_1 \Rightarrow {\frak C} \models (\exists \bar x)
[\dsize \bigwedge_{(i,\ell) \in \text{ Dom}(f)} \varphi_i(\bar x,
\bar a_{i,\ell})^{f(i,\ell)}]$
\sn
\item "{$(*)_2$}"  $f \in F_2 \Rightarrow {\frak C} \models \neg(\exists 
\bar x)[\dsize \bigwedge_{(i,\ell) \in \text{ Dom}(f)} \varphi_i(\bar x,
\bar a_{i,\ell})^{f(i,\ell)}]$.
\ermn
6) For $\Gamma \subseteq \Gamma^*$ we say that $\bar \varphi = \langle
\varphi_i(\bar x,\bar y_i,\bar c_i) \rangle$ has the $\Gamma$-property
in $T$ if every $(n,F_1,F_1) \in \Gamma$ is $\bar \varphi$-realized for $T$.
\nl
7) We say that $T$ has the $\Gamma$-property for $\Gamma \subseteq
\Gamma^*_\alpha$ if some $\bar \varphi = \langle \varphi_i(\bar x,\bar y_i):
i < \alpha \rangle$ has the $\Gamma$-property.
\enddefinition
\bn
What about DOP, OTOP, etc?
\definition{\stag{cus.25} Definition}  For a logic ${\Cal L},(T,{\Cal L})$
has any of the properties defined above if we allow the formulas $\varphi$
to be in $L$ but as ${\Cal L}$ possibly fail compactness we should like
large case so:
\mr
\item "{$(a)$}"  $(T,{\Cal L})$ has the order property if for some
$\varphi(\bar x,\bar y) \in {\Cal L}$ with $\gamma =: \ell g(\bar x) =
\ell g(\bar y)$, for every linear order $I$ there are a model $M$ of $T$ and
$\bar a_t \in {}^\gamma M$ for $t \in I$ such that, for any $t,s \in I$ we
have $M \models \varphi[\bar a_t,\bar a_s]$ iff $t <_I s$
\sn
\item "{$(b)$}"  $(T,{\Cal L})$ has the independence property if for some
$\varphi(\bar x,\bar y) \in {\Cal L}$, for every $\lambda$ there are a model
$M$ of $T,\bar a_\gamma \in {}^{\ell g(\bar y)} M$ for $\gamma < \lambda,
\bar a_u \in {}^{\ell g(\bar x)} M$ for $u \subseteq \lambda$ such that
$M \models \varphi[\bar a_u,\bar b_\gamma]$ iff $\gamma \in u$.
\ermn\enddefinition
But if we look at NDOP or NOTOP, (for superstable $T$, in the standard
definition) we do not fully use ${\Cal L} = 
L_{|T|^T,|T|^*}$ or ${\Cal L} = L_{|T|^+,\aleph_0},{\Cal L} =
L_{|T|^+,\aleph_0}(Q^{\text{dim}\ge \aleph_1})$ (if we use finite sequences,
sufficient for superstable $T$, see \cite[Ch.XII]{Sh:c}) or 
${\Cal L} = L_{|T|^+,|T|^+}$, we rather use formulas of specific form.
But the order property and independence property becomes equivalent, and
main gap tend to show equivalence of such versions.
\bn
\centerline{$* \qquad * \qquad *$}
\bn
The universality spectrum raises many problems both set theoretic and model
theoretic.  For the set theoretic side, we still do not know enough on
UvSp$(T_{\text{ord}})$ and also the universe UvSp for the theory of graphs 
(see \cite{DjSh:659}).
\bn
\ub{\stag{cus.30} Problem}:  Is it consistent that for some 
$\lambda,\lambda,\mu,\lambda < \mu < 2^\lambda$ and 
$\mu \in \text{ UvSp}(T)$ for every countable $T$? (equivalently true
arithmetic).
\bn
Still theories with SOP$_4$ look essentially maximal (as the results on
linear orders hold for them) \nl
\ub{\stag{cus.31} Problem}:  Does every NSOP$_4$ theory 
$T$ have consistently a non-trivial universality spectrum?
\bn
\ub{\stag{cus.32} Thesis}:  The way, (a good way), a reasonable way to develop
the theory of NSOP$_4$ and/or NSOP$_4$ first order theories is
\mr
\item "{$(*)$}"  start by asking which first order theories 
fall under the persuasion of \cite[\S4]{Sh:457}? 
\ermn
So we know that if $T$ has SOP$_4$ then not, whereas if $T$ is simple, 
then yes.  Lastly, for some $T$ which has SOP$_3 \and$ NSOP$_4$, the 
answer is yes.
\sn
The nice scenario is if those will be exactly SOP$_4$.  
If this succeeds, this will be
very good for investigating universality spectrum.  It may give something
on the theory of NSOP$_4$.  Maybe, a right parallel of non-forking.  If 
it fails, it still gives important information on universality.  May give 
information on NSOP$_3$.
\bn
\ub{\stag{cus.33} Discussion}:  One may pose the question: is 
universality just a tool toward classifying?
\bn
\ub{Answer}:  In some sense, yes.

But, I believe the right way to classify is to choose a worthwhile relevant
test problem (like number of non-isomorphic models).  So it is true that 
in a sense the classification is higher, real aim but still the 
universality spectrum and classifying 
$I(\lambda,T)$ are very important.  Reason for 
optimism concerning the universality spectrum is: 
the positive and negative answers
(guessing clubs) and \cite[\S4]{Sh:457} seems to speak on the same thing.
\bn
\ub{\stag{cus.34} Question}:  CON(in $\lambda^+$ there is universal
linear order $\and 2^\lambda > \lambda^+ \and \lambda = \lambda^{> \lambda} 
> \aleph_0$).  If this fails, we can look at the examples in \cite[\S2]{Sh:500}
(existentially complete directed graph with no $(\le k)$-cycle).
\sn
Of course: we would like to ask for which first order theories the proof 
in \cite{Sh:457}, \cite{DjSh:614} will work?
\nl
For PA (piano arithmetic)?  Conceivably for PA we can prove that: there is
no universal in more cardinalities than the obvious ones ($\lambda =
2^{< \lambda}$, where $\lambda > |T|$ for simplicity) or can try there all 
theories.  If we fail for linear order but succeed for some other $T$'s, 
it should be very illuminating, maybe revealing new dividing lines. 
\bn
I have not looked at \nl
\ub{\stag{cus.35} Question}:  Does 
all simple unstable countable theories have the same
universality spectrum? 
Or, do they have many possible spectrums?

The natural way: look at forcing for graphs and think of a non-trivial simple
theory such that if in the beginning we force many models of it in $\chi$,
there would not be co-habitation.
\mn
If we discover too fine a distinction, it will not be so exciting to
investigate.
\bn
Even so, a \ub{Major} question is \nl
\ub{\stag{cus.36} Problem}: Find the maximal class for UnSp, that is a 
dividing line in the sense that they behave like linear order (at present).  
\nl
\ub{If} for all first order $T$ we have the
consistency hoped for linear order, \ub{but} many such theories behave
differently \ub{and} there is no alternate proofs for ``there is no
universal" in ZFC (+ cardinal arithmetic), \ub{then} finer distinction among
such theories look not inviting.
\mn
My feeling: the dividing line of the proof in \cite[\S2]{Sh:500} is a 
major dividing line, the one for universality. 
\bn

To get semi-ZFC distinction \nl
\ub{\stag{cus.37} Question}:  Generalize \cite{Sh:457} to $\lambda^{++}$.
\nl
Clearly for NSOP$_4$ theories and probably more
this fails; i.e. we get some notion but the property required in \cite{Sh:457}
fails; but this may provide a theory of types to NSOP$_4$ theory (or a new
dividing line).
\bn
Of course, we may like to know more on simple theories \nl
\ub{\stag{cus.38} Problem}:  For 
which theories the consistency results on graphs
(\cite{Sh:175}, \cite{Sh:175a}, \cite{DjSh:659}?) can be generalized?
\bn
Even for graphs (but probably not hard): \nl
\ub{\stag{cus.39} Problem}:  Can we in 
\cite{Sh:175a} get the consistency for all regular
cardinals in the intervals? also for the singulars? 
\bn
\centerline {$* \qquad * \qquad *$}
\bn
\ub{\stag{cus.40} Discussion}: 
In the spectrum from in the one end finding the bare outlines, 
finding some order in the total chaos, to the other end, perfectly 
understanding on what we know not little, I prefer the first.  So 
though I was (and am still) sure 
that there is much more to be said on 
superstable/stable theories (in fact, this essentially follows from the belief
that it is an important dividing line) not to say on theories of finite
Morley rank, and on simple theories, I am more excited from starting 
new frameworks.

Of course, I believe that such general theorems of f.o. theories will have
meaningful application for specific theories (though I do not agree with A.
Robinson that this is the aim of model theory or a needed justification; 
but I agree it is a worthwhile
one), in fact, such applicability is highly suggestive from belief in the
meaningfulness of the dividing line (if the theory is serious).  Well, some
may argue that has
not simple theories proved to be the only one with reasonable non-forking
(by Kim and Pillay \cite{KiPi})?  Yes, but this had been done for stable, too,
and maybe trying to generalize is not the only way to find an understanding
of such theories.
For example, probably the theory of NSOP$_3$ theories will replace elements
by formulas, and we shall have to make parallel replacement moving from
NSOP$_n$ to NSOP$_{n+1}$.  E.g. consider: for a formula $\varphi(\bar x,
\bar a)$ in $M_1$ and $M_0 \prec M_1$ and type $p \in S(M_1)$ to which
$\varphi(x,\bar a)$ belongs, as in \scite{cus.10a} $q \subseteq p 
\restriction M_0,|q| \le |T|,\varphi(x,\bar a)$ reflect nicely in $M_0$.
However, in some sense having proved the main gap
for countable f.o. theory, I feel my task (on first order theories) was 
done, just like \cite{Sh:460} in cardinal arithmetic.
\bn
\centerline{$* \qquad * \qquad *$} 
\bn
In linear order, if $\langle a_t:t \in J \rangle$ is indiscernible ($\equiv$
monotonic) over $A,t_0 < s < t_1$ and $\{t_0,s,t_1\} \subseteq J$, \ub{then}
tp$(a_s,\{a_{t_0},a_{t_1}\}) \vdash \text{ tp}(a_s,A)$.
\bn
\ub{\stag{cus.41} Question}:  Can we prove a similar phenomena for 
NIP theories?

This cannot be literally true as for stable theories it is false.  Probably
we should ``divide" the works between stable like parts and the above idea.

On the other hand putting together intervals of length $|T|$ and adding we
can find $\langle b'_t:t \in I \rangle$ such that $\bar b_t \subseteq
\bar b'_t$ and for $t_0 < t_1 < t_2$, tp$(\bar b_{t_1},\bar b'_{t_0} 
\char 94 \bar b'_{t_2}) \vdash \text{ tp}
(\bar b_{t_1},\cup \{\bar b_s:s \notin (t_0,t_2)\}$.

In some sense, a model of a stable theory $M$ can be represented by a well
ordering and unary functions: \nl
\ub{\stag{cus.41a} Fact}:  If Th$(M)$ is 
stable and $|M| = \{a_\alpha:\alpha < \alpha^*\}$, we can find 
$f_{\varphi,\ell}:\alpha^* \rightarrow \alpha^*$ satisfying
$f_{\varphi,\ell}(\alpha) < \text{ Max}\{2,\alpha\}$ 
(for $\varphi = \varphi(x,\bar y) \in L(\tau_T),\ell < n_\varphi < \omega$) 
such that tp$(\langle a_{\alpha_1},\dotsc,a_{\alpha_n} \rangle,
\emptyset,M)$ can be reconstructed from equalities between composition of
$f_{\varphi,\ell}$ (the point being that tp$_\varphi(a_\alpha,\{a_\beta:
\beta < \alpha\})$ is definable by some $\psi(\bar y,\bar c),\bar c \subseteq
\{a_\beta:\beta < \alpha\}$).
\bn
\ub{\stag{cus.42} Problem}:  1) For NIP theories, does 
something parallel hold with equalities replace by some 
($\le |T|$) linear orderings of $\alpha^*$? \nl
2) Find parallel theories for other properties of $T$.
\bn
\ub{\stag{cus.43} Problem}:  Investigate first order $T$ which are NIP (i.e.
without the independence property).
\bn
\ub{\stag{cus.44} Question}:  For $T$ with NIP: \nl
1) If $A \subseteq B \subseteq {\frak C}_T,p \in S(A)$, does there exist
$q \in S(B)$ extending $p$ which does not fork over $A$? \nl
2) Do ordered groups play here a role similar to groups for stable theories?
\bn
\ub{\stag{cus.45} Question}:  For (complete) $T$ with the independence
property, $T_1 \supseteq T$, and $\theta$ and for simplicity $\lambda$ a
successor of regular $> 2^\theta$, are there $\theta$-resplened models $M_1$
of $T_1$ with $M_1 \restriction \tau_T$ has large $L_{\infty,\lambda}$-Karp
height?
\newpage

\head {\S6 Classifying non-elementary classes} \endhead  \resetall \sectno=6
 % \resetall 
\bigskip

I see this as the major problem of model theory.
Cherlin presses me to expand on this point; now in '69 Morley and Keisler told
me that model theory of first order logic is essentially done and the future
is the
development of model theory of infinitary logics (particularly fragments of
$L_{\omega_1,\omega}$).  By the eighties it was clearly not the case and
attention was withdrawn from infinitary logic (and generalized quantifiers,
etc.) back to first order logic.  Now, of course, it is better to prove
theorems in a wider context, also we may recall that algebraists are not
restricting their attention to elementary classes; but wider context may 
have a heavy price in content, it is not clear that there interesting 
theory left at all.  As the theory for the family of
first order theories has widened and deepened this attention was justified.
But, of course, it would be wonderful if we have at all a 
classification theory for nonelementary classes.  Just 
generalizing with changes here and there is
not so exciting, but clearly, if 
there is a theory at all, there are in it many
dividing lines of different character; the danger it is the other direction:
having too weak theory.

Of course , this is phrased too generally, e.g. I feel 
classes defined by $\psi \in L_{\aleph_1,\aleph_1}$ are
probably hopeless (we can easily code behaviour which are very set
theoretically sensitive).  So 
``non elementary" should be restricted to a reasonable class, and
there are choices.
The first case I considered was $(K_D,\prec)$
where
\definition{\stag{nec.1} Definition}  Let $T$ be a first order complete theory,
$D \subseteq D(T) = \cup \{D(M):M$ a model of $T\}$ where $D(M) =
\{\text{tp}(\bar a,\emptyset,M):\bar a \in {}^m M,m < \omega\}$,
(so $D$ codes $T$, well when $D \ne \emptyset$).  Let \nl
1) $K_D = \{M:M \text{ a model of } T \text{ (so } \tau(M)=\tau(T)) \text{ and
moreover } D(M) \subseteq D\}$ (and $\prec$ is the usual being elementary
submodel order). \nl
2) $M$ is $\lambda$-sequence-homogeneous (or just $\lambda$-homogeneous)
\ub{if} for every elementary map $f$ of $M$ (i.e. $f$ one to one from
Dom$(f) \subseteq M$ to Rang$(f) \subseteq M$ and $f$ preserve first 
order formulas)
of cardinality $< \lambda$ and $a \in M$ there is an elementary map $f'$
of $M$ satisfying $f \subseteq f' \and a \in \text{ Dom}(f')$. \nl
3) $M$ is $(\lambda,D)$-homogeneous \ub{if} $M$ is $\lambda$-homogeneous
and $D(M)=D$.
\enddefinition
\bn
The reason for considering $K_D$ was that ``$(\lambda,D)$-homogeneous" was
similar to ``$\lambda$-saturated".  The older notion of model
homogenous had not looked managable to me (see \scite{nec.2}(1),(2) below).
\definition{\stag{nec.2} Definition}  1) $M$ is 
$\lambda$-model-homogeneous if: for every isomorphism $f$ from
$M_1 \prec M$ onto $M_2 \prec M,M_1 \prec N_1 \prec M,\|N_1\| < \lambda$
there is $N_2,M_2 \prec N_2 \prec M$ and an isomorphism $f'$ from $N_1$ 
onto $N_2$ extending $f$. \nl
2)  $M$ is model homogenous if it is $\|M\|$-model homogeneous. \nl 
3) ${\Cal D}_\kappa(M) = \{N / \cong:N \prec M,\|N\| \le \kappa\}$. \nl
4) $M$ is $(\lambda,{\Cal D})$-homogeneous if $M$ is $\lambda$-model
homogeneous and ${\Cal D}_{|\tau({\Cal D})|+\aleph_0}(M) = {\Cal D}$. \nl
5) $K_{\Cal D} = \{M:{\Cal D}_{|\tau({\Cal D})|+\aleph_0}(M) \subseteq
{\Cal D}\}$.  
\enddefinition
\bn
Still we do not know the answer to \nl
\ub{\stag{nec.3} Question}:  1) Is there a ``reasonable" upper bound to

$$
\align
\mu^*_\kappa = \biggl\{ \text{Min}\{\lambda:&\text{there is no }
(\lambda,D) \text{-homogeneous model of cardinality} \\
  &\lambda\}:\text{ for some complete first order theory } T \text{ of
cardinality } \kappa, \\
  &D \subseteq D(T) \biggr\}.
\endalign
$$
\mn
2) Similarly for $(\lambda,{\Cal D})$-homogeneity. \nl
I think that it is known (by the Kazachstan school, under GCH) that 
$\mu^*_{\aleph_0} \ge \aleph_\omega$.
\bn
But more central for me is \nl
\ub{\stag{nec.4} Problem}:  1) How much of the theory on stable theories can
be generalized to $(K_D,\prec)$ for stable $D$?  \nl
2) Similarly for superstable; where
\bn
\definition{\stag{nec.4a} Definition}  1) $K_D$ is 
stable if (for every $\lambda$ there is a
$(\lambda,D)$-homogeneous model of cardinality $\ge \lambda$, and) for
arbitrarily large $\lambda,K_D$ is stable in $\lambda$ which means
$A \subseteq M \in K_D,|A| = \lambda \Rightarrow
S(A,M) = \{\text{tp}(a,A,M):a \in M\}$ has cardinality $\le \lambda$. \nl
2) $K_D$ is superstable if the stability holds for every large enough
$\lambda$. \nl
Investigation of $(K_D,\prec)$ have been carried, see the introduction
of \cite{HySh:676}. \nl
There is little on $(\lambda,{\Cal D})$-homogeneity (see \cite{Sh:237c},
\cite{Sh:300}).
The interest is mainly in $D$ such that for every $\lambda$ there is
$(\lambda,D)$-homogeneous model of cardinality $\lambda$, but anyhow in
definition \scite{nec.4a}, it suffices to deal with ``small $\lambda$", the
rest follows.
\enddefinition
\bn
\ub{\stag{nec.5} Problem}:  1) Prove the main gap for

$$
I(\lambda,K_D) = \{M / \cong:M \in K_D,\|M\| = \lambda\}.
$$
\mn
2) Prove the main gap for

$$
\align
I(\lambda,\{M \in K_D:M \text{ is }(\kappa,D) \text{-homogeneous}\}) = 
|\{M / \cong:&M \in K_D \text{ has cardinality } \lambda \\
  &\text{and is } (\kappa,D) \text{-homogeneous}\}|.
\endalign
$$
\mn
Certainly for first order classes I considered as the main case version (1)
(note: when $D = D(T)$ we get back the elementary classes as special cases).
However, here the interest started with $(D,\mu)$-homogeneous model so
probably part (2) is more natural.  However, the problem has not been
resolved even for countable first order $T,D = D(T)$; see \cite[\S0]{HySh:676}
on what was done.
\sn
What is lost in this context compared with the first order one?  Formulas
are not so interesting any more, except as part of a complete type.  There
is a remnant of compactness: there is $a \in {\frak C}_D$ realizing a type
$p \in S(A)$ iff for every finite $B$ the type $p \restriction B$ is
realized.  Also the Hanf numbers of omitting types is helpful and
$(D,\kappa)$-homogeneous is quite parallel to $\kappa$-saturated; large
parts of stability theory for such models has been generalized to this 
context and much more is still to be done.  
Note that it should not all be parallel
to the first order case, first there are new aspects (like 
$\lambda$-goodness), also some early work was done first in this context
(the stability spectrum $\kappa(D)$ and $\lambda(D)$) 
and lastly, something like
${\frak C}^{eq}$ may be better here in some respects.

Another direction has been universal class, where a class $K$ of $\tau$-models
closed under isomorphism is a universal class when $M \in K$ iff every
finitely generated submodel belongs to $K$ (see \cite{Sh:300},
\cite{Sh:h}).  This context is incomparable with first order; a universal
class is certainly not necessarily first order, and also the inverse 
implication fails.  Now there may be long sequences on which a quantifier 
free formula defines order, in
which case we have a strong nonstructure.  Otherwise we can define being a
submodel $M \le N$, axiomatize the setting and start developing
the parallel of \cite{Sh:c}, with types being defined by chasing arrows
rather than as a set of formulas, starting with the parallel of the theorem
``saturated $\equiv$ homogeneous universal", and having some new dividing 
lines, getting regular types, etc.
The idea was that assuming some possible reasons for strong nonstructure
does not hold, we can define a stronger notion of submodel $<_1$ (like
$\prec_{\Sigma_1}$) and prove that ${\frak K}^+ = ({\frak K},<_1)$ is 
inside our setting.
We think that after enough such strengthening, the intersection is similar
enough to the first order case to prove the main gap, but this was not done.
\bn
\ub{\stag{nec.x} Question}:  Does the main gap (of course with depth possibly
quite large) hold for universal classes?

Note that though first order formulas does not play a role, types, dimension
of indiscernible sets, prime models, orthogonality and regularity does.  
Also we believe
that the idea of changing inductively the context will be helpful (as it is
in \cite{Sh:600}).

We may rather look at classes defined say by $\psi \in L_{\omega_1,\omega}$,
here it is harder to begin.

Note that generally in this section I have thought that we should 
expect not just the situation
in cardinals $\lambda \le |T|$ to be different than in ``large enough
$\lambda$" (as was the case for first order) but say $\lambda <$ relevant 
Hanf number of $L_{\omega_1,\omega}$, so the small cardinal should have 
different behaviour.  The theory is not totally empty as we 
can prove some things:
\bigskip

\proclaim{\stag{nec.6} Theorem}  Assume $2^{\aleph_n} < 2^{\aleph_{n+1}}$ for
$n < \omega$. \nl
1) If $\psi \in L_{\omega_1,\omega}$ have ``few" models in $\aleph_1,\dotsc,
\aleph_n$ (essentially  $I(\aleph_m,\psi) < 2^{\aleph_m}$) but has an 
uncountable
model \ub{then} $\psi$ has a model in $\aleph_{n+1}$. \nl
2) If $\psi \in L_{\omega_1,\omega}$ have few models in $\aleph_1,\aleph_2,
\dotsc,\aleph_n,\dotsc(n < \omega)$
but has an uncountable model, \ub{then} $\psi$ has models on all
cardinalities. \nl
3) If $\psi \in L_{\omega_1,\omega}$ is categorical in $\aleph_1,\aleph_2,
\ldots,\aleph_n,\dotsc,(n < \omega)$, \ub{then} $\psi$ is categorical in 
every $\lambda > \aleph_0$; in fact under the assumption of part (2), $\psi$
is excellent, and for excellent classes categoricity is one $\lambda >
\aleph_0$ suffice here
(essentially \cite{Sh:87a}, \cite{Sh:87b} when ``few" is strengthened a
little, see more in \cite{Sh:600}, more on excellent class \cite{GrHa89}).
\endproclaim
\bn
We do not know: \nl
\ub{\stag{nec.7} Problem}:  If $\psi \in L_{\omega_1,\omega}$ (or even
$\psi \in L_{\kappa^+,\omega}$) is categorical in one $\lambda \ge
\beth_{\omega_1}$ (or $\lambda \ge \beth_{(2^\kappa)^+}$), \ub{then} $\psi$
is categorical in every such $\lambda$?
\mn
Some wonder why ``$\lambda \ge \beth_{\omega_1}$"?  Now $\lambda \ge
\aleph_\omega$ is necessarily as by \cite{HaSh:323}, $\psi$ may categorical
in $\aleph_0,\dotsc,\aleph_n$, but not in $\lambda$ if $2^{\aleph_n} <
2^\lambda$ (or so).
\mn
Others wonder why such modest question, isn't the main gap better?  Of
course it is, but I think it is more reasonable first to resolve the
categoricity.  But are ``a class of models of $\psi \in L_{\kappa^+,\omega}$"
the best context?  Thinking of putting \cite{Sh:87a} + \cite{Sh:87b}
together with results on $L_{\omega_1,\omega}(Q)$ in \cite{Sh:48},  I 
consider (\cite{Sh:88}) abstract elementary classes. 
I have preferred this context, certainly the widest I think has any chance
at all.  \nl
In \cite {Sh:87a}, \cite{Sh:87b}, \cite{Sh:88} it is proved:
\mr
\item "{$(*)_2$}"  catgoricity (of $\psi \in L_{\omega_1,\omega}(Q))$ in
$\aleph_1$ implies the existence of a model of $\psi$ of cardinality
$\aleph_2$;
\sn
\item "{$(*)_3$}"  if $n > 0,2^{\aleph_0} < 2^{\aleph_1} < \ldots <
2^{\aleph_n},\psi \in L_{\omega_1,\omega}$ and $1 \le I(\aleph_\ell,\psi)
< \mu_{\text{wd}}(\aleph_\ell)$ for $1 \le \ell \le n$, \ub{then}
$\psi$ has a model of cardinality $\aleph_{n+1}$
\sn
\item "{$(*)_4$}"  if $2^{\aleph_0} < 2^{\aleph_1} < \dotsc,\psi \in
L_{\omega_1,\omega}$ and $1 \le (\aleph_\ell,\psi) < \mu_{\text{wd}}
(\aleph_\ell)$ for $\ell = 1,2,\ldots$, \ub{then} $\psi$ has a model in
every infinite cardinal and is categorical in one $\lambda > \aleph_0$ iff
it is categorical in every $\lambda > \aleph_0$.
\ermn
Now the problems were: \nl
\ub{\stag{nec.8} Problem}:  Prove $(*)_3,(*)_4$ in the context of an
abstract elementary class ${\frak K}$ which is $PC_{\aleph_0}$.
\bn
\ub{\stag{nec.9} Problem}:  Parallel results in ZFC; e.g. prove $(*)_3$ 
when $n = 1,2^{\aleph_0} = 2^{\aleph_1}$.
By \cite[\S6]{Sh:88} there are classes categorical in $\aleph_1$ if MA, but
not so if $2^{\aleph_0} < 2^{\aleph_1}$ so really there is here a different
model theory involved.
\bn
\ub{\stag{nec.10} Problem}:   Construct examples; e.g. ${\frak K}$
(or $\psi \in L_{\omega_1,\omega}$), categorical in $\aleph_0,\aleph_1,\dotsc,
\aleph_n$ but not in $\aleph_{n+1}$ (see \cite{HaSh:323}).
\bn
\ub{\stag{nec.11} Problem}:   If ${\frak K}$ is $\lambda$-a.e.c. (abstract
elementary class), and is categorical in $\lambda$ and $\lambda^+$, does it
necessarily have a model in $\lambda^{++}$? assuming $2^\lambda <
2^{\lambda^+} < 2^{\lambda^{++}}$?  In \cite{Sh:576} we solve a somewhat
weaker version of \scite{nec.11}.
\bn
It is reasonable to be willing to assume 
large cardinal, if we can develop some interesting theory.
In \cite{MaSh:285} a version of Los Conjecture for $T \subseteq
L_{\kappa,\omega},\kappa$ compact cardinal was proved (starting for large
enough successor). \nl
\ub{\stag{nec.12} Question}:  1) If $T \subseteq L_{\kappa^+,\omega}$ is
categorical in one limit $\lambda > \beth_{(2^{\kappa+|T|})^+}$, \ub{then}
$T$ is categorical in every $\lambda \ge \beth_{(2^{\kappa+|T|})^+}$. \nl
2) Similarly for ${\frak K}$ a $\kappa$-a.e.c. with amalgamation. \nl
3) Similarly for ${\frak K}$ as $\kappa$-a.e.c. \nl
Note: that for (2) there are some results (\cite{Sh:394}).
\bn
Moreover \nl
\ub{\stag{nec.13} Problem}:  1) Develop classification (or at least stability) 
theory for $T \subseteq
L_{\kappa,\omega}$ at least if $\kappa$ is compact, or even just measurable.
\nl
In Kolman Shelah \cite{KlSh:362}, \cite{Sh:472} the parallel 
(downward part) is proved for $\kappa$-measurable.
\bn
Several cases lead to \nl
\ub{\stag{nec.14} Problem}:  
Classify $\Phi$ proper for linear order (more accurately
$(\Phi,\tau),\tau \subseteq \tau(\Phi))$ according to the function
$I(\lambda,K_{\Phi,\tau})$ where

$$
K_{\Phi,\tau} = \{EM_\tau(I,\Phi):I \text{ a linear order of cardinality }
\lambda\}.
$$
\mn
Probably as a first step we should consider generic $I \subseteq
({}^\lambda 2,<_{lex})$ of cardinality $\lambda$ (and then try to work in
ZFC).  Maybe it is reasonable to restrict ourselves to a dense family of
$\Phi$'s, see \cite{Sh:394}.
\bn
\ub{\stag{nec.15} Problem}:  More 
interesting classes to serve as index models.  We have considered linear 
orders, trees
with $\kappa + 1$ levels, ordered graphs (see \cite[Ch.III,end of \S2]{Sh:e},
\cite{LwSh:560}).

If $T \subseteq L_{\kappa,\omega},\kappa$ compact have compactness for
$L_{\kappa,\kappa}$-types and can prove (under categoricity or a failure
of a nonstructure assumption) that $\prec_{L_{\kappa,\omega}} =
\prec_{L_{\kappa,\kappa}}$.  But when we consider e.g. $\kappa$-a.e.c. with
amalgamation, we may have a formal description of a type $p \in S(M)$ having
$p \restriction N$ wherever $N \le_{\frak K} M$ has small cardinality,
neither knowing it there a $\le_{\frak K}$-extension of $M$ in which it is
realized; not knowing it is unique.  Remember the type was defined by chasing
$\le_{\frak K}$-embedding.

In \cite{Sh:576} we consider whether we can do anything
without any remnant of compactness (i.e. without E.M.-models, no large
cardinals, no omitting type theorems) with some success.  This is continued
in \cite{Sh:600}, where we look at an abstract version of superstability
(proved to occur in ``nature" relying on earlier work.
\bn
\centerline {$* \qquad * \qquad *$}
\bn
There may be, however, limitations.  First order logic was characterized
e.g. by Lowenheim Skolem to $\aleph_1 +$ compactness, now those are the first
step, and we may well have the parallel of the theory without having the
basic properties (Lowenheim Skolem and compactness).
\bn
\ub{\stag{nec.16} Problem}  Can we characterize what part of stability
theory are actually peculiar to first order?

We may consider generalizing the definitions and theorems on simple theories
(see \S5 particularly \scite{cus.3}, \scite{cus.4}).  
Now the context which seems less hostile is $(D,\lambda)$-homogeneous one 
(see the beginning of the section).
\bigskip

\definition{\stag{nec.17} Definition}  Assume $D$ is a finite diagram.
\nl
1) Let $\kappa_{\theta,\sigma}(D)$ be the first regular (for simplicity)
cardinal $\kappa$ such that there is no increasing continuous sequence
$\langle A_i:i \le \kappa \rangle$ of $D$-sets each of cardinality
$< \kappa + \theta$ and $p \in S^m_D(A_\kappa)$
such that for every $i < \kappa,p \restriction A_{i+1}$ does $(\theta,
\sigma)$-divide over $A_i$ (see below). \nl
2) We say that $p \in S_D(B)$ does $(\theta,\sigma)$-divide over $A$ if:

($A \subseteq B$ are $D$-sets and in some $D$-set $C \supseteq B$) \nl

there are $\bar b \in {}^{\theta >}B$ and sequence $\langle \bar b_t:
t \in I \rangle$ in which \nl

$\bar b$ appears, as $\bar b_{t^*},|I| = \sigma, 
C = B \cup \dbcu_t \bar b_t$ and there are no \nl

$D$-set $C_1 \supseteq C$ and $\bar d \in {}^m(C_1)$ such that:
\mr
\item "{$(*)$}"  letting $\gamma = \ell g(\bar b)$ and $q(x,\bar y) \in
S^{1 + \gamma}(B)$ be 
$\{\varphi(x,\bar y,\bar c):\bar c \subseteq B,\varphi(x,\bar b,\bar c) 
\in p\}$ we have: $\bar d$ realizes $q(x,\bar b_t)$ for at
least two $t \in I$.
\ermn
3)  If we omit $\sigma$ from (1) we mean $\beth((2^{\theta + \kappa
+ (\tau(D))})^+)$ and in (2) we mean $\beth((2^{|B|+(\tau(D))})^+)$.

The value of $|I|$ is to allow us to use the Hanf number for omitting type,
no point to increase further.  
\enddefinition
\bn
Of course
\proclaim{\stag{nec.17a} Claim}  1) If $p \in S_D(B)$ does $\theta$-divide
over $A$, \ub{then} $p$ does $\theta$-divide$^+$ over $A$ which means that we can
choose $\langle \bar b_t:t \in I \rangle$ to be an indiscernible sequence.
\nl
2) If $p \in S_D(B)$ does $\theta$-divide$^+$ over $A$ then for every
$\sigma$, we have $p \in S_D(B)$ does $(\theta,\sigma)$-divide over $A$. \nl
3) If $\kappa < \kappa_{\theta,\sigma}(D),\mu = \mu^{< \kappa} > \kappa +
\theta,D$ is $\mu$-good, \ub{then} we can find a $D$-set $A,|A| = \mu$, and
$p_i \in S_D(A_i),A_i \subseteq A,|A_i| = \kappa + \theta$ for $i < \mu
^\kappa$ such that for $i \ne j,p_i,p_j$ are contradictory, i.e. no
$p \in S_D(A_i \cup A_j)$ extend $p_i$ and $p_j$. \nl
4) If $\kappa_{\theta,\sigma}(D)$ is $\ge \beth((2^{\theta + \sigma +
|\tau(D)|})^+)$ then it is $\infty$. \nl
5) If $\kappa_\theta(D)$ is $\ge \beth((2^{\theta + |\tau(D)|})^+)$ then
it is $\infty$.
\endproclaim
\bn
In Definition \scite{nec.17}(2), we can demand, in
$(*)$ instead two, a fix $n < \omega$, we do real change.  If we ask $\mu
\ge \aleph_0$, the theorem on Hanf numbers are no longer helpful, but
weakened forms of the statement \scite{nec.17a}(3) holds.
\nl
We now may generalize the test problem from \cite{Sh:93}.
\proclaim{\stag{nec.18} Theorem}  Assume the axiom (Ax)$_\mu$ of \cite{Sh:80},
$2^\mu > \lambda > \mu,\lambda^{< \mu} = \lambda,\mu^{< \mu} = \mu$.

If $D$ is a good finite diagram and $\kappa_{\mu,\mu^+}(D) \le \mu$ and
$A$ is a $D$-set of cardinality $\le \lambda$ \ub{then} we can find a
$(D,\mu)$-homogeneous model $M$ into which $A$ can be embedded.
\endproclaim
\bn
However \nl
\ub{\stag{nec.19} Question}:  Is $\langle \kappa_{\theta,\sigma}(D):\theta,
\sigma \rangle$ characterized by few invariants?  Mainly, is $\langle
\kappa_\theta(D):\theta \rangle$ constant for $\theta$ large enough and
$\langle \kappa_{\theta,\theta^+}(D):\theta \rangle$.
\bn
\centerline{$* \qquad * \qquad *$}
\bn
This may be connected to the ${\Cal P}^{(-)}(n)$-diagram theme.  Looking at
the proof of Morley's theorem, it struck me as a phenomenal good luck that
categoricity could be gotten from a global property (saturation) rather
than by painstakingly analyzing the models.  A model $M$ of cardinality
$\lambda$, with vocabulary of cardinality $\mu$, can be represented by on an
increasing continuous elementary chain $\langle M_i:i < \lambda \rangle$ with
$M = \dbcu_{i < \lambda} M_i,M_i$ 
of cardinality $|i| + \mu$.  Now for each $i$, we have
to analyze $M_{i+1}$ over $M_i$, so we represent the model $\langle M_{i+1},
M_i \rangle$ as an increasing continuous elementary chain $\langle
(M_{i+1},j,M_{i,j}):j < \|M_{i+1}\| \rangle,\|M_{i+1,j}\| = \|M_{i,j}\| =
|j| + \mu$ and now our problem is to construct $M_{i+1,j+1}$ over $M_{i,j},
M_{i+1,j},M_{i,j+1}$, so we have to represent $(M_{i+1,j+1},M_{i+1},j,
M_{i,j+1},M_{i,j})$ by an increasing continuous elementary chain.  After
$n$ such stages we have a ${\Cal P}(n)$-diaigram $\langle M_u:u \in
{\Cal P}(n) \rangle$, for $n=0$ this is just $M$, i.e. $M_\emptyset = M$,
and for $\bar M = \langle M_u:u \in {\Cal P}(n+1) \rangle$ letting $\bar M^-
=: \langle M_u:u \in {\Cal P}(n) \rangle$ and $\bar M^+ =: \langle
M_{u \cup \{n-1\}}:u \in {\Cal P}(n) \rangle$, both are 
${\Cal P}(n)$-diagrams and $\bar M^- \prec \bar M^+$.  We can say $\bar M$
is a $(\lambda',{\Cal P}(n))$-diagram if in addition $\|M_u\| = \lambda'$ 
for $u \in {\Cal P}(n)$.

So to understand a model $M$ in $\lambda$, for each
$n < \omega$ and $\lambda' \in [\mu,\lambda)$ for each $(\lambda',{\Cal P}
(n))$-diagram $\langle M_n:u \in {\Cal P}(n) \rangle$ we have to understand
$M_n$ over $\bar M^* = \langle M_u:u \in {\Cal P}^-(n) \rangle$ where
${\Cal P}^-(n) = {\Cal P}(n) \backslash \{n\},\bar M^*$ is called a
$(\lambda',{\Cal P}^-(n))$-diagram.  So for categoricity, ``understand"
means in particular that it is essentially unique up to isomorphism (the
``essentially" hint that we may have ``time up to $\lambda$" 
to ``correct" some things).
What have we gained?  Just naturally we can prove statements by induction
on $\lambda'$:  a statement on ${\Cal P}^{(-)}(n)$-diagrams for 
all $n$ simultaneously (or for $\lambda = \mu^{+n}$, prove for
$(\mu^{+m},{\Cal P}(n-m))$!)   The gain is that the statement for 
$(\lambda',{\Cal P}(n))$ for $\lambda' > \mu$ naturally used 
$\lambda'' \in [\mu,\lambda')$ and ${\Cal P}(n+1)$.

To prove existence of a model in $\lambda$, we similarly prove by induction
on $\lambda' \in [\mu,\lambda)$ that a $(\lambda',{\Cal P}^-(n))$-diagram
can be completed to a $(\lambda',{\Cal P}(n))$-diagram.
\bn
Of course, we expect more conditions, complicating our induction.
\bn
\ub{\stag{nec.20} Thesis}:  For complicated problems (on say all cardinals)
we expect we need such a ${\Cal P}^-(n)$ analysis.

This scheme was used in \cite{Sh:87b} mentioned above, and also
\cite{SgSh:217}, \cite[Ch.XII]{Sh:c}, \cite{Sh:234}.  Returning to simple
finite diagrams, for proving goodness from good behaviour in small cardinals,
etc., this seems reasonable.  This also applies to the hopeful Pr$_n$ for
\scite{cus.5}.
\newpage

\head {\S7 Finite model theory \\ 
0-1 Laws} \endhead  \resetall 
 % \resetall 
\bn
Many were interested but hope is faint. \nl
\ub{\stag{fin.1} Problem}:  Find a logic with 0-1 law (or at least convergence
or at least with very weak 0-1 law) from which finite combinatorialist can
draw conclusion, novel for them.
\mn
But see \cite{Fri99}.  
We know that say for the random model $(n,<,{\Cal R}),{\Cal R}$ a random
2-place relation, the 0-1 law and even convergence fails (\cite{CHSh:245})
but the very weak 0-1 law holds (\cite{Sh:551}, a continuation with 
accurate estimates
Boppana Spencer \cite{BoSp}).  However, this positive result goes through
without telling us what first order formulas can define (in any random
enough such model).
\bn
\ub{\stag{fin.2} Question}:  Find the model theoretic content of the very
weak 0-1 laws for \nl
$(n,<,{\Cal R})$ and $(n,{\Cal F}),{\Cal F}$ a random 2-place
function.  \nl
We hope for a very weak ``elimination of quantifiers", saying
hopefully one which gives: first order formulas can
say much on ``small set", but little on the majority.
\bn
Let $G_{n,p}$ be the random graph with set of vertices $[n] = \{1,\dotsc,n\}$
and edge probability $p$.  It seems to me natural \nl
\ub{\stag{fin.3} Problem}:  1) Characterize the sequences $\langle
p_n:n < \omega \rangle$ of probabilities (that is reals in the interval
$[0,1)$) such that for every first order sentence $\psi$ in the language
of graphs we have
\mn
\ub{Possibility a}:  ($0-1$ law):

$\langle \text{Prob}(G_{n,p_n} \models \psi):n < \omega \rangle$ converge
to zero or converge to 1.
\mn
\ub{Possibility b}:  (convergence):

$\langle \text{Prob}(G_{n,p_n} \models \psi):n < \omega \rangle$ converge.
\mn
\ub{Possibility c}:  (very weak $0-1$ law):

$\langle \text{Prob}(G_{n+1,p_{n+1}} \models \psi) - \text{ Prob}
(G_{n,p_n} \models \psi):n < \omega \rangle$ converge to zero.
\mn
2) Like part (1) replacing $G_{n,p_n}$ by the $G_{n,\bar p}$, the random
graph with set of vertices $\{1,\dotsc,n\}$ and the probability of $\{i,j\}$
being an edge is $p_{(i-j)}$ (see \cite{LuSh:435}). \nl
3) Other cases (say random model on $\{1,\dotsc,n\}$ with vocabulary $\tau$).
\nl
A solution for \scite{fin.3}(2) case should be in \cite{Sh:581}.

In the cases of $0-1$ laws considered we usually get a dichotomy; say
${\Cal M}_n$ is the $n$-random structure, say on $\{1,\dotsc,n\}$; the
dichotomy has the form:  either (a) or (b) where
\mr
\item "{$(a)$}"  we have a complex case, i.e. we can define in ${\Cal M}_n$
(if $n$ large enough ${\Cal M}_n$ random enough) an initial segment of
arithmetic of size $k_{{\Cal M}_n}$, say of order of magnitude
$\sim \text{ log}(n)$ or at least log$_*(n)$ (or weakly complex:
${\underset\tilde {}\to k_n}$ going to infinity or at least for some
$\varepsilon > 0$ for every $k^*,\varepsilon < \text{ lim sup Prob}
(k_{\underset\tilde {}\to M_n} \ge k^*))$
\sn
\item "{$(b)$}"  we have a simple case, so $0-1$ law (or at least convergence)
(see \cite{Sh:550}).  So 
\endroster
\bn
\ub{\stag{fin.4} Problem}:  1) Prove for reasonable classes of 0-1
contexts $\langle {\Cal M}_n:n < \omega \rangle$ such dichotomies.
\nl
2) Investigate the family of $\langle {\Cal M}_n:n < \omega
\rangle$ which are nice (in the direction of having $0-1$ laws), like closure
under relevant operations. \nl
Concerning part (2), see \cite{Sh:550}, \cite{Sh:637}.
\bn
\ub{\stag{fin.5} Problem}:  In \S2 we discuss investigating reasonable
partial orders among generalized quantifiers.  Make a parallel investigating
on finite models.

See \cite{Sh:639} which try to do for the finite what \cite{Sh:171} do to
a large extent for the infinite case.
\newpage

\head {\S8 More on finite partition theorems} \endhead  \resetall \sectno=8
 % \resetall 
\bigskip

See discussion in \cite[\S8]{Sh:666}. 
\nl
\ub{\stag{app4.0} Question}:  What is the order of magnitude of the 
Hales-Jewitt numbers, $HJ(n,c)$ (see Definition \scite{app4.1}(3) below).
\bigskip

\definition{\stag{app4.1} Definition}  
1) Let $\Lambda$ be a finite nonempty alphabet, we define 
$f^{10}_\Lambda(m,c)$ where $m,c \in \Bbb N,|\Lambda|$ divide $m$, as the
first $k \le \omega$ divisible by $|\Lambda|$ such that:
\mr
\item "{$(*)$}"  if $d$ is a $c$-colouring of ${}^{[k]}\Lambda$, i.e. a
function from $\{\eta:\eta$ a function from $[k] = \{1,\dotsc,k\}$ 
into a set with $\le c$ members$\}$, \ub{then} we can find $\langle M_\ell:
\ell < m \rangle$ and $\eta^*$ such that:
{\roster
\itemitem{ $(a)$ }  $M_\ell \subseteq [k],\ell \ne m \Rightarrow M_\ell \cap
M_m = \emptyset$ and $\eta^*$ is a function from $M \backslash \dbcu_\ell
M_\ell$ into $\Lambda$
\sn
\itemitem{ $(b)$ }   $\|M_\ell\| = \|M_0\| > 0$ for $\ell < m$
\sn
\itemitem{ $(c)$ }  for $\nu_1,\nu_2 \in S = \{\eta:\eta \in {}^{[k]}\Lambda,
\eta^* \subseteq \eta$, and each $\eta \restriction M_\ell$ is constant$\}$
we have $d(\nu_1) = d(\nu_2)$ provided that for every 
$\alpha \in \Lambda$ for $i \in \{1,2\}$ we have \nl
$|\{\ell < m:\nu_i \restriction M_\ell$ is constantly $\alpha\}| = m/
|\Lambda|$
\sn
\itemitem{ $(d)$ }  if $\alpha,\beta \in \Lambda$ and $\nu \in S$ then \nl
$|\{a \in M \backslash \dbcu_{\ell < m} M_\ell:\eta^*(a) = \alpha\}| =
|\{a \in M \backslash \dbcu_{\ell < m} M_\ell:\eta^*(a) = \beta\}|$.
\endroster}
\ermn
2) Now $f^9_\Lambda(m,c)$ is defined similarly without clause (d). \nl
3) $HJ_\Lambda(m,c)$ is defined similarly omitting (d), and replacing
(b), (c) by:
\mr
\item "{$(b)'$}"  $M_\ell \ne \emptyset$
\sn
\item "{$(c)'$}"   $d \restriction S$ is constant.
\ermn
Lastly let $HJ(n,c) = HJ_{[n]}(1,c)$. \nl
4) Let $f^8_\Lambda(m,c)$ be defined as in part (2), replacing clause (b)
by (b)$'$ from part (3). \nl
5) We define $f^{10,*}_\Lambda(m,c)$ as in part (1) replacing clause (c)
by
\mr
\item "{$(c)^+$}"  for $\nu_1,\nu_2 \in S = \{\eta:\eta \in {}^{[k]} \Lambda,
\eta^* \subseteq \eta$ and each $\eta \restriction M_\ell$ is constant$\}$
we have $d(\nu_1) = d(\nu_2)$ provided that for every $\alpha \in \Lambda$
we have $|\{\ell < m:\nu_1 \restriction M_\ell$ is constantly $\alpha\}| =
|\{\ell < m:\nu_2 \restriction M_\ell$ is constantly $\alpha\}|$. 
\ermn
6) We define $f^{9,*}_\Lambda(m,c)$ as we have defined $f^{10,*}_\Lambda
(m,c)$ omitting clause (d). \nl
7) We define $f^{8,*}_\Lambda(m,c)$ as we have defined $f^{10,*}_\Lambda
(m,c)$ omitting clause (d) and replacing clause (b) by clause (b)$'$ from
part (3).
\enddefinition
\bigskip

\remark{Remark}  So 
$HJ_\Lambda(m,c)$ is the Hales-Jewett number for alphabet $\Lambda$,
getting $m$-dimensional subspace.
\endremark
\bn
\ub{\stag{app4.2} Fact}:  1) $f^{10}_\Lambda(m,c) \ge f^9_\Lambda(m,c) \ge
f^8_\Lambda(m,c)$. \nl
2) $f^{10,*}_\Lambda(m,c) \ge f^{9,*}_\Lambda(m,c) \ge
f^{8,*}_\Lambda(m,c)$. \nl
3) $f^{10}_\Lambda(m,c) \le f^{10,*}_\Lambda(m,c)$ and $f^{9}_\Lambda(m,c)
\le f^{9,*}_\Lambda(m,c)$ and $f^8_\Lambda(m,c) \le f^{8,*}_\Lambda(m,c)$.
\nl
4) $f^{8,*}_\Lambda(m,c) \le HJ_\Lambda(m,c)$.
\bigskip

\demo{Proof}  Read the definitions.

We can deal similarly with the density (like Szemeredi theorem) version
of those functions.
\enddemo
\bn
May those numbers be helpful for HJ-number?  First complimentarily to
\scite{app4.2}, clearly
\proclaim{\stag{app4.2a} Claim}  $HJ_\Lambda(m,c) \le f^{8,*}_\Lambda(m^*,c)$
if $m^*$ satisfies:
\mr
\item "{$\boxtimes$}"  assume that Par $= \{\bar \ell:\bar \ell = \langle
\ell_\alpha:\alpha \in \Lambda \rangle,\ell_\alpha \in [0,m^*)$ and
$\Sigma\{\ell_\alpha:\alpha \in \Lambda\} = m^*\}$ and $d$ is a $c$-colouring 
of Par; \ub{then} we can find $\bar \ell^\alpha \in$ Par for $\alpha \in 
\Lambda$ such that $d \restriction \{\bar \ell^\alpha:\alpha \in \Lambda\}$ is
constant and for some $\ell^* >0$ and $\langle \ell^*_\alpha:\alpha \in
\Lambda \rangle$ we have for any distinct $\alpha,\beta \in \Lambda:
\ell^\alpha_\beta = \ell^*_\beta,\ell^\alpha_\alpha = \ell^*_\alpha +
\ell^*$.
\endroster
\endproclaim
\bigskip

\remark{Remark}  We can choose $\alpha^* \in \Lambda$ let $\Lambda^* = 
\Lambda \backslash \{\alpha^*\}$ and restrict ourselves to Par$' = \{\bar \ell
\in \text{ Par}:\alpha \in \Lambda^* \Rightarrow \ell_\alpha \le m^{**} =:
m^*/|\Lambda|\}$ and let Par$'' = \{\bar \ell \restriction \Lambda^*:
\bar \ell \in \text{ Par}'\}$, now Par$'' = {}^\Lambda[0,m^{**})$, and 
$\bar \ell \mapsto \bar \ell \restriction \Lambda^*$ is a one-to-one map from
Par$'$ onto Par$''$.  So clearly it suffices to 
find a $d$-monocromatic $\{\bar \ell^*\} \cup \{\bar \ell^\alpha:
\alpha \in \Lambda^*\} \subseteq \text{ Par}''$ 
and $m > 0$ such that $\ell^\alpha_\beta = \ell^*
_\beta$ if $\beta \ne \alpha \in \Lambda^*,
\ell^\alpha_\beta = \ell^*_\alpha + m$ if
$\beta = \alpha \in \Lambda^*$.  Now this holds by $\boxtimes$ 
which is a case of the $|\Lambda^*|$-dimensional of v.d.W.
\hfill$\square_{\scite{app4.2a}}$
\endremark
\bigskip

\proclaim{\stag{app4.3} Claim}  1) $f^{10}_\Lambda(m,c) \le m \times 
HJ(|\Lambda|^m,c)$ so $f^{10}_\Lambda$ is not far from the Hales 
Jewett numbers. \nl
2) $f^{9,*}_\Gamma(m,c) \le m \times HJ(|\Gamma|^m,c)$.
\endproclaim
\bigskip

\demo{Proof}  1) Let $M_k$ be $\{0,\dotsc,k-1\}$.

Let $\Lambda_1$ be the set of function $\pi$ from $\{0,\dotsc,m-1\}$ to
$\Lambda$ such that $\alpha \in \Lambda \Rightarrow |\pi^{-1}\{\alpha\}| =
m/|\Lambda|$.

Let $n_1 = |\Lambda_1|$ so $n_1 \le |\Lambda|^m$ and $k_1 = 
HJ(|\Lambda|^m,c)$ and $k = m \times k_1$.

Let $d$ be a $c$-colouring of $V = {}^\Lambda(M_k)$.  Let
$V_1 = {}^{\Lambda_1}(M_{k_1})$ and we define a function $g$ from
$V_1$ onto $V$ as follows:

for $\eta \in V_1$, we have to define $\langle g(\eta)(a):a \in M_k \rangle,
g(\eta)(a) \in \Lambda$, now for $a \in \{0,\dotsc,k-1\}$ noting
$m[a/m] \le a < m[a/m] + m$ we define 
$g(\eta)(a) = (\eta([a/m]))([a/m] - m[a/m])$.
\sn
We define a $c$-colouring $d_1$ of $V_1$: $d_1(\eta) = d(g(\eta))$.  So
there is nonempty $N \subseteq M$ and $\rho^*_1$ a function from $M
\backslash N$ into $\Lambda_1$ such that $d_1 \restriction \{\rho \in
{}^{\Lambda_1}(M_{k_1}):\rho^*_1 \subseteq \rho$ and $\rho
\restriction N$ is constant$\}$ is constant.  
Let for $\ell < m,N_\ell = \{a:[a/m] \in N$
and $[a/m]-m[a/m] = \ell\}$ and $\rho^* \in
{}^\Lambda(M_k \backslash \dbcu_{\ell < m} N_\ell)$ be such that
$\rho^* \subseteq \rho \in {}^\Lambda(M_k) \Rightarrow \rho^*_1
\subseteq g(\rho)$.  Now check.  \nl
2) Similar proof.  ${{}}$  \hfill$\square_{\scite{app4.3}}$
\enddemo
\bn
\ub{\stag{app4.3A} Question}:  1) Can we give better bounds to
$f^\ell_\Lambda(m,c)$ than through $HJ$ for $\ell = 8,9,10$?  \nl
2) What is the order of magnitude of $f^8,f^9,f^{10}$? \nl
3) What about $f^{11}$ (see \scite{app4.10} below) and $f^{10,*},f^{9,*},
f^{8,*}$?
\bigskip

\definition{\stag{app4.4} Definition}  1) For a set $A$ let
\mr
\item "{$(a)$}"  seq$_\ell(A) = \{\eta:\eta$ is a sequence of length $\ell$
with no repetitions, of \nl

$\qquad \qquad$ elements of $A\}$
\sn
\item "{$(b)$}"  seq$(A) = \dbcu_{\ell \ge 1} \text{ seq}_\ell(A)$
\sn
\item "{$(c)$}"  seq$_{m,\ell}(A) = \{\bar \eta:\bar \eta = \langle \eta_i:i
< m \rangle,\eta_i \in \text{ seq}_\ell(A)$ and $i_1 < i_2 \Rightarrow
\text{ Rang}(\eta_{i_1}) \cap \text{ Rang}(\eta_{i_2}) = \emptyset$ 
\sn
\item "{$(d)$}"  seq$^*_m(A) = \cup\{\text{seq}_{m,\ell}(A):\ell \ge 1\}$
\sn
\item "{$(e)$}"  seq$^*(A) = \cup\{\text{seq}^*_m(A):m \ge 1\}$.
\ermn
2) For $\bar \eta \in \text{ seq}^*(A)$ let \nl
son$(\bar \eta) = \{\nu:\nu \in
\text{ seq}(A)$ and $\nu$ is a concatenation of some members of \nl
$\{\eta_\ell:\ell < \ell g(\bar \eta)\}$, in any order$\}$ \nl
legson$(\bar \eta) = \{\nu \in \text{ son}
(\bar \eta):\text{ Rang}(\nu) = \dbcu_i \text{ Rang}(\eta_i)\}$, \nl
dis$(\bar \eta) = \{\bar \nu \in \text{ seq}^*(A):\text{each } \nu_i$ is from
son$(\bar \eta)\}$, \nl
leg dis$(\bar \eta) = \{\bar \nu \in \text{ dis}(\bar \eta):
\dbcu_i \text{ Rang}(\nu_i) = \dbcu_i \text{ Rang}(\eta_i)\}$.
\nl
3) $f^{12}(m,c)$ is the first $k$ such that $k = \omega$ or $k < \omega$ and
for every $c$-colouring $d$ of seq$([0,k))$ there is $\bar \eta \in
\text{ seq}^*_m(A)$ such that the set son$(\bar \eta)$ is $d$-monocromatic.
\nl
There are other variants.
\enddefinition
\bn
\ub{\stag{app4.6} Question}:   Is $f^{12}(m,c)$ finite? 
\bigskip

\definition{\stag{app4.5} Definition}  1) For groups $G,H$ and subset $Y$ of
$H$ and cardinal $\kappa$ let $G \rightarrow (Y,H)_\kappa$ means that for any
$\kappa$-colouring $d$ of $G$ (i.e. $d$ is a function from $G$ into a set of
cardinality $\le \kappa$) there is an embedding $h$ of $H$ into $G$ such that
$d \restriction \{h(y):y \in Y\}$ is constant. \nl
2) $G \rightarrow (Y,H)_{\kappa,\theta}$ is defined similarly but
$d \restriction \{h(y):y \in Y\}$ has range with $< \theta$ members. \nl
3) If $Y=H$ we may omit it.
\enddefinition
\bn
\ub{\stag{app4.6a} Question}:  1) Investigate $G \rightarrow (Y,H)_c$ for
finite groups. \nl
2) Assume $H$ is a finite permutation group, $Y$ is one conjugacy class (say
permutation of order two) and $c$ finite, does $G \rightarrow (Y,H)_c$
exist?  (This is connected to \scite{app4.4}, just interpret $\eta \in
\text{ seq}(A)$ of even legnth $2n$ with the permutation of $A$ permuting
$\eta(i)$ with $\eta(n+1)$ for $i < n$. \nl
3) Similarly when we colour subgroups of $G$.
\bn
Similar problems to \scite{app4.3A} are
\definition{\stag{app4.10} Definition}  1) (See \cite{Sh:679} and the
notation there).
Let $f^1_{\bar \Lambda}(m,c)$ is the first $k$ such that $k=\omega$
(i.e. infinity) or
\mr
\item "{$(*)_k$}"  $k$ is divisible by $|\Lambda_{\text{id}}|$ and
letting $M = M^\tau_k$, we have: for every $c$-colouring $d$ of 
Space$_{\bar \Lambda}(M)$, we can find an $m$-dimensional subspace 
$S$ such that
{\roster
\itemitem{ $(a)$ }  letting $\langle M_\ell:\ell < m \rangle$ be a witness
for $S$ (see Definition \cite[1.7t(5)]{Sh:679} the dimension of 
$M_\ell,|P^{M_\ell}|$, is the same for $\ell = 0,\dotsc,m-1$
\sn
\itemitem{ $(b)$ }  let $K = K^\tau_{[0,m)},N = c \ell_M(\dbcu_{\ell < m}
M_\ell)$ and $f$ as in (c) of Definition \cite[1.7t(5)]{Sh:679}; if 
$\rho_1,\rho_2
\in$ Space$_{\bar \Lambda}(K)$ and $\nu_1,\nu_2 \in S$ are such that $b \in
N \Rightarrow \nu_1(b) = \rho_1(\hat f(b)) \and \nu_2(b) = \rho_2(\hat f
(b))$ and there is an automorphism of $K$ mapping $\rho_1$ to $\rho_2$ then
$d(\nu_1) = d(\nu_2)$.
\endroster}
\ermn
2)  If $\tau = \{\text{id}\},\Lambda = \Lambda_{\text{id}}$ we 
write $f^{11}_\Lambda(m,c)$ and above (b) means:
\mr
\item "{$(b)'$}"  for $\nu_1,\nu_2 \in S$ we have $d(\nu_1) = d(\nu_2)$ if
$\nu_1 \restriction (M \backslash \dbcu_{\ell < m} M_\ell) = \nu_2
\restriction (M \backslash \dbcu_{\ell < m} M_\ell)$ and for every $\alpha
\in \Gamma$ we have
$$
|\{\ell < m:\nu_1 \restriction M_\ell \text{ is constantly } \alpha\}| =
|\{\ell < m:\nu_2 \restriction M_\ell \text{ is constantly } \alpha\}|.
$$
\endroster
\enddefinition
\shlhetal

%% you may want to move the following lines up a bit
\newpage
    
REFERENCES.  
\bibliographystyle{lit-plain}
\bibliography{lista,listb,listx,listf,liste}

\enddocument%%

\bye